\documentclass[11pt]{amsart}
\usepackage{amsfonts}
\usepackage{amsmath}
\usepackage{fancyhdr,amssymb}
\usepackage{amsthm}
\usepackage{ color, ulem}
\usepackage{tkz-graph}
\input xy
\xyoption {all}

\def \CAT {{\mathsf{C}}}

\newcommand{\comment}[1]{}
\newtheorem{theorem}{Theorem}
\newtheorem {lemma}[theorem]{Lemma}
\newtheorem{conjecture}[theorem]{Conjecture}
\newtheorem{question}[theorem]{Question}
\newtheorem{definition}[theorem]{Definition}
\newtheorem {corollary}[theorem]{Corollary}
\newtheorem {proposition}[theorem]{Proposition}
\theoremstyle{definition}
\newtheorem {example}[theorem]{Example} 
\theoremstyle {definition}
\newtheorem{remark}[theorem]{Remark}

\begin{document}
\baselineskip=16pt
\title [Quot schemes of curves and surfaces]
{Quot schemes of curves and surfaces:  \\ virtual classes, integrals,
Euler characteristics}

\author{D. Oprea}
\address{Department of Mathematics, University of California, San Diego}
\email {doprea@math.ucsd.edu}
\author{R. Pandharipande}
\address{Department of Mathematics, ETH Z\"urich}
\email {rahul@math.ethz.ch}

\begin{abstract}
We compute tautological integrals over Quot schemes 
 on curves and surfaces. After obtaining several explicit formulas 
over Quot schemes of dimension
0 quotients on curves (and finding a new symmetry), 
we apply the results to 
tautological integrals against the virtual fundamental classes
of Quot schemes of dimension 0 and 1 quotients 
on surfaces (using also universality, torus localization, and cosection
localization). 
The virtual Euler characteristics
of Quot schemes of  surfaces, a new 
theory  parallel to the Vafa-Witten Euler characteristics of
the moduli of bundles, is defined and studied. Complete formulas
for the virtual Euler characteristics are found in the case of dimension
0 quotients on surfaces. Dimension 1 quotients are studied on
 $K3$ surfaces and surfaces
of general type with connections to the Kawai-Yoshioka formula and the
Seiberg-Witten invariants respectively. The dimension 1 
theory is completely solved
for minimal surfaces of general type admitting a nonsingular canonical curve.
Along the way, we find a new connection
between weighted tree counting and multivariate Fuss-Catalan numbers which is of independent
interest. 
\end{abstract}

\maketitle

\setcounter{tocdepth}{1}
\tableofcontents

\section{Introduction}
\subsection{Overview}
The main goal of the paper is to study the virtual 
fundamental classes
of Quot schemes  of surfaces. The parallel
study for 3-folds was undertaken in \cite{MNOP1, MNOP2} and
led to the MacMahon function for Hilbert schemes of points 
and the  GW/DT correspondence
 for Hilbert schemes of curves.
For the surface case, we use
several techniques: the universality results of \cite{EGL}, $\mathbb{C}^\star$-equivariant localization of the virtual class \cite{GP}, and cosection localization \cite{KL}. 
However, the most important input to the surface theory concerns
the parallel study of Quot schemes of curves 
of quotients with dimension 0 support, 
which we develop first. 
By applying the curve results to the surface theory, we
prove several basic results about the virtual fundamental
classes of Quot schemes of quotients with supports of dimension 0 and 1 
on surfaces. The subject is full of open questions.

\subsection{Curves} \label{cccvvv}
Let $C$ be a nonsingular projective curve. 
Let $\mathsf {Quot}_{\,C}(\mathbb C^N, n)$ parameterize short exact sequences $$0\to S\to \mathbb C^N\otimes \mathcal O_C\to Q\to 0\, ,$$ where $Q$ is a rank 0
sheaf on $C$ with $$\chi(Q)=n\, .$$ 
The scheme ${\mathsf {Quot}_{\,C} (\mathbb C^N, n)}$ was viewed in \cite{quotients} as the stable quotient compactification of {\it degree $n$ maps 
to the point}, where the point  is the degenerate Grassmannian ${\mathbb G}(N, N)$.
By analyzing the Zariski tangent space,
 $\mathsf {Quot}_{\,C} (\mathbb C^N, n)$ is easily seen
to be a nonsingular projective variety of dimension $Nn$, see 
\cite[Section 4.7]{quotients}.

For a vector bundle $V\to C$ of rank $r$, the assignment $$Q\mapsto H^0(C, V\otimes Q)$$ 
for $[\mathbb C^N\otimes \mathcal O_C\to Q] \in \mathsf{Quot}_{\,C}(\mathbb C^N,n)$ 
defines a  {tautological} vector bundle  $$V^{[n]}\to \mathsf {Quot}_{\,C} (\mathbb C^N, n)$$
of rank $rn$.
 The construction descends to $K$-theory via locally free resolutions. 
We define  generating series of Segre{\footnote{For a vector bundle $V$ on a scheme $X$, we write $$s_t(V)=1+ts_1(V)+t^2 s_2(V)+\ldots $$ for the total Segre class.}} classes
on Quot schemes of curves
as follows.

\begin{definition}\label{ssseg}
Let
$\alpha_1, \ldots, \alpha_\ell$ be $K$-theory classes on $C$. Let
\begin{equation*}
\mathsf Z_{C, N}(q, x_1, \ldots, x_\ell\,|\,\alpha_1, \ldots, \alpha_\ell)=\sum_{n=0}^{\infty} q^n \int_{\mathsf {Quot}_{\,C} (\mathbb C^N, n)} s_{x_1}(\alpha_1^{[n]})\cdots
 s_{x_\ell} (\alpha_\ell^{[n]})\, .
\end{equation*}

\end{definition}

Since the integrals in Definition \ref{ssseg} depend upon $C$
only through the genus $g$ of the curve, we will often write
$$\mathsf Z_{g, N}(q, x_1, \ldots, x_\ell\,|\,\alpha_1, \ldots, \alpha_\ell)=
\mathsf Z_{C, N}(q, x_1, \ldots, x_\ell\,|\,\alpha_1, \ldots, \alpha_\ell)\, .$$
By the arguments of \cite {EGL}, there exists a
factorization \begin{equation}\label{factor}\mathsf Z_{g, N}(q, x_1, \ldots, x_\ell\,|\,\alpha_1, \ldots, \alpha_\ell)= {\mathsf A}_1^{c_1(\alpha_1)} \cdots {\mathsf A}_\ell^{c_1(\alpha_\ell)} \cdot \mathsf B^{1-g},\end{equation} for universal series 
\begin{equation}
\mathsf A_1\, ,\, \ldots\, ,\, \mathsf A_\ell\, ,\, \mathsf B\, \in \mathbb Q[[q, x_1, \ldots, x_\ell]] \label{tter}
\end{equation}
which {\it do not} depend on the genus $g$
 or the degrees $c_1(\alpha_i)$.
However, the series \eqref{tter} {\it do}
depend on the ranks 
$$\mathsf{r}=(r_1,\ldots, r_\ell)\, , \ \ \ r_i=\text{rank } \alpha_i\,
$$
and $N$.
The complete 
notation for the series \eqref{tter} is 
\begin{equation}
\mathsf A_{1,\mathsf{r},N}\, ,\,  \ldots\, ,\,  \mathsf A_{\ell, \mathsf{r},N}\, ,\,  \mathsf B_{\mathsf{r},N}\, 
\in \mathbb Q[[q, x_1, \ldots, x_\ell]] \, ,
\end{equation}
but we will often use the abbreviated notation \eqref{tter} with
the ranks $r_i$ and $N$ suppressed.

\begin{question} \label{q1}
Find closed-form expressions for the series 
$\mathsf A_{i,\mathsf{r},N}$ and $\mathsf B_{\mathsf{r},N}$.
\end{question} 

Integrals over Quot schemes of curves were also studied in \cite {mo}
via equivariant localization. In particular,   
formulas of Vafa-Intriligator \cite{Ber, I, ST} were recovered and extended.

 \subsection{Symmetric products ($N=1$)} 
For curves, the symmetric product $C^{[n]}$ is the Quot scheme 
in the $N=1$ case,
$$C^{[n]} =\mathsf {Quot}_{\,C}(\mathbb C^1, n)\, .$$ 
We give a complete answer 
to Question \ref{q1} for $N=1$. The result  will later play an important
role in our study of Quot schemes of surfaces.

\begin{theorem} \label{thm1} Let $\alpha_1, \ldots, \alpha_\ell$ have ranks $r_1, \ldots, r_{\ell}$, and let $N=1$.
 Then $$\mathsf Z_{g,1}(q, x_1, \ldots, x_\ell\,|\,\alpha_1, \ldots, \alpha_\ell)={\mathsf A}_1(q)^{\,c_1(\alpha_1)} \cdots {\mathsf A}_\ell(q)^{\,c_1(\alpha_\ell)}\cdot {\mathsf B}(q)^{1-g}\, ,$$ where, for the change of variables 
\begin{equation}\label{cchh}
q=t(1-x_1t)^{r_1}\cdots (1-x_\ell t)^{r_\ell}\, ,
\end{equation}
 we set $$\mathsf A_i(q)=1-x_i\cdot t\, ,\,\,\,\, \ \ \mathsf B(q)= \left(\frac{q}{t}\right)^2\cdot \frac{dt}{dq}\, .$$
\end{theorem} 
\vspace{5pt}

To compute the series{\footnote{For Theorem \ref{thm1}, the complete
notation is 
$\mathsf A_i = \mathsf{A}_{i,{\mathsf r}, 1}$ and
$\mathsf B= \mathsf B_{{\mathsf r},1}$.}}
 $\mathsf A_i(q)$ and $\mathsf{B}(q)$,
the change of variables \eqref{cchh}
 must be inverted to write $t$ as a function of $q$  with
$x_1,\ldots,x_\ell$ viewed as parameters.
By Theorem \ref{thm1}, the series $\mathsf Z_{g,1}(q, x_1, \ldots, x_\ell\,|\,\alpha_1, \ldots, \alpha_\ell)$
is a function in $q$ which is {\it algebraic} over the field $\mathbb Q(x_1, \ldots, x_\ell).$

\begin{remark} Specializing to the case $\ell=1$, $x_1=1$, and $r_1=r$, and letting $V\to C$ be a rank $r$ vector bundle, we recover the result of \cite {MOP}: 
\begin{equation}\label{dxxd}
\sum_{n=0}^{\infty} q^n \int_{C^{[n]}} s_{n} (V^{[n]}) =
\exp \Big(c_1(V)\cdot \widehat{\mathsf A}(q)+(1-g)\cdot 
\widehat{\mathsf B}(q)\Big)
\end{equation} for the series
\begin{eqnarray}\label{sese}\widehat{{\mathsf A}}(t(1-t)^r)&=&\log(1-t)\, , \\ \nonumber
  \widehat{{\mathsf B}}(t(1-t)^r)&=&(r+1)\log(1-t)-\log(1-t(r+1))\, .
\end{eqnarray}
These expressions  confirmed and expanded predictions of \cite {Wang}. The
$r=1$ case is related to the counts of secants to projectively embedded curves \cite {C,LeB}. 
\end{remark}
\begin{remark} To go beyond numerical invariants, we consider  a flat family
$$\pi:C\to S\, $$
of nonsingular projective curves with line bundles
$L_1,\ldots, L_\ell\to C$.
We write $$\pi^{[n]}: C^{[n]}\to S$$ for the relative symmetric product. A more difficult question concerns the calculation of the push-forwards $$\sum_{n=0}^{\infty} q^n \pi^{[n]}_{\star} \left(s_{x_1} (L_1^{[n]})\cdots s_{x_\ell}(L_{\ell}^{[n]})\right)\in A^{\star}(S)$$ in terms of the classes $$\kappa[a_1, \ldots, a_\ell, b]=\pi_{\star} \left(c_1(L_1)^{a_1}\cdots c_1(L_\ell)^{a_{\ell}} \cdot c_1(\omega_{\pi})^b\right)\in A^{\star}(S)\, .$$
When $\pi$ is the universal family over the moduli space of curves,
such constructions play a role in the study of tautological classes
\cite{P,PP}. \end{remark}
\vskip.1in

\subsection {Higher $N$ (for $\ell=1$)} \label{hnl1}
Our second result concerns the case of arbitrary $N$, but we assume  $\ell=1$. The corresponding series is
$$\mathsf Z_{g, N}(q\,|\, V)=\sum_{n=0}^{\infty} q^n \int_{\mathsf {Quot}_{\,C} (\mathbb C^N, n)} s(V^{[n]})\, ,$$ where $V\to C$ is a rank $r$ vector bundle. 
\begin{theorem} \label{thm2}
The universal Segre series is $$\mathsf Z_{g, N}(q\,|\, V)={\mathsf A}(q)^{c_1(V)}\cdot {\mathsf B}(q)^{1-g}\, ,$$ where $$\log \mathsf A(q)=\sum_{n=1}^{\infty} (-1)^{(N+1)n+1} \binom{(r+N)n-1}{Nn-1}\cdot \frac{q^n}{n}\, .$$ 
\end{theorem}
 
\begin{remark} In case $N=1$, Theorem \ref{thm2} is a special case of Theorem \ref{thm1}. The agreement of the formulas follows from the identity
 $$-\log(1-t)= \sum_{n=1}^{\infty} \binom{(r+1)n-1}{n-1}\cdot \frac{q^n}{n}
\ \ \ \ \ \text{for}\ \ \ \ \ q=t(1-t)^r\, $$ 
which will be proven in Lemma \ref{identity} below.  
\end{remark}

Theorem \ref{thm2} identifies the $\ell=1$ series
$\mathsf A=\mathsf{A}_{1,r,N}$, but does not specify the series $\mathsf B=\mathsf B_{r, N}$. 
However, for rank $r=1$, {\it closed-form} expressions for 
the $\mathsf A$ and $\mathsf B$-series are determined by the following
result. 
\begin{theorem}\label{cor1} For
 $\text{rank }V=1$, after the change of variables $$q=(-1)^N t(1+t)^{N}$$ we have $$\mathsf A_{1,1,N}(q)=(1+t)^N\,\,\, \text{ and }\,\,\, \mathsf B_{1,N}(q)=\frac{(1+t)^{N+1}}{1+t(N+1)}\, .$$ 
\end{theorem} 
\noindent We also write an explicit power series expansion for the $\mathsf B$-series parallel to Theorem \ref{thm2}.
\begin{corollary} \label{c2} For $\text{rank }V=1$, we have 
$$\mathsf B_{1,N}(q)= \sum_{n=0}^{\infty} (-1)^{n(N+1)} \cdot \binom{(n-1)(N+1)}{n} \cdot q^n.$$
\end{corollary} 

\noindent By comparing the expressions of Theorem \ref{cor1} with those of equation \eqref{sese},
 we obtain the following new symmetry exchanging $N$ and the rank.
\begin{corollary} \label{cor2} 
For any line bundle $L\to C$, we have 
$$\int_{{\mathsf{Quot}}_{\,C} (\mathbb C^N, n)} s(L^{[n]})=(-1)^{n(N-1)} \int_{C^{[n]}} s(L^{[n]})^{N}\, .$$ In particular, for $C=\mathbb P^1$,
 we have $$\int_{{\mathsf{Quot}}_{\,\mathbb P^1} (\mathbb C^N, n)} s(L^{[n]})=(-1)^{Nn} \binom{N\deg L-N(n-1)}{n}\, .$$
\end{corollary} 
\vspace{3pt}

\subsection {Catalan numbers} By specializing Theorem \ref{thm1} to the case of an elliptic curve $C$ and using Wick expansion techniques, we are led 
to a combinatorial identity for Catalan numbers which appears to be 
new.{\footnote{There are many
realizations of the Catalan numbers! But we have asked several experts
and ours does not appear to be in the literature. If you
know a reference, please tell us.}} 

The $m^{th}$ Catalan number $$\CAT_m=\frac{1}{m+1}\binom{2m}{m}$$ 
is well-known to count unlabelled ordered trees with $m+1$ vertices \cite {S}. The {\it multivariate Fuss-Catalan} numbers were introduced and studied in \cite {A}. A special case of the definition is used here. For non-negative integers $p_1, \ldots, p_k$, the multivariate Fuss-Catalan number of interest to us is 
 $$\CAT(p_1, \ldots, p_k)=\frac{1}{p_1+\ldots+p_k+1}\binom{2p_1+p_2+\ldots+p_k}{p_1}\cdots 
 \binom{p_1+p_2+\ldots+2p_k}{p_k}\, .$$ The case $k=1$ corresponds to the usual Catalan number $\CAT(m)=\CAT_m$. The multivariate Fuss-Catalan numbers were shown to count certain 
$k$-Dyck paths or, alternatively, $k$-nary trees,
 and also arise in connection with algebras of $B$-quasisymmetric
polynomials \cite {A}. 

We interpret the Catalan
and multivariate Fuss-Catalan numbers as a {\it weighted} count of trees. Let non-negative integers
$p_1,\ldots,p_k$ be given. Let 
$$n= p_1+\ldots+p_k+1\, .$$
A {\it labelled $k$-colored tree of type $(p_1,\ldots, p_k)$} is
a tree $T$ with
\begin{enumerate}
\item[$\bullet$]
$n$ vertices labelled $\{1, 2, \ldots, n\}$, 
\item[$\bullet$] $n-1$ edges each painted with one of the $k$ different colors
such that exactly $p_j$ edges are painted with the $j^{th}$ color.
\end{enumerate}
For each vertex $v$, we write $$d_v^1, \ldots, d_v^k$$ for the 
{\it out-degrees}{\footnote{The term {\it out-degree} comes from
regarding $T$ as an oriented graph with each edge oriented
in the direction of {\it decreasing} vertex label.}}
 of $v$ corresponding to each of the $k$ colors. More precisely, 
$d_v^j$ counts edges $e$ incident to $v$, of color $j$, such that $e$
connects $v$ to a vertex $w$ 
 satisfying $$v>w\, .$$ 
We define the {\it weight} of $T$ as the product $$\mathsf{wt} (T)=\frac{1}
{(n-1)!} \prod_{v\, \text{vertex}} d_v^1! \cdots  d_v^k!\, .$$ 
\begin{theorem} \label{thm3}
The Fuss-Catalan number is the 
weighted count of ordered $k$-colored trees of type
$(p_1,\ldots,p_k)$:
$$\CAT(p_1,\ldots,p_k)= \sum_{T} \mathsf{wt}(T)\, .$$   
\end{theorem} 

\begin{example} Let us now specialize to the single color ($k=1$) case with 
$m=p_1$ and $n=m+1$.
The weights then take the form:
$$\mathsf{wt} (T)=\frac{1}{m!} \prod_{v\, \text{vertex}} d_v!=\binom{m}{d_1, \ldots, d_n}^{-1}\, ,$$ 
where for each $v$,  $d_v$ denotes the out-degree of $v$. 
We then obtain the standard $m^{th}$ Catalan number as a weighted count
of labelled trees with $m+1$ vertices:
\begin{equation}\label{ll223}
\CAT(m)=\sum_{T} \mathsf{wt}(T) \, .
\end{equation}
The result \eqref{ll223} should perhaps be compared with the realization of
$\CAT(m)$ as the 
unweighted count of unlabelled ordered trees with $m+1$ vertices (see \cite{S} for instance).
The following diagram shows the two counts for $\CAT(2)$:\vskip.1in
\begin{itemize}
\item [--] weighted count \vskip.2in
\begin{tikzpicture}[shorten >=1pt,->]
  \tikzstyle{vertex}=[circle,fill=black!25,minimum size=12pt,inner sep=2pt]
  \node[vertex] (G_1) at (1,0) {1};
  \node[vertex] (G_2) at (2,0)   {2};
  \node[vertex] (G_3) at (3,0)  {3};
  \draw (G_1) -- (G_2) -- (G_3) -- cycle;

  \tikzstyle{vertex}=[circle,fill=black!25,minimum size=12pt,inner sep=2pt]
  \node[vertex] (H_1) at (6,0) {3};
  \node[vertex] (H_2) at (7,0)   {1};
  \node[vertex] (H_3) at (8,0)  {2};
  \draw (H_1) -- (H_2) -- (H_3) -- cycle;
  \tikzstyle{vertex}=[circle,fill=black!25,minimum size=12pt,inner sep=2pt]
  \node[vertex] (K_1) at (11,0) {1};
  \node[vertex] (K_2) at (12,0)   {3};
  \node[vertex] (K_3) at (13,0)  {2};
  \draw (K_1) -- (K_2) -- (K_3) -- cycle;
 \text{weight 1}
\end{tikzpicture}\vskip.1in
\item [--] unweighted count \vskip.2in
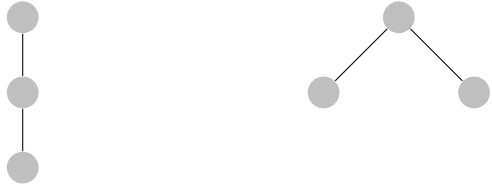
\begin{center}\begin{tikzpicture}[shorten >=1pt,->]
  \tikzstyle{vertex}=[circle,fill=black!25,minimum size=12pt,inner sep=2pt]
  \node[vertex] (G_1) at (1,0) {};
  \node[vertex] (G_2) at (1,1)  {};
  \node[vertex] (G_3) at (1,2) {};
  \draw (G_1) -- (G_2) -- (G_3) -- cycle;
 \tikzstyle{vertex}=[circle,fill=black!25,minimum size=12pt,inner sep=2pt]
  \node[vertex] (H_1) at (5,1) {};
  \node[vertex] (H_2) at (6,2)   {};
  \node[vertex] (H_3) at (7,1)  {};
  \draw (H_1) -- (H_2) -- (H_3) -- cycle;

\end{tikzpicture}
\end{center}
\end{itemize}
In the first count, the weights are $\frac{1}{2}, \frac{1}{2}$ and $1$ respectively and $$\CAT(2)=\frac{1}{2}+\frac{1}{2}+1.$$
\end{example}
\subsection{Surfaces: dimension 0 quotients} \label{sd0}
We can apply the above results for curves to the
calculation of tautological integrals over
Quot schemes of dimension 0 quotients 
of nonsingular projective surfaces $X$.

The Quot scheme $\mathsf {Quot}_{\,X}(\mathbb C^N, n)$ of 
short exact sequences $$0\to S\to \mathbb C^N\otimes \mathcal O_X\to Q\to 0,\,\,\,\ \  \chi(Q)=n\, ,\  c_1(Q)=0\, ,\ \text{rank}(Q)=0$$  
is known \cite {EL,Li}
 to be irreducible of dimension $n(N+1)$, but 
 may be singular.{\footnote{An example is given in Section \ref{vis} below.}}
 Since the higher obstructions for the
standard deformation theory lie in \begin{equation}\label{hio}\text{Ext}^2(S, Q)=\text{Ext}^0(Q, S\otimes K_X)^\vee=0,\end{equation} the Quot scheme carries a 
2-term perfect obstruction theory and a virtual fundamental cycle of dimension $$\text{Ext}^0(S, Q)-\text{Ext}^1(S, Q)=\chi(S, Q)=Nn\, .$$

\begin{question} \label{q2}Evaluate the integrals $$\mathsf Z_{X, N}(q, x_1, \ldots, x_\ell\,|\,\alpha_1, \ldots, \alpha_\ell)=\sum_{n=0}^{\infty} q^n 
\int_{\left[\mathsf {Quot}_{\,X}(\mathbb C^N, n)\right]^{\mathrm{vir}}} s_{x_1} (\alpha_1^{[n]}) \cdots s_{x_\ell} (\alpha_\ell^{[n]})\, $$
where $\alpha_1,\ldots,\alpha_\ell$ are K-theory classes on $X$. 
\end{question} 

By our next result, the surface series of Question \ref{q2}
 are obtained from the parallel curves series of Question \ref{q1}.
The relationship is not unlike the localization result for
the Gromov-Witten theory of surfaces of general type with respect to
a canonical  divisor \cite{KL,LP,MP}.

\begin{theorem} \label{thm4} 
Let the ranks of the
classes $\alpha_1,\ldots, \alpha_\ell$ be given by
$\mathsf{r}=(r_1,\ldots, r_\ell)$. Let the
series $\mathsf A_{1,\mathsf{r},N}, \ldots, 
\mathsf A_{\ell,\mathsf{r},N}, \mathsf B_{\mathsf{r},N}$ be 
defined by the curve integrals \eqref{factor}. Then, we have 
\begin{multline*}
\mathsf Z_{X, N}(q, x_1, \ldots, x_\ell\,|\,\alpha_1, \ldots, \alpha_\ell)=\\
 \mathsf A_{1,\mathsf{r},N}(-q)^{c_1(\alpha_i)\cdot K_X} \cdots
 \mathsf A_{\ell,\mathsf{r},N}(-q)^{c_1(\alpha_\ell)\cdot K_X}\cdot 
B_{\mathsf{r},N}(-q)^{-K_X^2}\, .
\end{multline*}
\end{theorem}

In case $X$ is a surface of general type with a nonsingular 
canonical divisor 
$$C\subset X\, ,$$
then $c_1(\alpha_i)\cdot K_X$ is the degree of the restriction of
$\alpha_i$ to $C$ and
$$-K_X^2 = 1-\text{genus}(C)$$
by adjunction. We may therefore write Theorem \ref{thm4} as
\begin{equation*}
\mathsf Z_{X, N}(q, x_1, \ldots, x_\ell\,|\,\alpha_1, \ldots, \alpha_\ell)=\\
\mathsf Z_{g(C), N}\left(-q, x_1, \ldots, x_\ell\,\Big|\,\alpha_1|_C, 
\ldots, \alpha_\ell|_C\right) \, .
\end{equation*}
However, Theorem \ref{thm4} holds for all $X$ (even if $X$
is not of general type).

For $N=1$, Theorems \ref{thm1} and \ref{thm4} together yield 
a complete answer for the virtual Segre integrals over the Hilbert scheme of points, $$X^{[n]} = \mathsf {Quot}_{\,X}(\mathbb C^1, n)\, .$$

\begin{corollary} 
Let $X$ be a nonsingular projective
surface. Then
 $$\sum_{n=0}^{\infty} q^n \int_{\left[X^{[n]}\right]^{\mathrm{vir}}} s_{x_1} (\alpha_1^{[n]}) \cdots s_{x_\ell} (\alpha_\ell^{[n]})= \mathsf A_1(q)^{c_1(\alpha_1)\cdot K_X}
\cdots \mathsf A_\ell(q)^{c_1(\alpha_\ell)\cdot K_X}\cdot \mathsf B(q)^{-K_X^2}$$ where,
 for the change of variable  $$q=-t(1-x_1t)^{r_1}\cdots (1-x_\ell t)^{r_\ell}\, ,$$  
we set $$\mathsf A_i(q)=1-x_i\cdot t\, ,\,\,\,\, \ \mathsf B(q)= -\left(\frac{q}{t}\right)^2\cdot \frac{dt}{dq}\, .$$
\end{corollary}

\noindent Similarly, for higher $N$, Theorems \ref{cor1} and \ref{thm4} yield the following evaluation. 
\begin{corollary} Let $L\to X$ be a line bundle on
a nonsingular projective surface. Then
$$\sum_{n=0}^{\infty} q^n \int_{\left[\mathsf {Quot}_X(\mathbb C^N, n)\right]^{\mathrm{vir}}} s(L^{[n]})=\mathsf A(q)^{c_1(L)\cdot K_X} \cdot \mathsf B(q)^{-K_X^2}$$ where,
for the change of variables
$$q=(-1)^{N+1}\,t\,(1+t)^{N}\, ,$$
we set
 $$\mathsf A(q)=(1+t)^N, \ \ \,\,\,\, \mathsf B(q)=\frac{(1+t)^{N+1}}{1+(N+1)t}\, .$$
\end{corollary}

\begin{remark} Question \ref{q2} is well-posed for integrals against
the actual fundamental class of dimension
$n(N+1)$ of $\mathsf {Quot}_{\,X}(\mathbb C^N, n)$ instead of
the virtual fundamental class of dimension $nN$.
The calculation for the actual fundamental class
 is more complicated. 
The  $N=1$ case is by far the most studied. 
Then, the series
 $$\mathsf Z_X(q, x_1, \ldots, x_\ell\,|\,\alpha_1, \ldots, \alpha_\ell)=\sum_{n=0}^{\infty} q^n \int_{X^{[n]}} s_{x_1} (\alpha_1^{[n]}) \cdots s_{x_r} (\alpha_\ell^{[n]})$$ are generalizations of the Segre integrals considered by Lehn \cite {L}. In fact, 
Lehn's case  corresponds to $\ell=1$ and $\text{rank }\alpha_1=1$, and was studied in \cite{segre, MOP, voisin}. 
The case $$x_1=\ldots=x_\ell=1$$ was studied in \cite {MOP2}, and a complete solution was given for $K$-trivial surfaces. The case $\ell=2$ was analyzed in \cite {WZ}, and the answer was found for all surfaces if
$$\text{ rank }\alpha_1=\text{rank }\alpha_2=-1$$ 
via connections to $K$-theory. 
\end{remark}

\subsection{Virtual Euler characteristics: dimension $0$ quotients} 
\label{vsd0} The topological Euler characteristics of the schemes $\mathsf {Quot}_{\,C}(\mathbb C^N, n)$ and $\mathsf {Quot}_{\,X}(\mathbb C^N, n)$ can be easily 
computed via equivariant localization: $$\sum_{n=0}^{\infty} q^n \mathsf e(\mathsf {Quot}_{\,C}(\mathbb C^N, n))=(1-q)^{N(2g-2)}\, ,$$ $$\sum_{n=0}^{\infty} q^n \mathsf e(\mathsf {Quot}_{\,X}(\mathbb C^N, n))=\prod_{n=1}^{\infty} (1-q^n)^{-N\chi(X)}\, .$$

More subtle is  the virtual Euler characteristic of $\mathsf {Quot}_{\,X}(\mathbb C^N, n)$ defined via the 2-term obstruction theory. 
A basic result for dimension $0$ quotients, proven
using a reduction to the Quot schemes of curves,
is the following rationality statement.
\begin{theorem}\label{rattt}
The generating series of virtual Euler characteristics
of $\mathsf {Quot}_{\,X}(\mathbb C^N, n)$ is a rational
function in $q$ which depends only upon $K_X^2$ and $N$,
$$\sum_{n=0}^{\infty} q^n \mathsf e^{\mathrm{vir}} (\mathsf {Quot}_{\,X}(\mathbb C^N, n))= \mathsf{U}_N^{K_X^2}\, , \ \ \ \ \mathsf{U}_N \in \mathbb{Q}(q)\,  .$$
\end{theorem}

We can calculate  $\mathsf{U}_1$ directly
using the evaluations given in Theorem \ref{thm1}:
$$
\mathsf{U}_1
=
\frac{(1-q)^2}{1-2q}\, .
$$ For higher $N$, a more involved computation in Section \ref{ves}
yields an exact expression in a different form:
\begin{equation} \label{mmsst}
 \mathsf{U}_N(q) =  \frac{(1-q)^{2N}}{(1-2^Nq)^{N}}\cdot \prod_{i<j} (1-(r_i-r_j)^2)\, ,
\end{equation}
 where $r_1(q), \ldots, r_N(q)$ are the $N$ distinct
roots of the polynomial equation 
 $$z^N-q(z-1)^N=0\, $$
in the variable $z$. The shape of the answer is reminiscent of the Vafa-Intriligator formulas for Quot schemes of curves \cite{Ber, I, mo, ST} which yield expressions depending on the roots of unity.

Using \eqref{mmsst}, we can easily calculate $\mathsf{U}_N$ as
a rational function of $q$. The next few cases are:
\begin{eqnarray}
\nonumber 
\nonumber \label{thm5}
\mathsf{U}_2
&=& \frac{(1-q)^2(1-6q+q^2)}{(1-4q)^2}\, ,\\
\nonumber 
\mathsf{U}_3
&=&\frac{(1-q)^2(1 - 22 q + 150 q^2 - 22 q^3 + q^4)}{(1-8q)^3}\, ,\\
\mathsf{U}_4 \nonumber
&=&\frac{(1-q)^2(1 - 62 q + 1407 q^2 - 15492 q^3 + 1407 q^4 - 62 q^5 + q^6)}{(1-16q)^4}\, .
\end{eqnarray}
Formula \eqref{mmsst} implies  
\begin{equation}\label{vrrt}
\mathsf{U}_N(q) = \frac{(1-q)^2}{(1-2^Nq)^N} \cdot \mathsf{P}_N(q)\, ,
\end{equation}
where $\mathsf{P}_N(q)\in \mathbb{Z}[q]$ is a palindromic polynomial of 
degree $2N-2$. A simple
functional equation holds for the transformation
 $q \leftrightarrow q^{-1}$.

\subsection{Surfaces: dimension $1$ quotients} \label{d1q}
Let $X$ be a nonsingular, simply connected{\footnote{The referee pointed out that our results hold under the weaker assumption $b_1(X)=0$.}}, projective surface,
and let
$D \in A^1(X)$
be a divisor class.
As observed in \cite {segre}, the Quot scheme $\mathsf {Quot}_{\,X}(\mathbb C^N, n,D)$ of 
short exact sequences $$0\to S\to \mathbb C^N\otimes \mathcal O_X\to Q\to 0,\,\,\,\ \  \chi(Q)=n\, ,\  c_1(Q)=D\, ,\ \text{rank}(Q)=0$$  
carries a 2-term perfect obstruction theory and a virtual fundamental
class of dimension $$\chi(S, Q)=Nn+D^2.$$ Indeed, the 
higher obstructions vanish
 \begin{equation*}
\text{Ext}^2(S, Q)=\text{Ext}^0(Q, S\otimes K_X)^\vee=0\, ,
\end{equation*}
since $Q$ is a torsion sheaf.
Using the above obstruction theory, we define 
generating series of virtual Euler characteristics.

\begin{definition}\label{ttt222} Let $X$ be a 
nonsingular,
  simply connected {\footnote{There is no difficulty to define the generating series
      in the non-simply connected
case, but then $D$ should be taken in $H^2(X,\mathbb{Z})$ instead of $A^1(X)$.}}, projective
 surface.
For a divisor class $D\in A^1(X)$  and an integer $N\geq 1$, let
$$\mathsf{Z}^{\mathcal{E}}_{X, N, D}(q)  = \sum_{n\in \mathbb{Z}} q^n
\mathsf e^{\mathrm{vir}}(\mathsf {Quot}_{\,X}(\mathbb C^N, n, D))\, .$$
\end{definition}

For fixed $N$ and $D$, 
 the Quot schemes $\mathsf {Quot}_{\,X}(\mathbb C^N, n, D)$
are empty for all $n$ sufficiently negative, so
$$ \mathsf{Z}^{\mathcal{E}}_{X,N, D}(q)\, \in \, \mathbb{Z}((q))\, .$$ 
The virtual Euler characteristic results described in Section \ref{vsd0} concern the 
generating series
$\mathsf{Z}^{\mathcal{E}}_{X,N, 0}(q)$ 
with vanishing divisor class $D$.
In case $D\neq 0$,  
exact calculations are  more difficult to obtain. 

\vspace{10pt}
\noindent{\bf (i) Rational surfaces}
\vspace{10pt}

A very rich theory arises for rational surfaces. In Section \ref{ves1}, we 
write general tautological integrals over Hilbert schemes of points which
compute the virtual Euler characteristics. The following result
provides an example of an exact solution.
 
\begin{proposition}\label{t6} Let $X$ be the blowup of a rational surface
with exceptional divisor $E$. We have
$$\mathsf{Z}^{\mathcal{E}}_{X,1, E}(q)= q\left(\frac{(1-q)^{2}}{1-2q}\right)^{K_X^2+1}\, .$$
\end{proposition}

The formula of  Proposition \ref{t6} concerns only the case $N=1$. The proof makes use again of Theorem \ref{thm1} for curves. Further exact calculations
for rational 
surfaces will require new techniques. However, we can calculate much more for
$K3$ surfaces and surfaces of general type.

\vspace{10pt}
\noindent {\bf (ii) $K3$ surfaces}
\vspace{+10pt}

For $K3$ surfaces, the standard obstruction theory contains a trivial factor 
which forces the virtual invariants to vanish. The natural generating
series
therefore concerns the virtual Euler characteristics of 
the {\it reduced} obstruction theory:
$$\mathsf{Z}^{\mathsf{r}\mathcal{E}}_{X, N, D}(q)  = \sum_{n\in \mathbb{Z}} q^n
\mathsf e^{\mathrm{red}}(\mathsf {Quot}_{\,X}(\mathbb C^N, n, D))\, .$$

In the $N=1$ case,
the reduced obstruction theory leads to expressions matching the curve counts on $K3$ surfaces. 
Specifically, let $N_{g, n}$ be defined by the Kawai-Yoshioka  \cite {KY}
formula: \begin{equation}\label{kawaiy}\sum_{g=0}^{\infty}\sum_{n=1-g}^{\infty} \,N_{g, n} \,y^n \,q^g=\left(\sqrt{y}-\frac{1}{\sqrt{y}}\right)^{-2}\prod_{n=1}^{\infty} \frac{1}{(1-q^n)^{20}(1-yq^n)^2(1-y^{-1}q^n)^2}\, .\end{equation} 
The Kawai-Yoshioka formula has played a central role in the Gromov-Witten and
the stable pairs theory of $K3$ surfaces \cite{MPT, PT22,PT}. For primitive
classes, we have complete results.

\begin{theorem}\label{t7} 
Let $X$ be a $K3$ surface, and 
let $D$ be a primitive divisor class of genus
$2g-2= D^2$
which is big and nef. We have 
$$\mathsf e^{\mathrm{red}} (\mathsf {Quot}_{\,X}(\mathbb C^1, n, D))= N_{g, n}\, .$$
\end{theorem} 

\noindent The argument matches the reduced virtual Euler characteristic integral of the Quot scheme to 
 the topological Euler characteristic integral
of the moduli space of stable pairs (the integrands however are
{\it not} the same).
\vskip.1in

\vspace{5pt}
\noindent {\bf (iii) Surfaces of general type}
\vspace{10pt}

Let $X$ be a 
simply connected surface of general type with $p_g>0$. 
In the class of the canonical divisor
$K_X$,
 we show the vanishing of the virtual Euler characteristics for $N=1$ 
in almost all cases.
 The single exception is significant: the
 Poincar\'e-Seiberg-Witten invariants of \cite {CK,DKO} are recovered,
 $$\mathsf e^{\mathrm{vir}} (\mathsf {Quot}_{\,X}(\mathbb C^1, n=-K_X^2, K_X))
=(-1)^{\chi(\mathcal O_X)}\, .$$
 For arbitrary $N$, a vanishing holds for minimal surfaces.{\footnote{We thank
M. Kool for very helpful discussions about Seiberg-Witten classes.}}

\begin{proposition}\label{pro22} Let $X$ be a simply connected minimal surface of general type with $p_g>0$. If $D$ is a curve class with $$\left[\mathsf {Quot}_{X} (\mathbb C^N, n, D)\right]^{\mathrm{vir}}\neq 0\, ,$$ then $D=\ell K_X$ for $0\leq \ell\leq N$. 

\end{proposition}

If we further assume the canonical class of $X$ is 
represented by a nonsingular curve, we can
calculate
$\mathsf Z^{\mathcal E}_{X, \,N, \,\ell K_X}(q)$ completely
in all cases. By Proposition \ref{pro22}, we need only
consider $\ell$ in the range
$$0\leq \ell \leq N\, .$$

\begin{theorem} \label{genttt}
Let $X$ be a simply connected minimal surface of general type with
a nonsingular canonical curve
of genus $g=K_X^2+1$. 
Then, 
$$\mathsf Z_{X, \,N, \,\ell K_X}^{\mathcal E}(q)=(-1)^{\ell\cdot\chi(\mathcal O_X)} \ q^{\ell(1-g)}\cdot \sum_{1\leq i_1<\ldots<i_{N-\ell}\leq N} \mathsf A(r_{i_1}, \ldots, r_{i_{N-\ell}})^{1-g}\, ,$$ 
where the sum is taken over all 
 $\binom{N}{N-\ell}$
choices of $N-\ell$ distinct roots $r_i(q)$ 
of the polynomial equation 
$$z^N-q(z-1)^N=0$$
in the variable $z$. The function $\mathsf{A}$ is
defined by
$$\mathsf A(x_1, \ldots, x_{N-\ell})= \frac{(-1)^{\binom{N-\ell}{2}}}{N^{N-\ell}} \cdot \prod_{i=1}^{N-\ell} \frac{(1+x_i)^N (1-x_i)}{x_i^{N-1}}\cdot \prod_{i<j} \frac{(x_i-x_j)^2}{1-(x_i-x_j)^2}.$$ 
 \end{theorem}

Since the answer of Theorem \ref{genttt} is
a symmetric function of the roots  $r_1(q),\ldots, r_N(q)$, we have
$$\mathsf Z^{\mathcal E}_{X, N, \ell K_X}(q)\, \in \, \mathbb{Q}(q)\, .$$
Theorem \ref{genttt} is the most advanced calculation of
paper. The proof uses essentially all of the ideas and methods that we 
have developed.

\begin{example}
Theorem \ref{genttt}  for $N=2$ and $\ell=1$
specializes to the following
formula:
\begin{multline}\label{jj559}
\mathsf Z^{\mathcal E}_{X, 2, K_X}(q)
=\\
(-1)^{\chi(\mathcal{O}_X)} \left(\frac{q}{2}\right)^{1-g}
\left(\left( 
\frac{(1+r_1)^2(1-r_1)}{r_1}\right)^{1-g}+ \left( \frac{(1+r_2)^2(1-r_2)}{r_2}\right)^{1-g}\right)\, ,
\end{multline}
 where  $r_1(q)$ and $r_2(q)$ are the two roots of the quadratic equation
 $$z^2-q(z-1)^2=0 $$
in the variable $z$.
For a minimal surface of general type $X$ with a
canonical curve of genus $2$, formula \eqref{jj559} yields:
$$\mathsf Z^{\mathcal E}_{{X}, 2, K_{X}}(q)= 
(-1)^{\chi(\mathcal{O}_{X})}\,
\frac{(16q-8)}{(1-4q)^2}\, .$$
For $X$ with a canonical curve of genus $3$, the answer is
$$\mathsf Z^{\mathcal E}_{X, 2, K_X}(q)= (-1)^{\chi(\mathcal{O}_X)}\,
\frac{(128q^4-64q^3+8q^2-16q+8)}{q(1-4q)^4}\, .$$
\end{example}

\subsection{Rationality}
By Theorem \ref{rattt}, the series $\mathsf{Z}^{\mathcal{E}}_{X,N,0}(q)$ 
is the expansion of a
 rational function in  $q$. Rationality also holds for
all the examples discussed in Section \ref{d1q} for Quot schemes
of quotients with dimension 1 support on surfaces. 

\begin{conjecture}\label{qrr} For a nonsingular, simply connected, projective
surface $X$,  $$\mathsf{Z}^{\mathcal{E}}_{X,N,D}(q) \in \mathbb{Q}(q)\, .$$
\end{conjecture}

A natural further direction is to study the
associated series in algebraic cobordism:
$$\mathsf{Z}^{\mathsf{Cobord}}_{X,N,D}\, = \, \sum_{n\in \mathbb{Z}} [\mathsf {Quot}_{\,X}(\mathbb C^N, n, D)]^{\mathrm{vir}} q^n
\, \in\,  \Omega_*(\mathsf{point})((q))\, .$$
The algebraic cobordism{\footnote{See \cite{LM} for a foundational
treatment of algebraic cobordism
and \cite{LevP} for applications to enumerative geometry.}}
class $$[\mathsf {Quot}_{\,X}(\mathbb C^N, n, D)]^\mathrm{vir} \in \Omega_*(\mathsf{point})$$
is well-defined by \cite{Sh}.  Are there formulas for
$\mathsf{Z}^{\mathsf{Cobord}}_{X,N,D}(q)$?

The parallel question for the virtual classes in algebraic cobordism of
the moduli spaces of stable pairs on 3-folds is
conjectured to have an affirmative answer, see \cite[Conjecture 0.3]{Sh}.
In the case of toric geometries, Shen is able to prove the rationality
of the cobordism series via the rationality results for the
descendent theory of stable pairs \cite{PP1,PP2}.

\subsection{Vafa-Witten theory}
There has been a series of recent papers studying the
virtual Euler characteristics of the moduli spaces of
stable bundles (and stable Higgs pairs) on surfaces \cite{GSY, GT, 
GK1, GK2, La, TT1, TT2}.
The outcome has been a clear mathematical proposal for the 
theory studied earlier by Vafa and Witten \cite{VW}.

Definition \ref{ttt222} here is motivated by the Vafa-Witten developments.
The Quot scheme geometry, with the associated obstruction theory, provides
a straightforward approach to sheaf counting on surfaces. The idea
is that given a stable bundle $B$ of rank $N$ on an algebraic surface $X$, we
can pick $N$ sections (assuming $B$ is sufficiently positive) which will
generically generate $B$:
\begin{equation}\label{ff556}
  0\to \mathbb C^N\otimes \mathcal O_X\to B \to F \to 0\, ,
\end{equation}
where $F$ is supported in dimension $1$.
By dualizing \eqref{ff556}, we obtain a quotient sequence
$$\left[0\to B^\vee \to \mathbb C^N\otimes \mathcal O_X\to Q \to 0\right]\,\in\,
\mathsf {Quot}_{\,X}(\mathbb C^N, \chi(Q), c_1(B))
 \, .$$
Of course, $\chi(Q)$
can be computed from the Chern classes of $B$ and $X$.

The calculations that
we have presented, which may be viewed as the beginning of the study of
the virtual Euler characteristics of Quot schemes of surfaces, 
already show some features of  Vafa-Witten theory: the appearance of
the Kawai-Yoshioka formula (in the $K3$
case) and the appearance of the Seiberg-Witten invariants (in the general type case). A difference is the rationality in the variable $q$ for the Quot scheme theory
versus
modularity in the variable 
$$q=\exp(2\pi i\tau)$$ for Vafa-Witten theory. 
A basic open question is the following.

\begin{question}\label{q24}
Formulate the precise relationship of the Quot scheme theory of surfaces
for all $N$ to
Vafa-Witten theory and Seiberg-Witten theory.
\end{question}

Moduli spaces of bundles on {\it curves} with sections have been considered by many authors, see \cite {BDW, Th}. 
Moreover, the relationship between the intersection theory
of Quot schemes and the moduli space of
stable bundles on curves has been   
successfully studied in \cite{MarianCr, MarOp}. 

\subsection{Higher rank quotients}
Let $X$ be a nonsingular projective surface, and consider
the Quot scheme
$\mathsf {Quot}_{\,X}(\mathbb C^N, n,D,r)$
of quotients with dimension 2 support,
$$0\to S\to \mathbb C^N\otimes \mathcal O_X\to Q\to 0,\,\,\,\ \  \chi(Q)=n\, ,\  c_1(Q)=D\, ,\ \text{rank}(Q)=r>0\, .$$ 
The existence of a 
virtual fundamental class of $\mathsf {Quot}_{\,X}(\mathbb C^N, n,D,r)$ 
for del Pezzo surfaces  was first noted in \cite {Sc}, but the
study can be pursued more generally.

As  in the cases of support of dimension $0$ and $1$, 
the higher obstructions of the standard deformation theory 
of $\mathsf {Quot}_{\,X}(\mathbb C^N, n,D,r)$
lie in $\text{Ext}^2(S,Q)$.
We have 
$$\text{Ext}^2(S,Q) = \text{Ext}^0(Q,S\otimes K_X)^\vee \ \ \
{\text \em and} \ \ \
\text{Ext}^0(Q,S\otimes K_X) 
 \hookrightarrow
  \text{Ext}^0(Q,\mathbb{C}^N \otimes K_X)\, .$$
Hence, if $\text{Ext}^2(S,Q)\neq 0$, then
$\text{Ext}^0(Q,\mathbb{C}^N \otimes K_X)\neq 0$.
Since $Q$ is generated by global sections, 
we conclude that $H^0(X,K_X)\neq 0$.

By the above logic, we obtain the following condition: {\it if
$X$ satisfies
\begin{equation}\label{ff9922}
H^0(X,K_X)=0\, ,
\end{equation}
then the standard deformation theory 
of $\mathsf {Quot}_{\,X}(\mathbb C^N, n,D,r)$ is
2-term and yields a virtual fundamental class of dimension
$\chi(S,Q)$}.

There are many surfaces which satisfy $H^0(X,K_X)=0$ including
rational surfaces, ruled surfaces, Enriques surfaces, and even
some surfaces of general type. 
The Quot scheme
virtual Euler characteristic
theory for such surfaces is well defined for all $r$, $D$, and $n$.
We leave the investigation for higher $r$ to a future paper.

\subsection{Plan of the paper}
We start by computing Segre integrals over the symmetric product $$C^{[n]}=\mathsf {Quot}_{\,C}(\mathbb C^1,n)$$
in Section \ref{ff445}. In particular, Theorem \ref{thm1} is proven in Section \ref{thm1pr}. 
Theorem \ref{thm3} about the
Fuss-Catalan numbers is obtained via Wick expansion
in Section \ref{thm3pr}. 
Segre integrals over Quot schemes of curves for higher $N$
are studied in Section \ref{hhh999} where the proofs of Theorem \ref{thm2} and
the first part of
Theorem \ref{cor1} are presented.

We then consider Quot schemes of surfaces. Section \ref{vis}
concerns the case of quotients with dimension 0 support. The second
part of Theorem \ref{cor1} as well as Theorems  \ref{thm4} and
\ref{rattt} are proven there by reducing surface integrals to  curve integrals.
Section \ref{vis1} concerns the case of dimension 1 support. The proofs of
Theorems \ref{t7} and \ref{genttt} are presented in Sections \ref{k3k3k3} and \ref{jjj999}
respectively.

\subsection{Subsequent developments} Further rationality results concerning the generating series of virtual $\chi_{-y}$-genera of Quot schemes of surfaces were obtained in \cite{Li}. Also in \cite{Li}, some of the assumptions made in this paper (the underlying surface being simply connected, smoothness of the canonical divisor, minimality) were removed by taking advantage of the connections with Seiberg-Witten theory. In parallel, in \cite{JOP}, series of descendant invariants were proven to be rational for several geometries, including the case of rational surfaces when $N=1$. Finally, the virtual $K$-theory of Quot schemes of surfaces is studied in \cite {AJLOP}.

\subsection {Acknowledgements.} 
Many of the ideas developed here are related to collaborations 
with Alina Marian in \cite {MarOp, mo, quotients, MOP2,  segre, MOP}. We are grateful 
to her for numerous discussions over the years related to Quot schemes of curves and surfaces. 

Our study of the virtual Euler characteristics of the
Quot scheme of surfaces was motivated in part by the
Euler characteristic 
calculations of L. G\"ottsche and M. Kool \cite{GK1, GK2} for the
moduli spaces of rank $2$ and $3$ stable sheaves on surfaces.
We thank I. Gessel, B. Rhoades, and
J. Verstraete for helpful
conversations about the Catalan numbers.
Discussions about several related topics with
 A. Okounkov and R. Thomas have been valuable.

D. O. was supported by the NSF through grant DMS 1802228. R.P. was supported by the Swiss National Science Foundation and
the European Research Council through
grants SNF-200020-182181, 
ERC-2017-AdG-786580-MACI, SwissMAP, and the Einstein Stiftung.  
We thank the Shanghai Center for Mathematical Science at Fudan
University for
a very productive visit in September 2018 at the start of the
project.
R.P. also thanks the
Hebrew University of Jerusalem where part of the research 
was pursued during a visit in October 2018.

The project has received funding from the European Research
Council (ERC) under the European Union Horizon 2020 Research and
Innovation Program (grant No. 786580).

\section{Symmetric products of curves} \label{ff445} 
\subsection{Overview} We first present
the proof of Theorem \ref{thm1}. Theorem \ref{thm3} will be obtained 
in Section \ref{thm3pr}
by specializing Theorem \ref{thm1} to genus $1$. 
In fact, {\it all} other main results of the paper 
(Theorems \ref{thm2}, \ref{cor1}, \ref{thm4}, \ref{rattt}  and \ref{genttt}) 
proven in later sections, rely either directly upon Theorem \ref{thm1} or 
upon the analysis of the integrals over $C^{[n]}$ developed here.

\subsection{Projective line} 
To begin the proof of Theorem \ref{thm1}, we observe that
the factorization \begin{equation}\label{fact} \sum_{n=0}^{\infty} q^n \int_{C^{[n]}} s_{x_1}(\alpha_1^{[n]})\cdots s_{x_\ell} (\alpha_\ell^{[n]})= {\mathsf A}_1^{c_1(\alpha_1)} \cdots {\mathsf A}_\ell^{c_1(\alpha_\ell)} \cdot \mathsf B^{1-g}\, ,\end{equation} allows us to specialize
the calculation to genus 0 where 
$$C\simeq \mathbb P^1 \ \ \ \text{and}\ \ \ C^{[n]}\simeq \mathbb P^n\, .$$
 We write $h$ for the hyperplane class on $\mathbb P^n$.  
 
\begin{lemma} \label{segr} For a $K$-theory class $\alpha$ on 
$\mathbb P^1$ of rank 
$r$ and degree $d=c_1(\alpha)$ we have $$s_{x} (\alpha^{[n]})= (1-xh)^{d-nr+r}\, .$$
 \end{lemma}
 \proof Both expressions are multiplicative in short exact sequences of vector bundles $$0\to V_1\to V\to V_2\to 0\, .$$
The claim is clear for the right hand side. 
For the left hand side, claim is a consequence
of the induced sequence $$0\to V_1^{[n]}\to V^{[n]}\to V_2^{[n]}\to 0\implies s_x(V^{[n]})=s_x(V_1^{[n]})\cdot s_x(V_2^{[n]})\,.$$ 
Since the $K$-theory of $\mathbb P^1$ is generated by line bundles,
we can restrict to $\alpha=\mathcal O_{\mathbb P^1}(d)$. 
By the proof of Theorem $2$ in \cite {MOP}, we have $$\text{ch }((\mathcal O_{\mathbb P^1}(d))^{[n]})=(d+1)-(d-n+1) \exp (-h)\, ,$$ which then gives $$s_x((\mathcal O_{\mathbb P^1}(d))^{[n]})=(1-xh)^{d-n+1}\, ,$$ 
completing the argument. \qed \vskip.1in
 
\subsection { Proof of Theorem \ref{thm1} (using $\mathbb P^1$)}
\label{thm1pr}
 Let $\alpha_1, \ldots, \alpha_\ell$ be $K$-theory classes of ranks $r_i$ and degree $d_i$. Using Lemma \ref{segr}, we obtain \begin{eqnarray*}\mathsf Z_{\mathbb P^1,1}(q, x_1, \ldots, x_\ell\,|\,\alpha_1, \ldots, \alpha_\ell)&=&\sum_{n=0}^{\infty} q^n \int_{\mathbb P^n} s_{x_1}(\alpha_1^{[n]})\cdots s_{x_\ell} (\alpha_\ell^{[n]})\\ &=& \sum_{n=0}^{\infty} q^n \int_{\mathbb P^n} (1-x_1h)^{d_1-nr_1+r_1}\cdots (1-x_\ell h)^{d_\ell-nr_\ell+r_\ell}\\ &=& \sum_{n=0}^{\infty} q^n \int_{\mathbb P^n} f(h)^n\cdot g(h)\\ &=& \sum_{n=0}^{\infty} q^n\cdot \left(\left[t^n\right] f(t)^{n}\cdot g(t)\right). \end{eqnarray*} 
In the third equality, 
\begin{eqnarray*}
f(t)&=& (1-x_1t)^{-r_1} \cdots (1-x_\ell t)^{-r_\ell}\, , \\
g(t)&=& (1-x_1t)^{d_1+r_1} \cdots (1-x_\ell t)^{d_\ell+r_\ell}\, .
\end{eqnarray*}
 The brackets denote the coefficient of the suitable power of $t$. 

We can evaluate such expressions using the Lagrange-B\"urmann formula \cite {WW}. Assuming $f(0)\neq 0$, for the change of variables $q=\frac{t}{f(t)}$,
 the following
general identity holds \begin{equation}\label{lb}\sum_{n=0}^{\infty} \left(\left[t^n\right] f(t)^{n}\cdot g(t)\right)\cdot q^n=\frac{g(t)}{f(t)}
 \cdot \frac{dt}{dq}\, .\end{equation} We will use the above identity
 repeatedly. 

In our case, the change of variables takes the form
 $$q=t(1-x_1t)^{r_1}\cdots (1-x_\ell t)^{r_{\ell}}\, ,$$ and the Segre series becomes $$\mathsf Z_{\mathbb P^1,1}(q, x_1, \ldots, x_\ell\,|\,\alpha_1, \ldots, \alpha_\ell)=\prod_{i=1}^{\ell} (1-x_i t)^{d_i}\cdot \prod_{i=1}^{\ell} (1-x_it)^{2r_i} \cdot \frac{dt}{dq}\, .$$ 
 Combined with the factorization \eqref{fact}, 
 $$\mathsf Z_{\mathbb P^1,1}(q, x_1, \ldots, x_\ell|\alpha_1, \ldots, \alpha_\ell)= {\mathsf A}_1^{d_1} \cdots  {\mathsf A}_\ell^{d_\ell} \cdot \mathsf 
B\, ,$$
the above calculation yields $$\mathsf A_i(q)=1-x_i t\, ,\,\,\ \ \  \mathsf {B}(q) = \prod_{i=1}^{\ell} (1-x_it)^{2r_i} \cdot \frac{dt}{dq}=\left(\frac{q}{t}\right)^2\cdot \frac{dt}{dq}.$$ This completes the proof of Theorem \ref{thm1}. 
\qed
\vskip.1in
For future use, we also record a formula for the {\it logarithms} of the functions $\mathsf A_i$. Of course, we may take $i=1$ without loss of generality. 
\begin{lemma} \label{loga} We have $$\log \mathsf A_1 = \sum_{n=1}^{\infty} (-1)^n\frac{q^n}{n} \cdot \mathsf a_n$$ where \begin{equation}\label{ll}\mathsf a_n (x_1, \ldots, x_{\ell}) = x_1\cdot \sum_{p_1+\ldots+p_{\ell}=n-1} \binom{-nr_1-1}{p_1}
\cdot \binom{-nr_2}{p_2}\cdots \binom{-nr_\ell}{p_\ell} \cdot x_1^{p_1}\cdots x_\ell^{p_{\ell}}\, .\end{equation}
\end{lemma} 

\proof The argument consists in another application of the Lagrange-B\"urmann formula \eqref{lb}. Indeed, write $\widetilde {\mathsf a}_n$ for the right hand side of \eqref{ll} and let 
$${\mathsf L}(q)=\sum_{n=1}^{\infty} (-1)^n\frac{q^n}{n} \cdot \mathsf {\widetilde a}_n\, .$$ 
We must prove $\log \mathsf A_1={\mathsf L}$.  Clearly, 
 \begin{eqnarray*}\widetilde{\mathsf a}_n(x_1, \ldots, x_{\ell}) &=& (-1)^{n-1} x_1 \cdot \left(\left [t^{n-1}\right] (1-x_1t)^{-nr_1-1} \cdot (1-x_2t)^{-nr_2}
\cdots (1-x_{\ell} t)^{-nr_\ell}\right)\\ &=& (-1)^{n-1} x_1 \cdot \left(\left[t^{n-1}\right] f(t)^{n-1}\cdot h(t)\right).\end{eqnarray*}
where we write as before $$f(t)= (1-x_1t)^{-r_1} \cdots (1-x_\ell t)^{-r_\ell}$$ $$h(t)=(1-x_1t)^{-r_1-1}\cdot (1-x_2t)^{-r_2}\cdots (1-x_\ell t)^{-r_{\ell}}\,.$$ 
We further compute \begin{eqnarray*}\frac{d\mathsf L}{dq}&=& \sum_{n=1}^{\infty} (-1)^nq^{n-1} \cdot \widetilde {\mathsf a}_n\\ &=& - x_1 \cdot \sum_{n=1}^{\infty} q^{n-1}\cdot \left(\left[t^{n-1}\right] f(t)^{n-1}\cdot h(t)\right)\\ &=& -x_1 \cdot \frac{h(t)}{f(t)} \cdot \frac{dt}{dq}\\ &=&-\frac{x_1}{1-x_1t}\cdot \frac{dt}{dq}\, ,\end{eqnarray*} where the Lagrange-B\"urmann formula \eqref{lb} was applied in the third equality, 
for the same change of variables $q=\frac{t}{f(t)}$ which we used previously. Therefore $$d \mathsf L= -\frac{x_1}{1-x_1t} dt\ \ \implies\ \  \mathsf L=\log (1-x_1t).$$ Combined with Theorem \ref{thm1}, we obtain 
$\mathsf L=\log \mathsf A_1$.
\qed

\subsection{Wick's formalism for an elliptic curve}
Let $C$ be a nonsingular genus $1$ curve. Let
 $L_1, \ldots, L_\ell$ be line bundles on $C$ of degrees 
$d_1, \ldots, d_{\ell}$. We lift the integrals over the symmetric product 
to the $n$-fold ordinary product via the morphism 
$$p_n: C^{\times n}=C\times \cdots \times C\to C^{[n]}\, .$$ 
We write $D_{ij}$ for the diagonals $$D_{ij}=\{x_i=x_j\}\subset C^{\times n}$$ and further set $$\Delta_i=D_{1, i}+D_{2, i}+\ldots+D_{i-1, i}\, .$$ 
We also write $\pi_i:C^{\times n}\to C$ for the canonical projections, $1\leq i\leq n$. 
From the exact sequence 
$$0\to \pi_n^{\star}L(-\Delta_n)\to p_n^{\star}L^{[n]}\to (\pi_1\times \cdots 
\times\pi_{n-1})^{\star} p_{n-1}^{\star}L^{[n-1]}\to 0\, ,$$ 
we inductively obtain $$p_n^{\star}s_{x}(L^{[n]})=\prod_{i=1}^{n} \frac{1}{1+x(\pi_i^{\star} c_1(L)-\Delta_i)}\, .$$
Consequently, $$\int_{C^{[n]}} s_{x_1}(L_1^{[n]})\cdots s_{x_{\ell}} (L_{\ell}^{[n]})=\frac{1}{n!} \int_{C^{\times n}} \prod_{j=1}^{\ell} \prod_{i=1}^{n} \frac{1}{1+x_j(\pi_i^{\star} c_1(L_j)-\Delta_i)}\, .$$ 
By Lemma \ref{loga}, we know  \begin{multline}\label{logb}\log \left(\sum_{n=0}^{\infty} \frac{q^n}{n!} \int_{C^{\times n}} \prod_{j=1}^{\ell} \prod_{i=1}^{n} \frac{1}{1+x_j(\pi_i^{\star} c_1(L_j)-\Delta_i)}\right)\\ =  \sum_{n=1}^{\infty} (-1)^n\frac{q^n}{n} \cdot \left(\sum_{j=1}^{\ell} \mathsf a_n^{(j)} \cdot \deg L_j\right)\, ,\end{multline} where $$\mathsf a^{(1)}_n (x_1, \ldots, x_{\ell}) = x_1\cdot \sum_{p_1+\ldots+p_{\ell}=n-1} \binom{-n-1}{p_1}\cdot \binom{-n}{p_2}\cdots \binom{-n}{p_\ell} \cdot x_1^{p_1}\cdots x_\ell^{p_{\ell}}\, $$
and $\mathsf a^{(j)}_n (x_1, \ldots, x_{\ell})$ is given by 
the correspondingly permuted formula.

We will expand the left hand side of
 \eqref{logb} using Wick's formalism. To connect with Theorem \ref{thm3}, write $$w_n(p_1, \ldots, p_\ell) = \sum_{T} \mathsf{wt}(T)$$ 
weighted count of ordered $\ell$-colored trees of type
$(p_1,\ldots,p_\ell)$ with
 $$n=p_1+\ldots+p_{\ell}+1\, .$$ 
Theorem \ref{thm3} is equivalent to the following claim:
 \begin{equation}\label{dll3}
w_n(p_1, \ldots, p_\ell) = \frac{(-1)^{n-1}}{n} \binom{-n}{p_1} 
\cdots \binom{-n}{p_\ell}\, . 
\end{equation}
To establish \eqref{dll3}, we set $$\mathsf W_n= \sum_{p_1+\ldots+p_{\ell}=n-1} w_n(p_1, \ldots, p_{\ell}) \cdot x_1^{p_1} \cdots x_\ell^{p_{\ell}}\, .$$ 
Define the differential operator $$\mathsf D_1  = 2x_1 \frac{\partial}{\partial x_1} + x_2 \frac{\partial}{\partial x_2}+\ldots +x_{\ell} \frac{\partial}{\partial x_\ell} + {\bf 1}$$ and define
$\mathsf D_2, \ldots, \mathsf D_{\ell}$
by the correspondingly permuted formulas.

\begin{lemma} \label{wick} The following identity holds $$\log \left(\sum_{n=0}^{\infty} \frac{q^n}{n!} \int_{C^{\times n}} \prod_{j=1}^{\ell} \prod_{i=1}^{n} \frac{1}{1+x_j(\pi_i^{\star} c_1(L_j)-\Delta_i)}\right)=-\sum_{n=1}^{\infty} \frac{q^n}{n} \cdot \left(\sum_{j=1}^{\ell} x_j \cdot \deg L_j\cdot \mathsf D_j \mathsf W_n\right).$$  
\end{lemma}
\proof
We refer the reader to Section 1.3 of \cite{PP} for a gentle
introduction to the Wick formalism in precisely the context which
we require here. 
By Wick, the logarithm on the left hand side
is given by 
\begin{equation}\label{gg999}
\sum_{n=1}^\infty \frac{q^n}{n!} S[n] \, ,
\end{equation}
where $S[n]$ is the {\it connected contribution} on $n$
vertices. We will 
match the connected contributions $S[n]$ with the
right hand side of Lemma \ref{wick}.

Consider an arbitrary monomial of degree $n$ in the diagonal classes. Such a monomial determines a graph with $n$ vertices, whose edges are given by the diagonal associations. Since $C$ is an elliptic curve, the squares of diagonals vanish
$$D_{i,j}^2= 0 \ \in H^*(C^{\times n},\mathbb{Z})\, .$$ 
Hence, a connected graph on $n$ vertices cannot have any cycles, thus it corresponds
exactly to a tree with $n-1$ edges determined by the
diagonals. The diagonals come from the expansions of the terms 
$$\prod_{j=1}^{\ell} \prod_{i=1}^{n} \frac{1}{1+x_j(\pi_i^{\star} c_1(L_j)-\Delta_i)}\, ,$$
and therefore may be considered as carrying colors
between $1, \ldots, \ell$ depending on the $j$ index. 

Let us first analyze the (simpler) connected contribution 
for the 
 terms
\begin{equation}\label{ff99t}
\prod_{j=1}^{\ell} \prod_{i=1}^{n} \frac{1}{1+x_j(-\Delta_i)}\, .
\end{equation}
We see that the
coefficient of $x_1^{p_1} \cdots x_\ell^{p_\ell}$ 
in the connected contribution  with
$$n= p_1+\ldots +p_\ell +1 $$
vertices
is exactly 
a sum over  
{\it labelled $\ell$-colored trees of type $(p_1,\ldots, p_\ell)$}.
The vertices of the trees $T$ are labelled by the $n$ ordered
factors of $C^{\times n}$.
To calculate the weight, we must expand \eqref{ff99t} as
\begin{equation}\label{ff99t9}
\prod_{j=1}^{\ell} \prod_{i=1}^{n} (1+ x_j\Delta_i +x_j^2 \Delta_i^2 +
x_j^3 \Delta_i^3 + \ldots)\, .
\end{equation}
If the $i^{th}$ vertex $v$ of $T$ has $d^j_v$ downward edges colored $j$, 
the weight
receives a factor of $d^j_v!$ since the coefficient of 
the monomial in the corresponding diagonal in
$$x_j^{d^j_v} \Delta_i^{d^j_v}$$
is exactly $d^j_v!$ . Hence, the full weight is
\begin{equation}\label{kk11kk}
 \prod_{v\, \text{vertex}} d_v^1! \cdots  d_v^k! =
\mathsf{wt} (T) \cdot (n-1)!
\ .
\end{equation}

The actual connected contribution $S[n]$ of \eqref{gg999},
which we must calculate,
also includes the insertions of $\pi_i^\star(c_1(L_j))$.
Since the diagonal edges already cut $C^{\times n}$ to
just a single elliptic curve $C$, exactly one
 insertion from the set
$$\big\{ \,  \pi_i^\star(c_1(L_j))\, \big\}_{1\leq i\leq n\, , \ 
1\leq j\ \leq \ell}$$ 
must be chosen.
We separate the contribution 
$$S[n]= \sum_{j=1}^\ell S[n,j]$$
by which $L_j$ is chosen as an insertion.
The connected contribution $S[n,j]$ will be matched with
$$- x_j \cdot \deg L_j\cdot \mathsf D_j \mathsf W_n\cdot (n-1)!$$
to complete the proof.

To this end, we calculate the effect of the insertion 
$L_j$ on the weight of a labelled $\ell$-colored tree of type 
$(p_1,\ldots, p_\ell)$ generated by the diagonals. 
The insertion $L_j$ can occur at any vertex $1\leq i \leq n$.
When $\pi_i^\star(c_1(L_j))$ is selected, the weight 
receives the factor $$-\deg L_j \cdot (d^j_v+1)!$$ since the coefficient 
of the corresponding monomial in
$$x_j^{d^j_v+1} \big(-\pi_i^\star(c_1(L_j))+\Delta_i\big)^{d^j_v+1}$$
is exactly $(d^j_v+1)!$. Since the insertion $L_j$ can be placed at any vertex $1\leq i \leq n$,
we must modify the weight \eqref{kk11kk} of $T$ by the 
prefactor
$$\sum_{i=1}^n (d^j_v +1)=p_j +n = 2p_j + \sum_{j' \neq j} p_{j'}+1\, . $$
This prefactor is achieved precisely by the action of the 
differential operator $\mathsf{D}_j$ on $\mathsf{W}_n$.

Collecting all terms, we obtain 
\begin{eqnarray*}
\sum_{n=1}^\infty \frac{q^n}{n!} S[n]  & = &
\sum_{n=1}^\infty\sum_{j=1}^\ell  \frac{q^n}{n!} S[n,j] \\
& = & \sum_{n=1}^\infty \sum_{j=1}^\ell
\left( -\frac{q^n}{n}  x_j \cdot \deg L_j\cdot \mathsf D_j \mathsf W_n\right)\, 
\end{eqnarray*}
which completes the calculation.
\qed
\vskip.1in

\subsection{Proof of  Theorem \ref{thm3} (using an elliptic curve)}
\label{thm3pr}
We prove Theorem \ref{thm3} here  
geometrically by specializing Theorem \ref{thm1} and Lemma \ref{loga} to genus $1$ and using the Wick result of Lemma \ref{wick}. Alternatively,
a direct combinatorial proof of Theorem \ref{thm3} is provided in
the Appendix.

By setting
the right hand side of \eqref{logb} equal to
the right hand side  of the formula of Lemma \ref{wick},
we obtain \begin{equation}\label{logc}x_1\mathsf D_1 \mathsf W_n = (-1)^{n-1} \mathsf a_n^{(1)}\, .\end{equation} The operator
 $\mathsf D_1$ acts on the monomial $x_1^{p_1}
\cdots x_{\ell}^{p_{\ell}}$ as multiplication by $(n+p_1)$.
By  matching coefficients of $x_1^{p_1}\cdots x_{\ell}^{p_{\ell}}$ on both sides of \eqref{logc}, 
we solve $$w_n(p_1, \ldots, p_{\ell})= \frac{(-1)^{n-1}}{n} \binom{-n}{p_1} 
\cdots \binom{-n}{p_\ell}\, ,$$ 
which completes the argument.\qed 

\section{Quot schemes of curves for higher $N$} \label{hhh999}
\subsection{Overview}
We prove here  Theorem \ref{thm2}, part of Theorem \ref{cor1}, and the
associated Corollaries \ref{c2} and \ref{cor2}.

We begin with Theorem \ref{thm2}. We specialize directly to the case of an elliptic curve $C$, seeking to show that
$$\sum_{n=0}^{\infty} q^n \int_{\mathsf {Quot}_{\,C} (\mathbb C^N, n)} s(V^{[n]})=\mathsf A(q)^{\deg V}\, ,$$ 
with the specified formula for $\mathsf A(q)=\mathsf{A}_{1,r,N}(q)$.

\subsection{Equivariant localization} \label{el}
The nonsingular projective variety
 ${\mathsf {Quot}_{\,C} (\mathbb C^N, n)}$ carries a
natural action of the algebraic torus $\mathbb C^{\star}$ defined
as follows.
Let $\mathbb C^{\star}$ act diagonally on $\mathbb C^N$ with weights 
$$w_1<w_2<\ldots<w_N\, .$$
The $\mathbb C^{\star}$-action on
$\mathsf {Quot}_{\,C}(\mathbb C^N, n)$ is then induced
via  the associated $\mathbb C^{\star}$-action on the middle term of the exact sequence $$0\to S\to \mathbb C^N\otimes \mathcal O_C\to Q\to 0\, .$$

We will prove Theorem \ref{thm2} by applying the
Atiyah-Bott $\mathbb C^{\star}$-equivariant localization formula
 to compute integrals over $\mathsf {Quot}_{\,C}(\mathbb C^N, n)$.
The fixed loci are indexed by partitions $n_1+\ldots+n_{N}=n$ 
where  
$$\mathsf F[n_1, \ldots, n_N]=C^{[n_1]}\times \ldots \times C^{[n_N]}$$ 
parameterizes tuples $(Z_1, \ldots, Z_N)$ of divisors on $C$ with 
$$\text{length}(Z_i)=n_i\, .$$ 
The inclusion $$j:\mathsf F[n_1, \ldots, n_N]\hookrightarrow \mathsf {Quot}_{\,C} (\mathbb C^N, n)$$ corresponds to the invariant sequences $$0\to S=\bigoplus_{i=1}^{N} \mathcal O_C(-{Z_i})\hookrightarrow \bigoplus_{i=1}^{N} \mathcal O_C\to Q=\bigoplus_{i=1}^{N} \mathcal O_{Z_i}\to 0\, .$$ 
The normal bundle to the fixed locus is found from the moving part of the tangent bundle: 
\begin{eqnarray*}\mathsf N[n_1, \ldots, n_N]&=&\text{Hom} (S, Q)^{\text{mov}}=\bigoplus_{i\neq j} \text{Hom}(\mathcal O(-Z_i), \mathcal O_{Z_j})[w_j-w_i]\\&=& \bigoplus_{i\neq j}H^{\bullet}(\mathcal O(Z_i)|_{Z_j})[w_j-w_i]\\&=& \bigoplus_{i\neq j}\left(H^{\bullet}(\mathcal O(Z_i))-H^{\bullet}(\mathcal O(Z_i-Z_j))\right)[w_j-w_i]\end{eqnarray*} with the brackets denoting the equivariant weights. 
We combine the mixed $(i, j)$ and $(j, i)$ terms by
setting $$\mathbb V_{ij}=H^{\bullet}(\mathcal O(Z_i-Z_j))[w_j-w_i]\oplus H^{\bullet} (\mathcal O(Z_j-Z_i))[w_i-w_j]\, .$$ 
Since $C$ is an elliptic curve, Serre duality yields
the $\mathbb C^\star$-equivariant isomorphism
 $$\mathbb V_{ij}\simeq \mathbb V_{ij}^{\vee}[-1]\, .$$ 
 Therefore, $$\mathsf e_{\mathbb C^{\star}}(\mathbb V_{ij})=(-1)^{ \chi(\mathcal O(Z_i-Z_j))}=(-1)^{n_i+n_j}\, .$$ 
For the remaining terms, we use the $K$-theoretic relation $$H^{\bullet}(\mathcal O(Z_i))=H^{1-\bullet}(\mathcal O(-{Z_i}))^{\vee}=-H^{\bullet}(\mathcal O)^{\vee}+H^{0}(\mathcal O_{Z_i})^{\vee}$$ obtained from the
 exact sequence $$0\to \mathcal O(-Z_i)\to \mathcal O\to \mathcal O_{Z_i}\to 0\, .$$ While the first summand is trivial, the second summand corresponds to the bundle $\left(\mathcal O^{[n_i]}\right)^{\vee}.$ We conclude 
 \begin{eqnarray*}\mathsf e_{\mathbb C^{\star}}(\mathsf N[n_1, \ldots, n_N])&=&\prod_{i<j} (-1)^{n_i+n_j} \prod_{i\neq j} \mathsf 
e_{\mathbb C^{\star}}\left((\mathcal O^{[n_i]})^{\vee}[w_j-w_i]\right)\\&=& \prod_{i\neq j} \mathsf e_{\mathbb C^{\star}} (\mathcal O^{[n_i]}[w_i-w_j]).\end{eqnarray*} 
Furthermore, the restriction of $V^{[n]}$ to $\mathsf F[n_1, \ldots, n_N]$ splits equivariantly as $$\iota^{\star} V^{[n]}=V^{[n_1]}[w_1] \oplus \ldots \oplus V^{[n_N]}[w_N]\, .$$ 
Atiyah-Bott localization then yields \begin{eqnarray*}\mathsf Z_{C, N}(q\,|\,V)&=&\sum_{n=0}^{\infty} q^n \int_{\mathsf {Quot}_{\,C} (\mathbb C^N, n)} s(V^{[n]})\\&=&\sum_{n_1+\ldots+n_N=n} q^{n} \int_{C^{[n_1]}\times \ldots \times C^{[n_N]}} \frac{\prod_{i} s (V^{[n_i]}[w_i])}{\prod_{i}\prod_{j\neq i} e_{\mathbb C^{\star}} (\mathcal O^{[n_i]}[w_i-w_j])}\, .\end{eqnarray*}For the expression on the right hand side, we also take the non-equivariant limit $w_1=\ldots=w_N=0$.

An important aspect of the above formula is that, in the genus $1$ case, 
the integral on the right hand side splits over the individual factors. 
For any tuple of equivariant weights $(a, b_1, \ldots, b_{N-1})$, we write $$\mathsf P_C(q\,|\,a\,|\,b_1, \ldots, b_{N-1})= \sum_{n=0}^{\infty} q^n \int_{C^{[n]}} \frac{s (V^{[n]}[a])}{e_{\mathbb C^{\star}} (\mathcal O^{[n]}[b_1])\cdots
 e_{\mathbb C^{\star}} (\mathcal O^{[n]}[b_{N-1}])}\, .$$
We can write the splitting explicitly as
 \begin{multline} \label{product} \mathsf Z_{C, N} (q\,|\,V)=\\ 
\mathsf P_C(q\,|\,w_1\,|\,w_1-w_2, \ldots, w_1-w_N)\cdots
 \mathsf P_C(q\,|\,w_{N}\,|\,w_N-w_1, \ldots, w_N-w_{N-1})\, .\end{multline} 
In fact, equation \eqref{product} holds equivariantly. 
To prove Theorem \ref{thm2}, we must take the non-equivariant limit: 
we must extract the free term with respect to the variables $w_1, \ldots, w_N$
on the right hand side.

\subsection {Symmetric products} Our next step
is to evaluate the expressions 
$$\mathsf P_C(q\,|\,a\,|\,b_1, \ldots, b_{N-1})$$
 by relating them to the integrals of Theorem \ref{thm1}. 
For convenience of notation, we write $$\widehat s_t(V)= t^{-\text{rank } V}\cdot s_{1/t}(V)= \prod_{i} \frac{1}{t+v_i}\, ,$$ where the $v_i$
are the roots of a vector bundle $V$ on a scheme $S$. 

Write $\mathsf R=H^{\star}_{\mathbb C^{\star}}(\text{pt})$ for the equivariant coefficient ring. For $\alpha, \beta_1, \ldots, \beta_{N-1}\in \mathsf R$, we introduce the function $$\mathsf Q_C(q \,| \,\alpha\, |\,\beta_1, \ldots, \beta_{N-1})= \sum_{n=0}^{\infty} q^n \int_{C^{[n]}} \widehat s_{\alpha} (V^{[n]}) \cdot \widehat s_{\beta_1} (\mathcal O^{[n]})\cdots \widehat s_{\beta_{N-1}} (\mathcal O^{[n]})\, .$$ Note that $$\mathsf Q_C(q \,| \,\alpha\, |\,\beta_1, \ldots, \beta_{N-1})\in \mathsf K[[q]]$$ where $\mathsf K$ denotes the fraction field of $\mathsf R$. The calculations below will take place in the power series ring
$\mathsf K[[q]]$.

For a scheme $S$ endowed with a trivial torus action, and a vector bundle $V\to S$ with nontrivial equivariant weight $t$, we have $$\widehat s_t(V)= \mathsf e_{\mathbb C^{\star}} (V[t])^{-1}\in H^{\star}(S)\otimes \mathsf K\, .$$ 
Applied to our setting, we obtain \begin{equation}\label{compare}\mathsf P_C(q\,|\,a\,|\,b_1, \ldots, b_{N-1})=
\mathsf Q_C(q \,| \,1+a \,|\, b_1, \ldots, b_{N-1})\, .\end{equation} The next result computes the logarithm of $\mathsf Q_C$. 

\begin{lemma} \label{refr}
For an elliptic curve $C$, we have $$\mathsf Q_C(q \,| \,\alpha \,|\,\beta_1, \ldots, \beta_{N-1})=\mathsf F(q \,| \,\alpha \,|\,\beta_1, \ldots, \beta_{N-1})^{\,\deg V}\, ,$$ where we define
 $$\log \mathsf F (q \,| \,\alpha\, |\,\beta_1, \ldots, \beta_{N-1})=\sum_{n=1}^{\infty} (-1)^n \frac{q^n}{n}\cdot \mathsf f_n (\alpha\, |\,\beta_1, \ldots, \beta_{N-1})\, ,$$ with 
\begin{multline}\label{fn} \mathsf f_n (\alpha\, |\,\beta_1, \ldots, \beta_{N-1})=\sum_{p+q_1+\ldots+q_{N-1}=n-1}\binom{-nr-1}{p} \binom{-n}{q_1}\cdots \binom{-n}{q_{N-1}}\\ \cdot
\alpha^{-nr-p-1}\beta_1^{-n-q_1}\cdots \beta_{N-1}^{-n-q_{N-1}}\, .\end{multline}
\end{lemma} 
\proof Using the definitions, we compute \begin{eqnarray*}&\mathsf Q_C&(q \,| \,\alpha\, |\,\beta_1, \ldots, \beta_{N-1})=
\sum_{n=0}^{\infty} q^n \int_{C^{[n]}} \widehat s_{\alpha} (V^{[n]}) \cdot \widehat s_{\beta_1} (\mathcal O^{[n]})\cdots \widehat s_{\beta_{N-1}} (\mathcal O^{[n]})\\&=&\sum_{n=0}^{\infty} (q\alpha^{-r}\beta_1^{-1}\cdots \beta_{N-1}^{-1})^n \int_{C^{[n]}} s_{\frac{1}{\alpha}} (V^{[n]}) s_{\frac{1}{\beta_1}} (\mathcal O^{[n]})
\cdots  s_{\frac{1}{\beta_{N-1}}} (\mathcal O^{[n]})\\&=& \mathsf Z_C(\widehat q, \alpha^{-1}, \beta_1^{-1}, \ldots, \beta_{N-1}^{-1}\,|\, V, \mathcal O, \ldots, \mathcal O )\,.\end{eqnarray*} Here we set
 $$\widehat q=q \alpha^{-r} \beta_1^{-1}\cdots \beta_{N-1}^{-1}\,, $$ and we remind the reader that the Segre series $\mathsf Z_{C}$ was introduced in Definition \ref{ssseg}. Since most of the bundles appearing are
trivial, only one universal function appears in the answer:
 $$\mathsf Z_C(\widehat q, \alpha^{-r}, \beta_1^{-1}, \ldots, \beta_{N-1}^{-1}\,|\, V, \mathcal O, \ldots, \mathcal O )=\mathsf F^{\,\deg V}\, .$$ 
The proof is completed by invoking Lemma \ref{loga} which gives an expression for $\log \mathsf F$ matching the one claimed here. \qed
\vskip.1in

\subsection{Proof of Theorem \ref{thm2}} 
By equation \eqref{product}, equation \eqref{compare}, and Lemma \ref{refr}, we obtain $$\mathsf Z_{C, N} (q\,|\, V) =\mathsf A(q)^{\deg V}\, ,$$ 
where  $\log \mathsf A (q)$ equals
\begin{multline*}
\log\, \mathsf F(q \,| \,1+w_1 \,|\,w_1-w_2, \ldots, w_1-w_N) + \ldots + \log\,\mathsf F(q \,| \,1+w_N \,|\,w_N-w_1, \ldots, w_N-w_{N-1})= \\ 
\sum_{n=1}^{\infty} \frac{(-q)^n}{n} \left[\mathsf f_n(1+w_1\, |\,w_1-w_2, \ldots, w_1-w_N) + \ldots + \mathsf f_n( 1+w_N\, |\,w_N-w_1, \ldots, w_N-w_{N-1})\right]\, .
\end{multline*}
Our goal is to prove
$$\log \mathsf A(q)=\sum_{n=1}^{\infty} (-1)^{(N+1)n+1} \binom{(r+N)n-1}{Nn-1}\cdot \frac{q^n}{n}\, .$$ 
Equivalently, we will show
that the free term, with respect to the variables $w_1, \ldots, w_N$, in the expression{\footnote{We use here that taking the free term can be done before or after taking the logarithm.}} 
$$\mathsf f_n(1+w_1\, |\,w_1-w_2, \ldots, w_1-w_N) + 
\ldots + \mathsf f_n(1+w_N\, |\,w_N-w_1, \ldots, w_N-w_{N-1})$$ 
equals $$(-1)^{Nn+1} \binom{(r+N)n-1}{Nn-1}\, .$$ 

To establish the last claim,
we will use the expression for $\mathsf f_n$ provided by equation \eqref{fn}. Each monomial in the formula contributes the following sum to the final answer 
$$(1+w_1)^{-nr-p-1} (w_1-w_2)^{-n-q_1}\cdots (w_1-w_N)^{-n-q_{N-1}}+\ldots $$
 $$+ (1+w_N)^{-nr-p-1}(w_N-w_1)^{-n-q_1}\cdots (w_N-w_{N-1})^{-n-q_{N-1}}\, .$$ 
By Lemma \ref{x1} below, the free term of the sum equals 
$$(-1)^{(n+q_1)+\ldots+(n+q_{N-1})}\cdot \binom{(nr+p)+(n+q_1)+\ldots+(n+q_{N-1})}{nr+p}$$ $$=(-1)^{n(N-1)+q} \cdot \binom{(N+r)n-1}{nr+p}\, ,$$ where $$q=q_1+\ldots+q_{N-1}\implies p+q=n-1\, .$$ Therefore, the free term we seek is 
$$\sum_{p+q_1+\ldots+q_{N-1}=n-1}(-1)^{n(N-1)+q}\binom{-nr-1}{p} \binom{-n}{q_1}
\cdots \binom{-n}{q_{N-1}}\cdot \binom{(N+r)n-1}{nr+p}\, .$$ By the
 Vandermonde identity the middle binomials can be summed: 
$$ \sum_{p+q=n-1} (-1)^{n(N-1)+q}\binom{-nr-1}{p} \binom{-n(N-1)}{q} \binom{(N+r)n-1}{nr+p}\, .$$ After substituting $p=n-1-q$ and rearranging the factorials, 
we obtain $$(-1)^{nN+1}\cdot \frac{(nN-1)!}{(n-1)!(n(N-1)-1)!} \binom{(r+N)n-1}{Nn-1}\cdot \sum_{q=0}^{n-1} \binom{n-1}{q}\cdot \frac{(-1)^q}{n(N-1)+q}\, .$$ Lemma \ref{x2} in case $x=n(N-1)$ evaluates the final sum as
$$(-1)^{nN+1}\binom{(r+N)n-1}{Nn-1}\, ,$$
which completes the proof of Theorem \ref{thm2} .\qed

\begin{lemma} \label{x1} Let $x_1, \ldots, x_N$ be fixed positive integers. Set $$\mathsf S(w_1, \ldots, w_N)=(1+w_1)^{-x_1}\cdot (w_1-w_2)^{-x_2}\cdots 
(w_1-w_N)^{-x_N}+\text { all symmetric combinations}.$$ 
Expand $\mathsf S(w_1, \ldots, w_N)$ in the region 
 $$w_1<<w_2<<\ldots<<w_N\, .$$ 
The free term of this expansion equals $$(-1)^{x_2+\ldots+x_N} \binom{x_1+\ldots+x_N-1}{x_1-1}\, .$$
\end{lemma} 

\proof We have 
\begin{multline*}(1+w_N)^{-x_1}\cdot (w_N-w_1)^{-x_2}\cdots (w_N-w_{N-1})^{-x_N}=\\
w_N^{-x_2-\ldots-x_N} \cdot (1+w_N)^{-x_1} \cdot \left(1-\frac{w_1}{w_N}\right)^{-x_2}\cdots \left(1-\frac{w_{N-1}}{w_N}\right)^{-x_N}\, .
\end{multline*} 
To extract the free term, we need the coefficient of 
$w_1^0\cdots w_{N-1}^0\cdot w_N^{x_2+\ldots+x_N}$ in 
$$(1+w_N)^{-x_1} \cdot \left(1-\frac{w_1}{w_N}\right)^{-x_2}\cdots
 \left(1-\frac{w_{N-1}}{w_N}\right)^{-x_N}\, .$$ 
This coefficient equals $$\binom{-x_1}{x_2+\ldots+x_N}=(-1)^{x_2+\ldots+x_N}\cdot \binom{x_1+\ldots+x_N-1}{x_1-1}\, .$$ An entirely parallel computation shows that the remaining terms 
$$(1+w_j)^{-x_1}\cdot (w_j-w_1)^{-x_2}\cdots (w_j-w_{N})^{-x_N}\,,$$ for $j\neq N$, do {\it not} contribute. 
\qed

\begin{lemma} \label{x2} For positive integers $x$ and $n$, we have 
$$\sum_{q=0}^{n-1} \frac{(-1)^q}{x+q}\binom{n-1}{q} =
\frac{(x-1)!(n-1)!}{(x+n-1)!}\, .$$
\end{lemma}
\proof We induct on $n$. For the inductive step, we compute \begin{eqnarray*} \sum_{q=0}^{n} \frac{(-1)^q}{x+q}\binom{n}{q} &=& \sum_{q=0}^{n} \frac{(-1)^q}{x+q}\binom{n-1}{q-1} +\sum_{q=0}^{n} \frac{(-1)^q}{x+q}\binom{n-1}{q}\\ &=& -\frac{x!(n-1)!}{(x+n)!}+\frac{(x-1)!(n-1)!}{(x+n-1)!}\\ & =& \frac{(x-1)!n!}{(x+n)!}
\, .
\end{eqnarray*} The first line is Pascal's identity, while the second line uses the induction hypothesis.
\qed

\subsection{Binomial identities} \label{binnum} We prove here part of Theorem \ref{cor1} stated in Section \ref{hnl1} 
together with Corollaries \ref{c2} and \ref{cor2}. \vskip.1in

\noindent {\it Proof of first half of Theorem \ref{cor1}.} The first statement in Theorem \ref{cor1} is purely combinatorial. In case $\text{rank }V=1$, the expression of Theorem \ref{thm2} 
simplifies: \begin{eqnarray*}\log \mathsf A_{1,1,N}(q)&=&\sum_{n=1}^{\infty} (-1)^{(N+1)n+1} \binom{(N+1)n-1}{Nn-1}\cdot \frac{q^n}{n}\\&=&N \cdot \sum_{n=1}^{\infty} (-1)^{Nn} \binom{-Nn-1}{n-1} \frac{q^n}{n}\, .\end{eqnarray*} The result
 can be rewritten in the form $$\mathsf A_{1,1,N}(q)=(1+t)^N \ \ \text{ for }\ \ 
 q= (-1)^N t(1+t)^{N}$$ using  Lemma \ref{identity} (ii) below. 

The main point of Theorem \ref{cor1}, however, is the formula \begin{equation}\label{expb}\mathsf B_{1,N}(q)=\frac{(1+t)^{N+1}}{1+(N+1)t} \ \ \ \ \text{ for }  \
 q=(-1)^N t(1+t)^N \end{equation} proven in Section \ref{later} below. 
The calculation uses a specialization to genus $0$ and is similar in spirit to
 the computations carried out in Section \ref{vis}.
 \qed
\vskip.1in

\noindent {\it Proof of Corollary \ref{c2}.} Similarly,  \begin{eqnarray*}\mathsf B_{1,N}(q)&=& \frac{(1+t)^{N+1}}{1+t(N+1)}\,\\&=& \sum_{n=0}^{\infty} (-1)^{nN} \cdot \binom{-Nn+N}{n}\cdot q^n\\
&=&\sum_{n=0}^{\infty} (-1)^{n(N+1)} \cdot \binom{(n-1)(N+1)}{n} \cdot q^n\end{eqnarray*} where we have used
Lemma \ref{identity}(i) with $d=0$ on the second line. \qed

\begin{lemma}\label{identity} For the change of variables $q=t(1+t)^r$, we have 
\begin{itemize}
\item [(i)] $$\sum_{n=0}^{\infty} \binom{d-rn+r}{n} \cdot q^n = \frac{(1+t)^{d+r+1}}{1+t(r+1)}\, ,$$
\item [(ii)] $$\log(1+t)=\sum_{n=1}^{\infty} \binom{-rn-1}{n-1}\cdot \frac{q^n}{n} \, .$$
\end{itemize}
\end{lemma}

\proof Part (i) is the content of \cite[Lemma 3]{MOP}. 
For part (ii), the identity to be established is \begin{equation}\label{eq3}\sum_{n=1}^{\infty} \frac{1}{n} \binom{-rn-1}{n-1}\cdot t^n (1+t)^{rn}= \log(1+t)\, .\end{equation} For the proof, 
we set $d=-2r-1$ in equation (i) $$\sum_{n=0}^{\infty} \binom{-r(n+1)-1}{n} \cdot t^n(1+t)^{rn}=\frac{(1+t)^{-r}}{1+t(r+1)}\, ,$$ which we  rewrite as $$\sum_{n=1}^{\infty} \binom{-rn-1}{n-1} \cdot t^{n-1}(1+t)^{rn-1}\cdot {(1+t(r+1))}=\frac{1}{1+t}\,. $$ The identity \eqref{eq3} is obtained by integration.\qed
\vskip.1in

\noindent {\it Proof of Corollary \ref{cor2}.} The first statement in  
Corollary \ref{cor2} follows by directly comparing Theorem \ref{cor1} and equation \eqref{sese}. Indeed, up to signs, the two universal functions $\mathsf A$
and  $\mathsf B$ agree for both sides. For the second statement, 
we observe  $$\int_{(\mathbb P^1)^{[n]}} s(L^{[n]})^{N}=\int_{\mathbb P^n} (1-h)^{N(\deg L-n+1)}=(-1)^n\binom{N(\deg L-n+1)}{n},$$ where Lemma \ref{segr} has been used in the first identity.  \qed

\section{Virtual invariants of surfaces: dimension 0 quotients} \label{vis}
\subsection{Overview} We prove here Theorems \ref{cor1}, \ref{thm4}, and \ref{rattt}. In particular, we study the virtual intersection theory of the 
Quot scheme $\mathsf {Quot}_{X}(\mathbb C^N, n)$ of short exact sequences 
$$0\to S\to \mathbb C^N\otimes \mathcal O_X\to Q\to 0,\,\,\,\ \  \chi(Q)=n\, ,\  c_1(Q)=0\, ,\ \text{rank}(Q)=0$$ on nonsingular projective surfaces $X.$
As noted in Section \ref{sd0},
 $\mathsf {Quot}_{X}(\mathbb C^N, n)$ carries a virtual fundamental class $$\left[\mathsf {Quot}_{X}(\mathbb C^N, n)\right]^{\mathrm{vir}}$$ of dimension $Nn.$ 
Our basic technique is to relate
integrals against the virtual class of Quot schemes of 
surfaces to integrals over Quot schemes of curves which
we have already studied.
Theorem \ref{thm4} is the first outcome.

The idea of dimensional reduction plays a central role in
the proof of Theorem \ref{rattt}. In the $N=1$ case, the
integrals over the Quot schemes of curves which arise  are
covered by Theorem \ref{thm1}. For higher $N$, a more
delicate analysis of the curve integrals is required.
A similar analysis is  used to 
complete the proof of Theorem \ref{cor1} in Section \ref{later}
(and appears also in the proof of 
Theorem \ref{genttt} for  surfaces of general type in 
Section \ref{jjj999}).

\subsection{Virtual integrals}
\subsubsection{Strategy} We first prove Theorem \ref{thm4}. The argument requires the following two steps: 
\begin{itemize}
\item [(i)] We show a universality statement allowing us to reduce to the 
case of a surface with nonsingular canonical curve $C\subset X$. 
\item [(ii)] The claim will then be obtained by direct comparison of the obstruction 
theories of the Quot schemes of $X$ and of $C$.  
\end{itemize} 

\subsubsection{Universality} \label{uniuni}
We will use equivariant localization to compute the series 
$$\mathsf Z_{X, N}(q, x_1, \ldots, x_\ell|\alpha_1, \ldots, \alpha_\ell)=\sum_{n=0}^{\infty} q^n \int_{\left[\mathsf {Quot}_{\,X}(\mathbb C^N, n)\right]^{\mathrm{vir}}} 
s_{x_1} (\alpha_1^{[n]}) \cdots s_{x_\ell} (\alpha_\ell^{[n]})\, .$$ 
The Quot scheme $\mathsf {Quot}_{X}(\mathbb C^N, n)$ 
carries torus action via the diagonal  $\mathbb C^{\star}$-action 
on the middle term of the sequence $$0\to S\to \mathbb C^N\otimes \mathcal O_X\to Q\to 0\, .$$ 
We write $w_1, \ldots, w_N$ for the equivariant weights. Just as in Section \ref{el}, 
the fixed loci are products of Hilbert schemes 
$$\mathsf F[n_1, \ldots, n_N]=X^{[n_1]}\times \cdots \times X^{[n_N]}$$ indexed by partitions $n_1+\ldots+n_{N}=n.$ We 
write $$S=\oplus_{i=1}^{N} I_{Z_i}\, ,\ \ \ Q=\oplus_{i=1}^{N} \mathcal O_{Z_i}
\, ,\ \ \  \text{ length }(Z_i)=n_i$$ for the fixed kernel and quotient. 
Furthermore, the induced obstruction theory of $\mathsf F[n_1, \ldots, n_N]$ splits: $$\text{Ext}^\bullet(S, Q)^{\text{fix}}=\oplus_{i=1}^{N} \text{Ext}^\bullet(I_{Z_i}, \mathcal O_{Z_i})\, .$$ In fact, the $\mathbb C^\star$-fixed
 obstruction sheaf is locally free with obstruction
 bundle
\begin{equation}\label{obbb} \left(K_X^{[n_1]}\oplus \ldots \oplus K_X^{[n_N]}\right)^{\vee}.\end{equation} This is a consequence of equation \eqref{obunn} below. 
The equivariant virtual normal bundle is the moving part of the tangent-obstruction theory
\begin{eqnarray*} \mathsf N[n_1, \ldots, n_N]^{\mathrm{vir}}=\text{Ext}^{\bullet}(S, Q)^{\text{mov}}&=& \bigoplus_{i\neq j} \text{Ext}^{\bullet} (I_{Z_i}, \mathcal O_{Z_j})[w_j-w_i]\, .\end{eqnarray*} 

Using the virtual localization theorem of \cite{GP}, the integral $$\int_{\left[\mathsf {Quot}_{\,X}(\mathbb C^N, n)\right]^{\mathrm{vir}}} s_{x_1} (\alpha_1^{[n]}) \cdots s_{x_\ell} (\alpha_\ell^{[n]})$$ 
can be rewritten as  
$$\sum_{n_1+\ldots+n_N=n} \int_{X^{[n_1]}\times \cdots \times X^{[n_N]}}\prod_{i=1}^{N} \mathsf e\left( (K_X^{[n_i]})^{\vee}\right) \,\cdot\, \prod_{i=1}^{\ell} s_{x_i} (\alpha_i^{[n_i]}[w_i]) \,\cdot \,\prod_{i\neq j} \mathsf e (\text{Ext}^{\bullet}(I_{Z_i}, \mathcal O_{Z_j})[w_j-w_i])^{-1}\, .$$ 
As in \cite [Theorem 5.1]{GNY}, we regard the above expression as a tautological integral over the Hilbert scheme of the disconnected surface $Y=X\sqcup X\sqcup \ldots \sqcup X,$ so that $$Y^{[n]}= \bigsqcup_{n_1+\ldots+n_N=n} X^{[n_1]}\times \cdots \times X^{[n_N]}.$$ The answer depends solely on the Chern numbers of the data involved: monomials in the Chern classes of $\alpha_i$ and Chern classes of the surface $X$. In the absence of better notation, we write $\mathsf m_k$ for these monomials enumerated in some order. Thus $$\mathsf Z_{X, N}(q, x_1, \ldots, x_\ell|\alpha_1, \ldots, \alpha_\ell)=\text{universal function of }{\mathsf m_k}.$$ Splitting the surface $X=X'\sqcup X''$ and the classes $\alpha_i=\alpha_i'\sqcup \alpha_i''$ one sees that $$\text{Quot}_X(\mathbb C^N, n)=\bigsqcup_{n'+n''=n} \text{Quot}_{X'}(\mathbb C^N, n')\times \text{Quot}_{X''}(\mathbb C^N, n''),$$ and the tangent-obstruction theory and the tautological elements $\alpha_i^{[n]}$ split as well. We then conclude the multiplicative form of the generating series $$\mathsf Z_{X, N}(q, x_1, \ldots, x_\ell|\alpha_1, \ldots, \alpha_\ell)=\prod \mathsf A_{k}^{\mathsf m_k}.$$ As usual, $\mathsf A_k$ are universal functions in the variables $q, x_1, \ldots, x_{\ell}$ that may depend on the ranks of the $\alpha$'s and $N$. 

To complete the proof of Theorem \ref{thm4}, 
we may assume  $X$ admits a nonsingular canonical curve 
$$C\subset X\, ,$$ 
since such surfaces $X$ separate all the monomials $\mathsf m_k$. 

\subsubsection {Quot schemes of curves and surfaces} 
For all nonsingular curves $C\subset X$, there is a 
natural embedding $$\iota: \mathsf{Quot}_{C} (\mathbb C^N, n)\hookrightarrow \mathsf{Quot}_{X} (\mathbb C^N, n)\, ,\ \ \ \  \left[\mathbb C^N\otimes \mathcal O_C\to Q\right]\mapsto \left[\mathbb C^N\otimes \mathcal O_X\to Q\right]\, .$$ In the case of canonical curves, the following result relates the obstruction theories of the Quot schemes above and plays a crucial role in
the proof of  Theorem \ref{thm4}. 

\begin{lemma}\label{l8} If $C$ is a nonsingular canonical curve, we have
 $$\iota_{\star} \left[\mathsf {Quot}_{C} (\mathbb C^N, n)\right]=(-1)^n  \left[\mathsf{Quot}_{X} (\mathbb C^N, n)\right]^{\mathrm{vir}}$$
 in the localized $\mathbb C^\star$-equivariant Chow theory of
$ \mathsf{Quot}_{X} (\mathbb C^N, n)$.
\end{lemma} 

\vskip.1in
\proof We first consider the case $N=1$. While the Hilbert scheme of points $X^{[n]}$ is smooth, the virtual fundamental class studied here does not equal the usual fundamental class. Indeed, $X^{[n]}$ carries the locally free obstruction sheaf $\left(K_X^{[n]}\right)^{\vee}$. The obstruction sheaf is obtained from the following sequence of
canonical isomorphisms: 
\begin{eqnarray}\label{obunn} 
\text{Ext}^1(I_{Z}, \mathcal O_{Z})&=&\text{Ext}^2(\mathcal O_Z, \mathcal O_Z)
\\ \nonumber &=&\text{Ext}^0(\mathcal O_Z, \mathcal O_Z\otimes K_X)^{\vee}\\
\nonumber &=&\text{Ext}^0(\mathcal O, \mathcal O_Z\otimes K_X)^{\vee}\\
\nonumber &=&\left(K_X^{[n]}\right)^{\vee}\Big|_{Z}\, .
\end{eqnarray} 
The defining equation $s$ of the canonical curve $C\subset X$ yields a section $s^{[n]}$ of $K_X^{[n]}$ via the assignment $$Z\mapsto s|_{Z}\in H^0(K_X\otimes \mathcal O_Z)\, .$$ 
The section $s^{[n]}$ vanishes precisely along\begin{equation}\label{emb} \iota: C^{[n]}\hookrightarrow X^{[n]}\, .\end{equation} 
Using that $X^{[n]}$ is smooth, we find \begin{equation}\label{emb1}\left[X^{[n]}\right]^{\mathrm{vir}}=\mathsf e\left((K_X^{[n]})^{\vee}\right)\cap X^{[n]}=(-1)^n\,  \iota_{\star} \left[C^{[n]}\right]\, ,\end{equation} 
which completes the proof of Lemma \ref{l8} in case $N=1$. 

Now let $N$ be arbitrary. We apply $\mathbb C^\star$-equivariant localization to both Quot schemes over $X$ and $C$ using the same weights for the two torus actions. The fixed loci are $$\mathsf F_C[n_1, \ldots, n_N]=C^{[n_1]}\times \cdots \times C^{[n_N]}\, ,
\ \ \  \mathsf F_X[n_1, \ldots, n_N]=X^{[n_1]}\times \cdots \times X^{[n_N]}$$ respectively. Parallel to \eqref{emb}, there is a natural embedding $$\iota: \mathsf F_C[n_1, \ldots, n_N]\hookrightarrow \mathsf F_X[n_1, \ldots, n_N].$$ We noted in \eqref{obbb} that the obstruction bundle of $\mathsf F_X[n_1, \ldots, n_N]$ splits as $$\left((K_X)^{[n_1]}\oplus \ldots \oplus (K_X)^{[n_N]}\right)^{\vee}.$$
Using \eqref {emb1}, we find that \begin{equation}\label{cl} \iota_{\star} \left[\mathsf F_C[n_1, \ldots, n_N]\right]= (-1)^n \mathsf e\left(\left((K_X)^{[n_1]}\oplus \ldots \oplus (K_X)^{[n_N]}\right)^{\vee}\right) \cap \left[\mathsf F_X[n_1, \ldots, n_N]\right]\end{equation}  $$=(-1)^n \left[\mathsf F_X[n_1, \ldots, n_N]\right]^{\mathrm{vir}}.$$ We furthermore claim  \begin{equation}\label{no} \iota^{\star} \mathsf e(\mathsf N_X[n_1, \ldots, n_N]^{\mathrm{vir}})=\mathsf e(\mathsf N_C[n_1, \ldots, n_N])\end{equation} where $\mathsf N_X^{\mathrm{vir}}$ and $\mathsf N_C$ are two normal bundles of the fixed loci.

The proof of \eqref{no} requires several steps.
 First, the difference $$\iota^{\star} \mathsf N_X[n_1, \ldots, n_N]^{\mathrm{vir}}-\mathsf N_C[n_1, \ldots, n_N]$$ equals \begin{eqnarray*} \bigoplus_{i\neq j} \text{Ext}^{\bullet}_X(I_{Z_i/X}, \mathcal O_{Z_j})[w_j-w_i] -\bigoplus_{i\neq j} \text{Ext}^{\bullet}_C(I_{Z_i/C}, \mathcal O_{Z_j})[w_j-w_i]\, .\end{eqnarray*} The latter expression can be further simplified using
\begin{eqnarray*} \text{Ext}_X^{\bullet} (I_{Z_i/X}, \mathcal O_{Z_j})-\text{Ext}^{\bullet}_C(I_{Z_i/C}, \mathcal O_{Z_j})&=&
-\text{Ext}_X^{\bullet} (\mathcal O_{Z_i}, \mathcal O_{Z_j}) +\text{Ext}_C^{\bullet} (\mathcal O_{Z_i}, \mathcal O_{Z_j})\\&=&-\text{Ext}^{\bullet}_C( \mathcal O_{Z_i}, \mathcal O_{Z_j}\otimes \Theta)[-1]\, ,\end{eqnarray*} where $\Theta=\mathcal O_C(C)$ is the theta characteristic of $C$. For the first equality, we have
expressed the ideal sheaves in terms of structure sheaves in $K$-theory.
The second equality follows from the exact sequence $$\ldots \to \text{Ext}_C^i(\mathcal O_{Z_i}, \mathcal O_{Z_j})\to \text{Ext}_X^i(\mathcal O_{Z_i}, \mathcal O_{Z_j})\to \text{Ext}_C^{i-1}(\mathcal O_{Z_i}, \mathcal O_{Z_j}\otimes \Theta)\to \ldots $$ proven, for instance, in \cite[Lemma $3.42$]{T}. Next, in the difference of the normal bundles, 
we group the terms corresponding to the pairs $(i, j)$ and $(j, i)$. 
We define $$\mathbb V_{ij}=\text{Ext}^{\bullet}_C( \mathcal O_{Z_i}, \mathcal O_{Z_j}\otimes \Theta)[w_j-w_i] \oplus \text{Ext}^{\bullet}_C( \mathcal O_{Z_j}, \mathcal O_{Z_i}\otimes \Theta)[w_i-w_j]\, ,$$ and write $$\iota^{\star} \mathsf N_X[n_1, \ldots, n_N]^{\mathrm{vir}}-\mathsf N_C[n_1, \ldots, n_N]=\bigoplus_{i<j} \mathbb V_{ij}\, .$$ By Serre duality, making use of the fact that $\Theta$ is a theta characteristic, we obtain $$\mathbb V_{ij}^{\vee}=\mathbb V_{ij}[1]\, .$$ Therefore, $$\mathsf e_{\mathbb C^{\star}}(\mathbb V_{ij})=(-1)^{\chi(\mathcal O_{Z_i}, \mathcal O_{Z_j}\otimes \Theta)}=1\, ,$$ which proves  \eqref{no}. 

Finally, by the virtual localization formula \cite{GP}, we have $$\left[\mathsf{Quot}_{X} (\mathbb C^N, n)\right]^{\mathrm{vir}}=\sum_{n_1+\ldots+n_N=n} (j_X)_{\star} \left( \frac{1}{\mathsf e(\mathsf N_X[n_1, \ldots, n_N]^{\mathrm{vir}})} \cap \left[\mathsf F_X[n_1, \ldots, n_N]\right]^{\mathrm{vir}}\right)$$ $$\left[\mathsf{Quot}_{C} (\mathbb C^N, n)\right]=\sum_{n_1+\ldots+n_N=n} (j_C)_{\star} \left( \frac{1}{\mathsf e(\mathsf N_C[n_1, \ldots, n_N])}\cap \left[\mathsf F_C[n_1, \ldots, n_N]\right]\right).$$ Using equations \eqref{cl} and \eqref{no} we obtain $$\iota_{\star} \left[\mathsf {Quot}_{C} (\mathbb C^N, n)\right]=(-1)^n  \left[\mathsf{Quot}_{X} (\mathbb C^N, n)\right]^{\mathrm{vir}}\, ,$$  which proves the Lemma. 
\qed

\begin{remark} The result of Lemma \ref{l8} should be expected. In fact, the canonical curve $C$ gives a cosection $$\text{Ob}\to \mathcal O_{\mathsf{Quot}}$$ of the obstruction sheaf of $\mathsf{Quot}_{X} (\mathbb C^N, n)$ via the composition $$\text{Ext}^1(S, Q)\to \text{Ext}^2(Q, Q)\stackrel{\text{Trace}}{\to} H^2(\mathcal O_X)=
H^0(K_X)^{\vee}\to \mathbb C\, .$$ A careful analysis shows that the cosection vanishes along the quotients supported on $C$. By \cite {KL}, the virtual fundamental cycle is localized along such quotients. 
However, the precise determination of the cycle still requires a calculation. The known techniques require stronger smoothness assumptions than what we can prove in our case, so we have given a
different argument for the proof of Lemma \ref{l8}. 

For example, $\mathsf{Quot}_{X} (\mathbb C^N, n)$ is singular for every
 $N\geq 2$ and $n\geq 2$ even at quotients of the form $$Q=\mathcal O_Z\oplus \mathcal O_Z\,, \ \ \  \text{length}(Z)=\frac{n}{2}\, .$$ 
Indeed, the Zariski tangent space $$\text{Hom}(S, Q)=\text{Hom}(I_Z\,\oplus\, I_Z\,\oplus\, \mathbb C^{N-2}\,\otimes\, \mathcal O_X, \mathcal O_Z\,\oplus \,\mathcal O_Z)$$ has dimension $(N+2)n$ which is 
higher than the actual dimension $(N+1)n$. \end{remark}

\subsubsection{Proof of Theorem \ref{thm4}} We argued in Section \ref{uniuni} that it suffices to consider the case when $X$ admits a nonsingular canonical curve $C$. Let $\alpha_i$ be classes on $X$ and set $\beta_i=\alpha_i|_C$. By Lemma \ref{l8}, we have $$ \int_{\left[\mathsf {Quot}_{\,X}(\mathbb C^N, n)\right]^{\mathrm{vir}}} s_{x_1} (\alpha_1^{[n]}) 
 \cdots s_{x_\ell} (\alpha_\ell^{[n]})= (-1)^n \int_{\mathsf {Quot}_{\,C}(\mathbb C^N, n)} s_{x_1} (\beta_1^{[n]}) \cdots s_{x_\ell} (\beta_\ell^{[n]})\, .$$ Theorem \ref{thm4} follows immediately \begin{eqnarray*}\mathsf Z_{X, N}(q, x_1, \ldots, x_{\ell}\,|\,\alpha_1, \ldots, \alpha_{\ell})&=& \mathsf Z_{g, N} (-q, x_1, \ldots, x_{\ell}, \beta_1, \ldots, \beta_{\ell})\\&=&\mathsf A_1(-q)^{c_1(\alpha_i)\cdot K_X}  \cdots \mathsf A_{\ell}(-q)^{c_1(\alpha_{\ell})\cdot K_X}\cdot \mathsf B(-q)^{1-g}\,.\end{eqnarray*}
 \qed

\subsection{Virtual Euler characteristics}\label{ves} 
Theorem \ref{rattt} will be proven next. Before presenting the argument, we review general statements regarding virtual Euler characteristics. 

\subsubsection {Generalities} Let $Z$ be a scheme admitting a 
2-term perfect obstruction theory $$\mathbb E^{\bullet}=\left[E_{-1}\to E_0\right]\to \tau^{[-1, 0]}\mathbb L_{Z},$$ and a virtual fundamental class $\left[Z\right]^{\mathrm{vir}}$ of dimension $$d=\text{rank }E_0-\text{rank }E_{-1}\, .$$ 
The virtual tangent bundle $T^{\mathrm{vir}}Z$ is defined in the $K$-theory of $Z$ as the difference $$(E_0)^{\vee}-(E_{-1})^{\vee}\, .$$ 
We define the virtual Euler characteristic 
\begin{equation}\label{ff99}
\mathsf e^{\mathrm{vir}} (Z)=\int_{\left[Z\right]^{\mathrm{vir}}} c_d (T^{\mathrm{vir}}Z)\, ,
\end{equation} see also \cite{GF}. 
Virtual Euler characteristics are deformation invariants.

In particular, if $Z$ is nonsingular 
with a locally free obstruction bundle $B$, 
then  $$\left[Z\right]^{\mathrm{vir}}=\mathsf e(B)\cap \left[Z\right]\, $$ 
and the virtual tangent bundle is the difference $TZ-B$. 
By definition, we obtain \begin{equation}\label{smooth}\mathsf e^{\mathrm{vir}} (Z)=\int_{Z} \mathsf e(B)\cdot \frac{c(TZ)}{c(B)}\, .\end{equation}

\subsubsection{Proof of Theorem \ref{rattt} for $N=1$} \label{jjj222}
We must prove
\begin{equation}\label{vireu}
\sum_{n=0}^{\infty} q^{n} \cdot \mathsf e^{\mathrm{vir}}\left(X^{[n]}\right)=
\left(\frac{(1-q)^2}{1-2q}\right)^{K_X^2}\, .\end{equation}

\proof We observed in Lemma \ref{l8} that the Hilbert schemes $X^{[n]}$ have locally free obstruction sheaves $\left(K_X^{[n]}\right)^{\vee}$. By \eqref{smooth}, the virtual Euler characteristics are $$\mathsf e^{\mathrm{vir}} (X^{[n]})=\int_{X^{[n]}} \mathsf e\left((K_X^{[n]})^{\vee}\right)\cdot \frac{c(TX^{[n]})}{c \left((K_X^{[n]})^{\vee}\right)}\, .$$ 
The above rewriting of the virtual Euler characteristic shows, via \cite[Theorem 4.5] {EGL}, that expression \eqref{vireu} takes the universal form $$\mathsf U(q)^{K_X^2}\cdot \mathsf V(q)^{c_2(X)}\, .$$ 
To prove $$\mathsf U(q)= (1-q)^{2}\cdot (1-2q)^{-1}\, , \ \ \ \
 \mathsf V(q)=1\, ,$$ we may specialize to surfaces $X$ which admit a 
nonsingular canonical curve $$C\subset X\, .$$ 
By \eqref{emb}, we have the embedding   
$$\iota: C^{[n]}\hookrightarrow X^{[n]}\, $$ and furthermore, by \eqref{emb1},
we have
  $$\left[X^{[n]}\right]^{\mathrm{vir}}=\mathsf e\left((K_X^{[n]})^{\vee}\right)\cap X^{[n]}=(-1)^n\,  \iota_{\star} \left[C^{[n]}\right]\, .$$
We conclude 
 $$\mathsf e^{\mathrm{vir}} (X^{[n]})=(-1)^n \int_{C^{[n]}} \iota^{\star}  \frac{c(TX^{[n]})}{c \left((K_X^{[n]})^{\vee}\right)}\, .$$ 
Going further, let $\Theta=\mathcal O_C(C)$ be the theta characteristic of $C$. If $Z\subset C$, consider the exact sequence $$0\to \mathcal O_X(-C)\to I_{Z/X}\to \iota_{\star} I_{Z/C}\to 0.$$ Taking $\text{Hom} (, \mathcal O_Z)$ we find $$0\to TC^{[n]}\to \iota^{\star} TX^{[n]} \to \Theta^{[n]}\to 0\implies \iota^{\star} c(TX^{[n]}) = c(\Theta^{[n]}) \cdot c(TC^{[n]})\, .$$ 
Moreover, we have $$\iota^{\star} K_X^{[n]}=\Theta^{[n]}\, .$$
We conclude
\begin{equation}\label{qpp22}
\mathsf e^{\mathrm{vir}} (X^{[n]})=(-1)^n \int_{C^{[n]}} \frac{c(\Theta^{[n]}) \cdot c(TC^{[n]})}{c \left((\Theta^{[n]})^{\vee}\right)}\, .
\end{equation}

There are now several ways to evaluate the integral
\eqref{qpp22}, but the most direct path is to use Theorem \ref{thm1}. 
We observe $$TC^{[n]}= \left(K_C^{[n]}\right)^{\vee}\, .$$ Then, we have 
\begin{eqnarray*}
\mathsf e^{\mathrm{vir}} (X^{[n]}) &=&(-1)^n \int_{C^{[n]}} \frac{c(\Theta^{[n]}) \cdot c\left((K_C^{[n]})^{\vee}\right)}{c \left((\Theta^{[n]})^{\vee}\right)}\\
&= &(-1)^n\int_{C^{[n]}} s_1((-\Theta)^{[n]})\cdot s_{-1} ((-K_C)^{[n]}) \cdot s_{-1} (\Theta^{[n]})\, .
\end{eqnarray*}
Invoking Theorem \ref{thm1}, we find
 $$\sum_{n=0}^{\infty} q^{n} \cdot \mathsf e^{\mathrm{vir}}\left(X^{[n]}\right)=\mathsf Z_{C,1} (-q, x_1=1, x_2=-1, x_3=-1 \,|\, \alpha_1=-\Theta, \alpha_2=-K_C, \alpha_3=\Theta)$$ $$=\mathsf A_1^{\deg \alpha_1}\cdot 
\mathsf A_2^{\deg \alpha_2}\cdot \mathsf A_3^{\deg \alpha_3}\cdot \mathsf B^{1-g}\, .$$ 
The change of variables specified by
Theorem \ref{thm1} takes the simple form $$-q=\frac{t}{1-t}\, ,$$ and
 the universal functions are $$\mathsf A_1= 1-t=(1-q)^{-1}\, ,\ \ \ 
 \mathsf A_2= \mathsf A_3=1+t=(1-2q)(1-q)^{-1}\, , \ \ \ 
\mathsf B=1\, .$$ 
We conclude $$\sum_{n=0}^{\infty} q^{n} \cdot \mathsf e^{\mathrm{vir}}\left(X^{[n]}\right)=\left((1-q)^2\cdot (1-2q)^{-1}\right)^{K_X^2}\, ,$$ which
completes the proof of the $N=1$ case of Theorem \ref{rattt}. \qed

\begin{remark} Using the same techniques, we can also compute the virtual $\chi_{-y}$ genera: $$\sum_{n=0}^{\infty} q^n \cdot \chi^{\mathrm{vir}}_{-y}(X^{[n]}) = \left(\frac{(1-q)\cdot (1-yq)}{1-q-qy}\right)^{K_X^2}\, .$$
Theorem \ref{rattt} is then recovered in the limit $y\to 1$. 
\end{remark} 

\begin{remark} For future reference, we record the following slight generalization of the above calculations. For any nonsingular projective surface $X$ and $M\to X$ a line bundle, set $$\mathsf Z_{X, M}=\sum_{n=0}^{\infty} q^n\int_{X^{[n]}} \mathsf e\left(\left(M^{[n]}\right)^{\vee}\right) \frac{c(TX^{[n]})}{c\left(\left(M^{[n]}\right)^{\vee}\right)}\, .$$ Without the duals placed on tautological bundles, such integrals also appear in the work \cite {K} on stable pair invariants of local surfaces. The above calculations yield the following result.

\begin{corollary}\label{pro2} We have
 \begin{equation}\label{xm}\mathsf Z_{X, M} = \mathsf U(q)^{c_1(M)^2} \cdot \mathsf V(q)^{c_1(M)\cdot K_X}\end{equation} where $$\mathsf U(q) = 1-q\, ,\ \ \ \mathsf V(q) = (1-2q)^{-1}\cdot (1-q)\, .$$\end{corollary}
\end{remark}

\subsubsection{Proof of Theorem \ref{rattt} for higher $N$} \label{prt5} Theorem \ref{rattt} concerns the generating series
\begin{equation}\label{jj23}
\mathsf{Z}^{\mathcal{E}}_{X,N,0}=\sum_{n=0}^\infty q^n \mathsf e^{\mathrm{vir}}(\text{Quot}_X(\mathbb C^N, n))\, .
\end{equation}
For notational convenience, we will denote the series \eqref{jj23}
by $\mathsf E_{X}(q)$.
We will follow a strategy similar to that of the proof of Theorem \ref{thm4}:
\begin{itemize}
\item [(i)] We will
first show the factorization $$\mathsf E_{X}(q)=\mathsf A (q)^{K_X^2}\cdot \mathsf B (q)^{\chi(\mathcal O_X)}$$ holds for universal power series $\mathsf A, \mathsf B\in \mathbb Q[[q]]$.
\item [(ii)] To identify the series $\mathsf A, \mathsf B$, we will
use Theorem \ref{thm4} to localize the calculation to a nonsingular
 canonical curve $$C\subset X\, .$$ 
\item [(iii)] The evaluation $\mathsf B=1$ will follow
for formal reasons. 
\item[(iv)]
To determine $\mathsf A$, we will 
use equivariant localization on $\mathsf {Quot}_C(\mathbb C^N, n)$ for
$C=\mathbb{P}^1$.
We will find closed form expressions for the localization sums which will furthermore prove the rationality
of Theorem \ref{rattt}.
\end{itemize}

\begin{remark} \label{cc445} We warn the reader that both the statement
and the proof of the
 torus equivariant localization formula for virtual Euler characteristics stated  in
 \cite[Corollary 6.6 (3)]{GF} are wrong. In particular, application
of \cite[Corollary 6.6 (3)]{GF} to the diagonal $\mathbb C^{\star}$-action
on $\mathbb{C}^N$ to calculate
$\mathsf{Z}^{\mathcal{E}}_{X,N,0}$ in terms of $\mathsf{Z}^{\mathcal{E}}_{X,1,0}$ 
will give incorrect results.{\footnote{B. Fantechi and L. G\"ottsche
agree with Remark \ref{cc445} about the error in part (3), but
confirm that parts (1) and (2) of \cite[Corollary 6.6]{GF} are
correct.}}

\end{remark} 

\vskip.1in
\noindent {\it Step (i).} We first apply the virtual localization formula
 to prove that the series $\mathsf E_{X}(q)$ depends only upon $K_X^2$ and $\chi(\mathcal O_X)$. By definition, 
$$\mathsf e^{\mathrm{vir}}(\text{Quot}_X(\mathbb C^N, n))=\int_{\left[\text{Quot}_X(\mathbb C^N, n)\right]^{\mathrm{vir}}} c(T^{\mathrm{vir}} \text{Quot}_X(\mathbb C^N, n))$$ where $$T^{\mathrm{vir}}\text{Quot}_X(\mathbb C^N, n)=\text{Ext}^0(S, Q)-\text{Ext}^1(S, Q)$$ is the virtual tangent bundle. By the virtual
localization formula of \cite{GP}, we obtain $$\mathsf e^{\mathrm{vir}}(\text{Quot}_X(\mathbb C^N, n))=\sum_{n_1+\ldots+n_N=n} \int_{\left[X^{[n_1]}\times \cdots \times X^{[n_N]}\right]^{\mathrm{vir}}} \frac{\iota^{\star} c(T^{\mathrm{vir}} \text{Quot}_X(\mathbb C^N, n))}{\mathsf e_{\mathbb C^{\star}} (\mathsf N[n_1, \ldots, n_N]^{\mathrm{vir}})}\, .$$ 
Using  $$\iota^{\star} T^{\mathrm{vir}} \text{Quot}_X(\mathbb C^N, n)=\bigoplus_{i, j} \text{Ext}^{\bullet}(I_{Z_i}, \mathcal O_{Z_j})[w_j-w_i]$$ and $$\mathsf N[n_1, \ldots, n_N]^{\mathrm{vir}}=\bigoplus_{i\neq j} \text{Ext}^{\bullet}(I_{Z_i}, \mathcal O_{Z_j})[w_j-w_i]\, ,$$ we rewrite the right hand side of the
virtual localization as
$$\sum_{n_1+\ldots+n_N=n} \int_{X^{[n_1]}\times \cdots \times X^{[n_N]}} \prod_{i=1}^{N}\mathsf e\left((K_X^{[n_i]})^{\vee}\right) \,\cdot \,c(\text{Ext}^{\bullet}(I_{Z_i}, \mathcal O_{Z_i})) \,\cdot\, \prod_{i\neq j} \frac{c(\text{Ext}^{\bullet}(I_{Z_i}, \mathcal O_{Z_j})[w_j-w_i])}{\mathsf e(\text{Ext}^{\bullet}(I_{Z_i}, \mathcal O_{Z_j})[w_j-w_i])}\, .$$
As in \cite [Theorem 5.1]{GNY}, each Hilbert scheme integral depends solely on the Chern numbers of the surface $X$, so $\mathsf E_X(q)$ is a function of $$K_X^2 \text { and } \chi(\mathcal O_X)\, .$$ 
By splitting the surface $X=X'\sqcup X''$, we see
 $$\text{Quot}_X(\mathbb C^N, n)=\bigsqcup_{n'+n''=n} \text{Quot}_{X'}(\mathbb C^N, n')\times \text{Quot}_{X''}(\mathbb C^N, n'') $$ with
a splitting also of  the obstruction theory. 
We therefore conclude 
$$\mathsf E_X(q)=\mathsf E_{X'}(q)\cdot \mathsf E_{X''}(q)\, ,$$ which implies the factorization $$\mathsf E_X(q)=\mathsf A(q)^{K_X^2} \cdot \mathsf B(q)^{\chi(\mathcal O_X)}\, .$$ 

\vskip.1in
\noindent {\it Step (ii).} When $C\subset X$ is a nonsingular canonical curve, we can apply the result of Lemma \ref{l8} to write $$\iota_{\star} \left[\mathsf {Quot}_{C} (\mathbb C^N, n)\right]=(-1)^n  \left[\mathsf{Quot}_{X} (\mathbb C^N, n)\right]^{\mathrm{vir}}.$$ Here $$\iota: \mathsf {Quot}_{C} (\mathbb C^N, n)\to \mathsf{Quot}_{X} (\mathbb C^N, n)$$ is the natural inclusion $$\left[\mathbb C^N\otimes \mathcal O_C\to Q\right]\mapsto \left[\mathbb C^N\otimes \mathcal O_X\to Q\right].$$ As a consequence, we obtain \begin{eqnarray*}\mathsf e^{\mathrm{vir}}(\mathsf{Quot}_{X} (\mathbb C^N, n))&=&\int_{\left[\text{Quot}_X(\mathbb C^N, n)\right]^{\mathrm{vir}}} c(T^{\mathrm{vir}} \text{Quot}_X(\mathbb C^N, n))\\ &=& (-1)^n \int_{\text{Quot}_C(\mathbb C^N, n)} c(\iota^{\star} T^{\mathrm{vir}} \text{Quot}_X(\mathbb C^N, n))\\&=& (-1)^n \int_{\text{Quot}_C(\mathbb C^N, n)} c(T \text{Quot}_C(\mathbb C^N, n))\cdot c(\mathcal T_n).
\end{eqnarray*}
Here, $\mathcal T_n\to \text{Quot}_C(\mathbb C^N, n)$ is the virtual bundle given pointwise by $$\mathcal T_n=\text{Ext}^{\bullet}_C(Q, Q\otimes \Theta)\, ,$$ 
where $\Theta=N_{C/X}$ is the associated theta characteristic. The last line follows from the $K$-theoretic decomposition \begin{equation}\label{splt}\iota^{\star} T^{\mathrm{vir}} \text{Quot}_X(\mathbb C^N, n)=T \text{Quot}_C(\mathbb C^N, n)+
\mathcal T_n\, .\end{equation} 

To prove \eqref{splt}, let $S_C$ denote the kernel of the surjection $$\mathbb C^N\otimes \mathcal O_C\to Q\to 0$$ on the curve $C$, 
and let $S$ denote the kernel of the similar surjection $$\mathbb C^N\otimes \mathcal O_X\to Q\to 0$$ on the surface $X$. 
The splitting \eqref{splt} is a consequence of the following computation: \begin{eqnarray*} \text{Ext}_X^{\bullet}(S, Q)-\text{Ext}_C^{\bullet}(S_C, Q)&=& -\text{Ext}_X^{\bullet}(Q, Q)+ \text{Ext}_C^{\bullet}(Q, Q)\\ &=& 
-\text{Ext}_C^{\bullet}(Q, Q\otimes \Theta)[-1]\, .\end{eqnarray*}
For the first equality, we have
expressed $S, S_C$ in terms of $Q$ in the $K$-theory of $X$ and $C$. 
The second equality follows from the exact sequence $$\ldots \to \text{Ext}_C^i(Q, Q)\to \text{Ext}_X^i(Q, Q)\to \text{Ext}_C^{i-1}(Q, Q\otimes \Theta)\to \ldots $$ provided by \cite[Lemma $3.42$]{T}.

\vskip.1in
\noindent {\it Step (iii).} By (ii), we are now left to evaluating the generating series $$\mathsf E_C(q)= \sum q^n (-1)^n \cdot \int_{\text{Quot}_C(\mathbb C^N, n)} c(T \text{Quot}_C(\mathbb C^N, n))\cdot c(\mathcal T_n).$$ By the argument in Step (i), the answer takes the form $$\mathsf E_C(q)=\mathsf A(q)^{1-g}$$ with $g$ the genus of $C$. The second series $\mathsf B(q)=1$ since there is no $\chi(\mathcal O_X)$-dependence in the curve integral above.

\vskip.1in
\noindent {\it Step (iv).}
To determine the series $\mathsf A$, we specialize first to the $N=2$ case. We prove $$\mathsf A(q)=\frac{(1-4q)^2}{(1-q)^2\cdot (1-6q+q^2)}\, .$$ 

The problem at hand is now purely 
a curve calculation. We can therefore  discard the surface $X$ and concentrate
on the curve $C$. To find $\mathsf A$, we take 
$$C=\mathbb P^1\, .$$
 Our goal is then to prove the second equality in
the equation
 \begin{eqnarray}\mathsf A(q) &=& 
\label{sumeq}\sum_{n=0}^{\infty} q^n (-1)^{n} \cdot \int_{\text{Quot}_{\mathbb P^1}(\mathbb C^2, n)} c(T \text{Quot}_{\mathbb P^1}(\mathbb C^2, n))\cdot c(\mathcal T_n)\\
\nonumber &=&\frac{(1-4q)^2}{(1-q)^2\cdot (1-6q+q^2)}\, .
\end{eqnarray}

We will apply $\mathbb{C}^\star$-equivariant localization on $\text{Quot}_{\mathbb P^1}(\mathbb C^2, n)$. We write $$\mathbb C^2=\mathbb C[w_1]\oplus \mathbb C[w_2]$$ for the weights of the diagonal $\mathbb C^{\star}$-action on $\mathbb C^2$. 
The fixed loci are $${\mathsf F} [n_1, n_2]=C^{[n_1]}\times C^{[n_2]}=\mathbb P^{n_1}\times \mathbb P^{n_2}\stackrel{\iota}{\hookrightarrow} \text{Quot}_{\mathbb P^1}(\mathbb C^2, n)\, .$$ The fixed points correspond to the exact sequences $$0\to I_{Z_1}\oplus I_{Z_2}\to \mathbb C^2\otimes \mathcal O_{\mathbb P^1}\to \mathcal O_{Z_1}\oplus \mathcal O_{Z_2}\to 0\, .$$
Thus, by Atiyah-Bott localization, we find  \begin{equation}\label{loceqsum}\int_{\text{Quot}_{\mathbb P^1}(\mathbb C^2, n)} c(T \text{Quot}_{\mathbb P^1}(\mathbb C^2, n))\cdot c(\mathcal T_n)=\sum_{n_1+n_2=n} \int_{\mathbb P^{n_1}\times \mathbb P^{n_2}} \mathsf {Contr}(n_1, n_2)\, .\end{equation} Here, we set $$\mathsf {Contr}(n_1, n_2)=\frac{c(\iota^{\star} T\text{Quot}_{\mathbb P^1}(\mathbb C^2, n))\cdot c(\iota^{\star}\mathcal T_n)}{\mathsf e_{\mathbb C^{\star}}(\mathsf N[n_1, n_2])}$$ for the contribution of the $(n_1, n_2)$-fixed locus, where  $\mathsf N[n_1, n_2]$ denotes the normal bundle. We will evaluate \eqref{loceqsum} explicitly. 

For the analysis of $\mathsf {Contr}(n_1, n_2)$, the notation
$w=w_2-w_1$
will be convenient. We compute \begin{eqnarray*}\iota^{\star} T\text{Quot}_C(\mathbb C^2, n)=T\mathbb P^{n_1}+T\mathbb P^{n_2} + \text{Ext}^{\bullet}(I_{Z_1}, \mathcal O_{Z_2})[w]+\text{Ext}^{\bullet}(I_{Z_2}, \mathcal O_{Z_1})[-w]\, . \end{eqnarray*} The last two terms come from the normal bundle $$\mathsf N[n_1, n_2]=\text{Ext}^{\bullet}(I_{Z_1}, \mathcal O_{Z_2})[w]+\text{Ext}^{\bullet}(I_{Z_2}, \mathcal O_{Z_1})[-w]\, .$$ Similarly, $\iota^{\star} \mathcal T_n$ can be written
 as $$\text{Ext}^{\bullet}(\mathcal O_{Z_1}, \mathcal O_{Z_1}(-1))+\text{Ext}^{\bullet}(\mathcal O_{Z_2}, \mathcal O_{Z_2}(-1))$$ $$+\text{Ext}^{\bullet}(\mathcal O_{Z_1}, \mathcal O_{Z_2}(-1))[w]+\text{Ext}^{\bullet}(\mathcal O_{Z_2}, \mathcal O_{Z_1}(-1))[-w].$$ 

We now explicitly compute the various tautological structures appearing above. The arguments follow the proof of \cite[Theorem 2]{MOP}. 
We observe that the universal subschemes $$\mathcal Z_1\subset \mathbb P^1\times \mathbb P^{n_1}\, ,\ \ \  \mathcal Z_2\subset \mathbb P^{1}\times \mathbb P^{n_2}$$ take the form $$\mathcal O(-\mathcal Z_1)= \mathcal O_{\mathbb P^1}(-n_1)\boxtimes \mathcal O_{\mathbb P^{n_1}}(-1)\,, \ \ \  \mathcal O(-\mathcal Z_2)= \mathcal O_{\mathbb P^1}(-n_2)\boxtimes \mathcal O_{\mathbb P^{n_2}}(-1)\, .$$ 
We require the following three calculations:
\begin{multline*}
\hspace{-10pt}
\text{Ext}^{\bullet}(\mathcal O_{\mathcal Z_1}, \mathcal O_{\mathcal Z_1}(-1))=
\text{Ext}^{\bullet} (\mathcal O-\mathcal O(-\mathcal Z_1), \mathcal O_{\mathbb P^1}(-1)-\mathcal O(-\mathcal Z_1)\otimes \mathcal O_{\mathbb P^1}(-1))\hspace{50pt} \\=
\text{Ext}^{\bullet} (\mathcal O-\mathcal O_{\mathbb P^1}(-n_1)\boxtimes \mathcal O_{\mathbb P^{n_1}}(-1), \mathcal O_{\mathbb P^1}(-1)-\mathcal O_{\mathbb P^1}(-n_1-1)\boxtimes \mathcal O_{\mathbb P^{n_1}}(-1))\\=
\mathbb C^{n_1} \otimes \mathcal O_{\mathbb P^{n_1}}(-1)-\mathbb C^{n_1} \otimes \mathcal O_{\mathbb P^{n_1}}(1)\, , \hspace{200pt}\end{multline*}
\vspace{-5pt}
$$\text{Ext}^{\bullet} (I_{\mathcal Z_1}, \mathcal O_{\mathcal Z_2}) = \,\mathbb C^{n_1+1}\, \otimes\, \mathcal O_{\mathbb P^{n_1}}(1)\,-\,\mathbb C^{n_1-n_2+1}\,\otimes\, \mathcal O_{\mathbb P^{n_1}}(1)\,\otimes \,\mathcal O_{\mathbb P^{n_2}}(-1)\, , \hspace{100pt}$$ 
\vspace{-5pt}
 $$\text{Ext}^{\bullet} (\mathcal O_{\mathcal Z_1}, \mathcal O_{\mathcal Z_2}(-1)) =- \,\mathbb C^{n_1}\, \otimes\, \mathcal O_{\mathbb P^{n_1}}(1)\,+\,\mathbb C^{n_2}\,\otimes\, \mathcal O_{\mathbb P^{n_2}}(-1)\,+\,\mathbb C^{n_1-n_2}\,\otimes\, \mathcal O_{\mathbb P^{n_1}}(1)\,\otimes \,\mathcal O_{\mathbb P^{n_2}}(-1)\, .$$

\vspace{5pt}
\noindent As a consequence, we find $$\iota^{\star} T\text{Quot}_C(\mathbb C^2, n)+ \iota^{\star} \mathcal T_n$$ can be calculated as $$T\mathbb P^{n_1}+T\mathbb P^{n_2}+\left(\mathbb C^{n_1} \otimes \mathcal O_{\mathbb P^{n_1}}(-1)-\mathbb C^{n_1} \otimes \mathcal O_{\mathbb P^{n_1}}(1)\right)+\left(\mathbb C^{n_2} \otimes \mathcal O_{\mathbb P^{n_2}}(-1)-\mathbb C^{n_2} \otimes \mathcal O_{\mathbb P^{n_2}}(1)\right)$$ $$+\left(\mathcal O_{\mathbb P^{n_1}}(1)+\mathbb C^{n_2}\otimes \mathcal O_{\mathbb P^{n_2}}(-1) - \mathcal O_{\mathbb P^{n_1}}(1)\,\otimes \,\mathcal O_{\mathbb P^{n_2}}(-1)\right)[w]$$ $$+\left(\mathcal O_{\mathbb P^{n_2}}(1)+\mathbb C^{n_1}\otimes \mathcal O_{\mathbb P^{n_1}}(-1) - \mathcal O_{\mathbb P^{n_1}}(-1)\,\otimes \,\mathcal O_{\mathbb P^{n_2}}(1)\right)[-w]\, .$$ We also have \begin{equation}\label{nn1n2}\mathsf N[n_1, n_2]=\left(\mathbb C^{n_1+1}\otimes \mathcal O_{\mathbb P^{n_1}}(1)  - \mathbb C^{n_1-n_2+1} \otimes \mathcal O_{\mathbb P^{n_1}}(1)\,\otimes \,\mathcal O_{\mathbb P^{n_2}}(-1)\right)[w]  $$ $$+ \left(\mathbb C^{n_2+1} \otimes \mathcal O_{\mathbb P^{n_2}}(1)- \mathbb C^{n_2-n_1+1} \otimes  \mathcal O_{\mathbb P^{n_1}}(-1)\,\otimes \,\mathcal O_{\mathbb P^{n_2}}(1)\right) [-w].\end{equation}

We write $h_1$ and $h_2$ for the hyperplane classes on $\mathbb P^{n_1}$
and $\mathbb P^{n_2}$ respectively. 
After substituting the last equation into \eqref{loceqsum},
 we find $$\int_{\text{Quot}_{\mathbb P^1}(\mathbb C^2, n)} c(T \text{Quot}_{\mathbb P^1}(\mathbb C^2, n))\cdot c(\mathcal T_n)=\sum_{n_1+n_2=n} \int_{\mathbb P^{n_1}\times \mathbb P^{n_2}} \mathsf {Contr}(n_1, n_2)\, ,$$ where $\mathsf {Contr}(n_1, n_2)$ is given by 
$$\frac{(1-h_1)^{n_1}(1+h_1)(1-h_2)^{n_2}(1+h_2)(1+h_1+w)(1-h_2+w)^{n_2}(1+h_2-w)(1-h_1-w)^{n_1}}{(1+h_1-h_2+w)(1-h_1+h_2-w)}$$ $$\cdot \frac{(w+h_1-h_2)^{n_1-n_2+1}(-w-h_1+h_2)^{n_2-n_1+1}}{(h_1+w)^{n_1+1}(h_2-w)^{n_2+1}}\, .$$ 

While the expression may seem unwieldy, nonetheless, we will be able to sum the localization contributions explicitly via the Lagrange-B\"urmann formula.
We  write \begin{equation}\label{phipsi}\Phi_1(h_1)=(1-h_1)\cdot (1-h_1-w)\cdot (h_1+w)^{-1}\end{equation} $$\Phi_2(h_2)=(1-h_2)\cdot(1-h_2+w)\cdot (h_2-w)^{-1}$$ $$\Psi(h_1, h_2)=(1+h_1)\cdot (1+h_2)\cdot (1+h_1+w)\cdot (1+h_2-w)\cdot (1+h_1-h_2+w)^{-1}\cdot (1-h_1+h_2-w)^{-1}$$ $$(h_1+w)^{-1}\cdot (h_2-w)^{-1}\cdot (w+h_1-h_2)^{2}.$$ We obtain $$\mathsf {Contr}(n_1, n_2)=(-1)^{n+1} \cdot \Phi_1(h_1)^{n_1} \cdot \Phi_2(h_2)^{n_2} \cdot \Psi(h_1, h_2)\, .$$ The sign in the last equality is a consequence of rewriting the numerator of the normal bundle: 
$$(w+h_1-h_2)^{n_1-n_2+1}(-w-h_1+h_2)^{n_2-n_1+1}=(-1)^{n+1} (w+h_1-h_2)^{2}\, .$$ 
Therefore, we have \begin{eqnarray*}\mathsf A(q)&=&\sum_{n=0}^{\infty} q^n (-1)^n \cdot \int_{\text{Quot}_{\mathbb P^1}(\mathbb C^N, n)} c(T \text{Quot}_{\mathbb P^1}(\mathbb C^N, n))\cdot c(\mathcal T_n)\\&=&-\sum_{n=0}^{\infty} q^n \sum_{n_1+n_2=n} \int_{\mathbb P^{n_1} \times \mathbb P^{n_2}}  \Phi_1(h_1)^{n_1} \cdot \Phi_2(h_2)^{n_2} \cdot \Psi(h_1, h_2)\\ &=&-\sum_{n=0}^{\infty} \sum_{n_1+n_2=n} q^n \cdot \left[h_1^{n_1}\cdot h_2^{n_2}\right] \left( \Phi_1(h_1)^{n_1} \cdot \Phi_2(h_2)^{n_2} \cdot \Psi(h_1, h_2) \right)\, .
\end{eqnarray*} 
As before, the brackets indicate taking the suitable coefficient of the expression following it. Omitted from the notation is the fact that we also need to take the $w$-free term at the end.

The multivariable Lagrange-B\"urmann formula of   \cite[Theorem 2 (4.4)]{G} is:
 \begin{equation}\label{mvli}\sum_{n_1, n_2\geq 0} t_1^{n_1} t_2^{n_2} \cdot \left[h_1^{n_1}\cdot h_2^{n_2}\right] \left( \Phi_1(h_1)^{n_1} \cdot \Phi_2(h_2)^{n_2} \cdot \Psi(h_1, h_2) \right)=\frac{\Psi}{K}(h_1, h_2)\end{equation} for the change of variables $$t_1=\frac{h_1}{\Phi_1(h_1)}\, ,\  \ \ \ t_2=\frac{h_2}{\Phi_2(h_2)}$$ and for $$K(t_1, t_2)=\left(1-\frac{t_1}{\Phi_1(t_1)}\cdot \Phi_1'(t_1)\right)\cdot \left(1-\frac{t_2}{\Phi_2(t_2)}\cdot \Phi_2'(t_2)\right).$$ In our case, by \eqref{phipsi}, we have \begin{equation}\label{h1h2q}
t_1=\frac{h_1(h_1+w)}{(1-h_1)(1-h_1-w)}\,,\  \ \ \ 
 t_2=\frac{h_2(h_2-w)}{(1-h_2)(1-h_2+w)}\, .\end{equation} 
 Using \eqref{phipsi} again, by direct calculation,
 we find $\frac{\Psi}{K}(h_1, h_2)$ equals
 $$\frac{(1 - h_1^2) \cdot (1 - (w+ h_1)^2) \cdot  (1 - h_2^2)\cdot (1 -(w - 
     h_2)^2) \cdot(w + h_1 - h_2)^2 }{(2h_1^2+2h_1(w-1)+w(w-1))\cdot (2h_2^2-2h_2(w+1)+w(w+1))\cdot (1 -(w + h_1 - h_2)^2)}\, .$$ 
     \vskip.05in
\noindent We set $t_1=t_2=q$ and use the above equations \eqref{h1h2q} to solve 
$$h_1=-\frac{q}{1-q}-\frac{w}{2}+\sqrt{\frac{q}{(1-q)^2}+\frac{w^2}{4}},\,\,\,\, h_2=-\frac{q}{1-q}+\frac{w}{2}-\sqrt{\frac{q}{(1-q)^2}+\frac{w^2}{4}}.$$
A direct computation then shows that 
$$\frac{\Psi}{K} \left(h_1(q), h_2(q)\right)=-\frac{(1 - w^2) -4 q (2-w^2) + 4 q^2 (4 - w^2)}{(1 -q)^2 (1 - w^2 - 2 q (3-w^2) +   q^2 (1 - w^2))}$$ so that 
 $$\frac{\Psi}{K}\left(h_1(q), h_2(q)\right)\big{|}_{w=0}=     -\frac{(1-4q)^2}{(1-q)^2(1-6q+q^2)}\, .$$
 Therefore, $$\mathsf A(q)=-\sum_{n=0}^{\infty}\sum_{n_1+n_2=n} q^n \cdot \left[h_1^{n_1}\cdot h_2^{n_2}\right] \left( \Phi_1^{n_1} \cdot \Phi_2^{n_2} \cdot \Psi\right)=\frac{(1-4q)^2}{(1-q)^2(1-6q+q^2)}\, .$$
We have completed the proof of 
the $N=2$ case of Theorem \ref{rattt}. \qed
\subsubsection{The case $N>2$}. 
The calculation of $\mathsf{Z}_{X,N=2,0}^{\mathcal{E}}$
presented above can be exactly followed for all higher $N$.
The universal series $\mathsf{U}_N$ of Theorem \ref{rattt} is determined by
the equation
\begin{equation}\label{ddkk33}
  \mathsf{U}_N^{-1} = \sum_{n=0}^{\infty} q^n (-1)^{n} \cdot \int_{\text{Quot}_{\mathbb P^1}(\mathbb C^N, n)} c(T \text{Quot}_{\mathbb P^1}(\mathbb C^N, n))\cdot c(\mathcal T_n)\, ,
  \end{equation}
where $\mathcal{T}_n$ is the bundle
$$\mathcal T_n=\text{Ext}^{\bullet}_{\mathbb P^1}(Q, Q\otimes \mathcal{O}(-1) )\, .$$
Localization with respect to the diagonal $\mathbb{C}^{\star}$-action on $\mathbb C^N$
yields
\begin{equation*}\int_{\text{Quot}_{\mathbb P^1}(\mathbb C^N, n)} c(T \text{Quot}_{\mathbb P^1}(\mathbb C^N, n))\cdot c(\mathcal T_n)=\sum_{n_1+\ldots+n_N=n} \int_{\mathbb P^{n_1}\times\cdots \times \mathbb P^{n_N}}
  \mathsf {Contr}(n_1, \ldots, n_N)\, .\end{equation*}

By an explicit analysis of $\mathsf {Contr}(n_1, \ldots, n_N)$,
we can write
\begin{equation}\label{pxpxpx}
\mathsf {Contr}(n_1, \ldots,n_N)=(-1)^{n(N-1)+\binom{N}{2}} \cdot \Phi_1(h_1)^{n_1} \cdots \Phi_N(h_N)^{n_N} \cdot \Psi(h_1, \ldots, h_n)\, 
\end{equation}
for rational functions 
$$
\Phi_i(h_i)=\prod_{j=1}^{N} (1-h_i+w_i-w_j)\cdot \prod_{j\neq i} (h_i+w_j-w_i)^{-1}\, ,$$ 
\begin{multline*}\Psi=\prod_{i} (1+h_i)\cdot \prod_{i<j} (h_i-h_j+w_j-w_i)^2 \\
\cdot \prod_{j\neq i} (1+h_i+w_j-w_i)\cdot (1+h_i-h_j+w_j-w_i)^{-1}\cdot (h_i+w_j-w_i)^{-1}\, ,
\end{multline*}
 which depend upon $N$.
After applying the Lagrange-B\"urmann formula with $$t_i=\frac{h_i}{\Phi_i(h_i)}=h_i \cdot \prod_{j=1}^{N} (1-h_i+w_i-w_j)^{-1}\cdot \prod_{j\neq i} (h_i+w_j-w_i)\, ,$$ we find $$\sum_{n_1, \ldots, n_N}  t_1^{n_1}\cdots t_N^{n_N}\cdot \left([h_1^{n_1}\cdots h_N^{n_N}] \,\Phi_1(h_1)^{n_1} \cdots \Phi_N(h_N)^{n_N} \cdot \Psi(h_1, \ldots, h_n)\right) =\frac{\Psi}{K}(h_1, \ldots, h_N)\, .$$ After setting $$t_1=\ldots=t_N=q(-1)^N$$ the series  \eqref{ddkk33}  becomes \begin{equation}\label{d44} \mathsf U_N^{-1}=(-1)^{\binom{N}{2}} \cdot \frac{\Psi}{K} (h_1, \ldots, h_N)\,\end{equation} where $h_i$ solves the equation  $$q(-1)^N=  \prod_{j=1}^{N} \frac{h_i+w_j-w_i}{1-h_i+w_i-w_j}\, .$$ 
We must select the analytic solution $h_i(q)$ with $$h_i|_{q=0}=0.$$ 

We prefer however to work with a single equation. Let $H_1, \ldots, H_N$ be all solutions to the $i=1$  equation 
$$q(-1)^N=\prod_{j=1}^{N} \frac{h+w_j-w_1}{1-h+w_j-w_1}$$ with initial values $H_j(q=0)=w_1-w_j$. Then, by direct computation, we see that
 $$h_i=H_i+w_i-w_1$$ solves the $i^{\text{th}}$ equation. By \eqref{d44}, we obtain \begin{equation}\label{d55} U_N^{-1}=(-1)^{\binom{N}{2}} \cdot \frac{\Psi}{K} \left(H_1, H_2+w_2-w_1, \ldots, H_N+w_N-w_1\right)\,.\end{equation}
Using the explicit expressions of $\Psi$ and $K$, we see that the right hand side of \eqref{d55} is symmetric in $H_1, \ldots, H_N$.
 Since  symmetric functions in $H_1, \ldots, H_N$ are rational functions in $w$ and $q$ (with possible poles at $q=1$), the same is true of $\mathsf U_N^{-1}$. 

In fact, there are no poles of $\mathsf U_N^{-1}$ at $w=0$. Indeed, after 
setting the equivariant weights to $0$, 
the series \eqref{d55} is expressed as a symmetric rational function in the
$N$ roots $h_i=r_i$ of the polynomial equation
\begin{equation}\label{poly}q(-1)^N=h^{N}(1-h)^{-N}.\end{equation} 

A direct computation shows that the expression \eqref{d55} becomes $$\mathsf U_N^{-1}= (-1)^{\binom{N}{2}}\cdot \frac{\prod_{i=1}^{N} \left((1-r_i)\cdot (1+r_i)^{N} \right) \cdot \prod_{i<j} (r_i-r_j)^2}{N^N (r_1\cdots r_N)^{N-1}}\cdot \prod_{i<j} (1-(r_i-r_j)^2)^{-1}.$$ We write $$f(h)= \frac{h^N-(h-1)^Nq}{1-q}=\prod_{i=1}^N (h-r_i)$$ for the normalized equation \eqref{poly}. Then, 
\begin{eqnarray*}
\prod_{i=1}^N (1+r_i)&=& (-1)^Nf(-1)=\frac{1-2^Nq}{1-q}\, , \\ 
\prod_{i=1}^{N} (1-r_i)\cdot \frac{\prod_{i<j} (r_i-r_j)^2}{N^N(r_1\cdots r_N)^{N-1}}&=& (-1)^{\binom{N}{2}} \prod_{i=1}^N \frac{(1-r_i)\cdot f'(r_i)}{Nr_i^{N-1}}=(-1)^{\binom{N}{2}}\cdot \prod_{i=1}^{N}\frac{1}{1-q}\, .
\end{eqnarray*}
 Therefore, we find \begin{equation}\label{al}\mathsf U_N=\frac{(1-q)^{2N}}{(1-2^Nq)^{N}} \cdot \prod_{i\neq j} (1-(r_i-r_j)^2)\, .\end{equation}

We can easily calculate $\mathsf{U}_N$ for each $N$
from formula \eqref{al} by elementary algebra. For instance
\begin{eqnarray}
\nonumber 
\mathsf{U}_3
&=&\frac{(1-q)^2(1 - 22 q + 150 q^2 - 22 q^3 + q^4)}{(1-8q)^3}\, ,\\
\mathsf{U}_4 \nonumber
&=&\frac{(1-q)^2(1 - 62 q + 1407 q^2 - 15492 q^3 + 1407 q^4 - 62 q^5 + q^6)}{(1-16q)^4}\, .
\end{eqnarray}
Moreover, since \eqref{al} is clearly a symmetric
rational function of  the roots $r_1, \ldots, r_N$, the series
 $\mathsf{U}_N$
is a rational function in the elementary symmetric functions 
of the roots 
and hence a rational function of $q$. \qed

\subsection {Proof of Theorem \ref{cor1}.}\label{later} The methods of 
Section \ref{ves} can also be used to give a proof of 
the second part of Theorem \ref{cor1}: $$\mathsf B_{1,N}(q)=\frac{(1+t)^{N+1}}{1+(N+1)t} \ \ \ \ \text{ for } \ q=(-1)^N t(1+t)^N.$$ The first part of
 Theorem \ref{cor1}  was proven in Section \ref{binnum}. 

Recall the  $\mathsf A$ and $\mathsf B$-series
 defined by $$\sum_{n=0}^{\infty} q^n \int_{\mathsf {Quot}_{\,C} (\mathbb C^N, n)} s(L^{[n]})=\mathsf A^{\deg\, L}_{1,1,N}\cdot \mathsf B_{1,N}(q)^{1-g}\, ,$$ 
for a line bundle $L\to C$. After specializing to $C=\mathbb P^1$ and $L=\mathcal O_{\mathbb P^1}$,
 we obtain $$\mathsf B_{1,N}(q)=\sum_{n=0}^{\infty} q^n \int_{\mathsf {Quot}_{\mathbb P^1}(\mathbb C^N, n)} s(\mathcal O^{[n]})\, .$$ As usual, we set $\mathsf{B}(q)=\mathsf{B}_{1,N}(q)$ for notational convenience.

Consider the standard
$\mathbb{C}^{\star}$-action on $\mathsf {Quot}_{\mathbb P^1}(\mathbb C^N, n)$ with weights $w_1, \ldots, w_N$. In order to keep the notation manageable, we
specialize to $N=2$  (the argument for arbitrary $N$ is exactly parallel). 
By localizing, we obtain $$\mathsf B(q)=\sum_{n=0}^{N} q^n 
 \sum_{n_1+n_2=n} \mathsf {Contr}(n_1, n_2)\, ,$$ 
where each fixed locus $\mathbb P^{n_1}\times \mathbb P^{n_2}$ contributes $$\mathsf {Contr}(n_1, n_2)=\int_{\mathbb P^{n_1}\times \mathbb P^{n_2}} \frac{s(\mathcal O^{[n_1]}[w_1]) \cdot s(\mathcal O^{[n_2]}[w_2])}{\mathsf e_{\mathbb C^{\star}} (\mathsf N[n_1, n_2])}\, .$$ Using Lemma \ref{segr}, we find $$s(\mathcal O^{[n_i]}[w_i])=(1-w_i)^{-1}\cdot (1-h_i+w_i)^{-n_i+1}\, .$$ 
For the normal bundle, we use equation \eqref{nn1n2}:
 $$\mathsf e_{\mathbb C^{\star}} (\mathsf N[n_1, n_2])=(-1)^{n+1} \cdot (h_1+w_2-w_1)^{n_1+1} \cdot (h_2+w_1-w_2)^{n_2+1} \cdot (h_1-h_2+w_2-w_1)^{-2}\, .$$ 
We define $$\Phi_1(h_1)=\left(1-h_1+w_1\right)^{-1}\cdot (h_1+w_2-w_1)^{-1}\, ,$$ 
$$\Phi_2(h_2)=\left(1-h_2+w_2\right)^{-1}\cdot (h_2+w_1-w_2)^{-1}\, ,$$
\begin{multline*}
\Psi= \left(1-w_1\right)^{-1}\cdot \left(1-w_2\right)^{-1}\cdot\left(1-h_1+w_1\right) \cdot \left(1-h_2+w_2\right)\\
\cdot (h_1+w_2-w_1)^{-1} \cdot (h_2+w_1-w_2)^{-1}
\cdot (h_1-h_2+w_2-w_1)^2\, .
\end{multline*}
Therefore, $$\mathsf {Contr}(n_1, n_2)=(-1)^{n+1}\cdot \int_{\mathbb P^{n_1}\times\mathbb P^{n_2}} \Phi_1(h_1)^{n_1}\cdot \Phi_2(h_2)^{n_2} \cdot\,\Psi(h_1, h_2)$$
which gives  
$$\mathsf B(q)= - \sum_{n=0}^{\infty}\sum_{n_1+n_2=n} (-q)^n \cdot \left([h_1^{n_1}\cdot h_2^{n_2}] \,\,\Phi_1(h_1)^{n_1}\cdot \Phi_2(h_2)^{n_2} \cdot\,\Psi(h_1, h_2)\right).$$ Using Lagrange-B\"urmann inversion, we find $$\mathsf B(q)= - \frac{\Psi}{K}(h_1, h_2)$$ for the change of variables $$-q=\frac{h_1}{\Phi_1(h_1)}=h_1\cdot \left(1-h_1+w_1\right)\cdot(h_1+w_2-w_1)\, ,$$ 
$$-q=\frac{h_2}{\Phi_2(h_2)}=h_2\cdot \left(1-h_2+w_2\right)\cdot (h_2+w_1-w_2)\,
.$$
The first of the two equations $$-q=h\cdot (1-h+w_1)\cdot (h+w_2-w_1)$$ has two solutions $H_1(q)$ and $H_2(q)$ with $$H_1(0)=0\, ,\ \ \  H_2(0)=w_1-w_2\,.$$ 
 The root of the second equation $$-q=\frac{h_2}{\Phi_2(h_2)}=h_2\cdot \left(1-h_2+w_2\right)\cdot (h_2+w_1-w_2)$$ with initial value $0$ at $q=0$ is then 
\begin{equation}\label{mm9933}
\widetilde H_2(q)=H_2(q)+w_2-w_1\, .
\end{equation} Equation \eqref{mm9933}
 is easily seen by direct substitution. We conclude $$\mathsf B(q)=-\frac{\Psi}{K}\left(H_1(q), \widetilde H_2(q)\right)=-\frac{\Psi}{K}\left(H_1(q), H_2(q)+w_2-w_1\right).$$ Further direct calculation shows  
$$\mathsf B(q)=-\frac{\Psi}{K} \left(H_1, H_2+w_2-w_1\right)$$ equals $$-\frac{(1 - H_1 + w_1)^2\cdot (1 - H_2 + w_1)^2\cdot (H_1 - H_2)^2}{\prod_{i=1}^{2}(1-w_i)\cdot (3 H_i^2 - 2H_i \cdot (1 + 2 w_1 -w_2)+(1 + w_1) \cdot (w_1 - w_2))}\, .$$
We finally take the limit $w_1, w_2\to 0$.  
Write $\mathsf h_1, \mathsf h_2$ for the two roots of the equation 
$$-q=h^2(1-h)\, ,\ \ \  \mathsf h_1(0)=\mathsf h_2(0)=0\, .$$ These are power series in $q^{1/2}$. In the limit $w_i\to 0$, we obtain 
 $$\mathsf B(q)= -\frac{(1-\mathsf h_1)^2\cdot (1-\mathsf h_2)^2\cdot (\mathsf h_1-\mathsf h_2)^2}{(3\mathsf h_1-2)\cdot (3\mathsf h_2-2)\cdot (\mathsf h_1\mathsf h_2)}\, .$$ 

For general $N$, a similar analysis yields $$\mathsf B(q)=(-1)^{\binom{N+1}{2}} \cdot \frac{\prod_{i<j} (\mathsf h_i-\mathsf h_j)^2\cdot (\mathsf h_1\cdots \mathsf h_N)^{-(N-1)} \cdot \prod_{i} (1-\mathsf h_i)^2} {\prod_i ((N+1)\mathsf h_i-N)}$$ where $\mathsf h_1, \ldots, \mathsf h_N$ solve the equation \begin{equation}\label{qnn}(-1)^{N-1}q=h^N(1-h)\, ,\ \ \  \mathsf h_i(0)=0\, .\end{equation}

Equation \eqref{qnn} has an additional solution $\mathsf h(q)$ with $\mathsf h(0)=1$, which we can express in simple form. Indeed, if $$q=(-1)^N t(1+t)^N\, ,$$ then by direct verification $$\mathsf h(q)=1+t\, .$$ 
To complete the proof of Theorem \ref{cor1}, we must
show $$\mathsf B(q)=\frac{\mathsf h^{N+1}}{(N+1)\mathsf h-N}=\frac{(1+t)^{N+1}}{1+(N+1)t}\, .$$ Equivalently, we prove the identity
 \begin{equation}\label{ee}(-1)^{\binom{N+1}{2}} \cdot \frac{\prod_{i<j} (\mathsf h_i-\mathsf h_j)^2\cdot (\mathsf h_1\cdots \mathsf h_N)^{-(N-1)} \cdot \prod_{i} (1-\mathsf h_i)^2} {\prod_i ((N+1)\mathsf h_i-N)}=
\frac{\mathsf h^{N+1}}{(N+1)\mathsf h-N}\, .\end{equation}

The identity \eqref{ee} is straightforward to check. Let $$f(h)=h^N(h-1)-q(-1)^N=(h-\mathsf h) \prod_{i=1}^{N} (h-\mathsf h_i)\, .$$ After setting $h=1$,
 we obtain 
$$\prod_{i=1}^{N} (1-\mathsf h_i)=-\frac{q(-1)^N}{1-\mathsf h}\, .$$ 
We compute $$f'(h)=h^{N-1}((N+1)h-N)\, .$$ 
In particular, we find $$f'(\mathsf h_i)=(\mathsf h_i-\mathsf h)\cdot \prod_{j\neq i} (\mathsf h_i-\mathsf h_j)=\mathsf h_i^{N-1} ((N+1)\mathsf h_i-N)$$ $$f'(\mathsf h)=\prod_{i=1}^{N} (\mathsf h-\mathsf h_i)=\mathsf h^{N-1}((N+1)\mathsf h-N).$$ Therefore, $$(-1)^{\binom{N+1}{2}} \prod_{i<j} (\mathsf h_i-\mathsf h_j)^2=\frac{\prod_{i=1}^{N} f'(\mathsf h_i)}{f'(\mathsf h)}=\frac{\prod_{i=1}^{N} \mathsf h_i^{N-1}((N+1)\mathsf h_i-N)}{\mathsf h^{N-1}((N+1)\mathsf h-N)}\, .$$ 
After substitution, the left hand side of equation \eqref{ee} becomes $$\left(-\frac{q(-1)^N}{1-\mathsf h}\right)^2\cdot \frac{1}{\mathsf h^{N-1}((N+1)\mathsf h-N)}=\frac{\mathsf h^{N+1}}{(N+1)\mathsf h-N}\, ,$$ where equation \eqref{qnn} was used in the last step. \qed
\vskip.1in

The same method can be used to determine the series $\mathsf B_{r,N}$ 
for arbitrary values of $r=\text{rank }(V)$. While in general the formulas are less explicit, for $\text{rank }(V)=2$ and $N=2$, we obtain $$\mathsf B_{2,2}(-t^2)= \frac{\left(1+\sqrt{1-4t}\right)^4\cdot \left(1+\sqrt{1+4t}\right)^4 \cdot \left(1-\sqrt{1-16t^2}\right)} {2048 t^2\cdot \sqrt{1-16t^2}}\, .$$
\section{Virtual invariants of surfaces: dimension 1 quotients} \label{vis1}

\subsection{Overview} Let $X$ be a nonsingular, simply connected,
projective surface, and let $D$ an effective divisor on $X$. We compute here invariants associated to the scheme $\mathsf {Quot}_{\,X}(\mathbb C^N, n,D)$ of 
short exact sequences $$0\to S\to \mathbb C^N\otimes \mathcal O_X\to Q\to 0,\,\,\,\ \  \chi(Q)=n\, ,\  c_1(Q)=D\, ,\ \text{rank}(Q)=0\, .$$ In particular, we will prove Proposition \ref{t6}, Theorem \ref{t7}, and Theorem \ref{genttt}.
 
\subsection {Tangent-obstruction theory} \label{ves1}

Since $X$ is simply connected, the Hilbert scheme of curves is isomorphic to $$\mathsf {Quot}_{\,X}(\mathbb C^1, n,D)\simeq X^{[m]}\times \mathbb P$$
where $\mathbb P=|D|$. Here $$m=n+\frac{D(D+K_X)}{2}=n+(g-1)\, ,$$
where $g$ is the genus of a nonsingular curve in the linear series $|D|$. Indeed, to each pair $(Z, C)$ with $C\in |D|$, we can associate the sequence $$0\to I_Z(-C)\to \mathcal O_X\to Q\to 0\, .$$
While the actual dimension is $2m+ h^0(D)-1$, the expected dimension of the Hilbert scheme equals $$m+\frac{D(D-K_X)}{2}\, .$$ The first term $m$ comes from the Hilbert scheme of points, while the second is the virtual dimension of $|D|$ endowed with its natural obstruction theory as a Hilbert scheme. 

We calculate the tangent-obstruction theory, following \cite {segre}, in case $m>0$. Let
$$\mathcal L\to |D|  \ \text{ and }\ \mathcal Z\subset X^{[m]}\times X$$ denote
the tautological bundle $\mathcal{O}_{|D|}(1)$
and the universal subscheme of the Hilbert scheme
respectively. Over $\mathsf {Quot}_{\,X}(\mathbb C^1, n,D)$, we compute \begin{eqnarray*}\text{Tan}-\text{Obs}&=&\text{Ext}^{\bullet}(\mathcal S, \mathcal Q)\\&=&\text{Ext}^{\bullet}(\mathcal I_{\mathcal Z}\otimes \mathcal O_X(-D)\otimes \mathcal L^{-1}, \mathcal O-\mathcal I_{\mathcal Z}\otimes\mathcal O_X(-D)\otimes \mathcal L^{-1})\\ &=&\text{Ext}^{\bullet}(\mathcal I_{\mathcal Z}\otimes \mathcal O_X(-D)\otimes \mathcal L^{-1}, \mathcal O) - \text{Ext}^{\bullet}(\mathcal I_{\mathcal Z}, \mathcal I_{\mathcal Z})\\&=&\text{Ext}^{\bullet}(\mathcal I_{\mathcal Z}\otimes \mathcal O_X(-D), \mathcal O)\otimes \mathcal L- \text{Ext}^{\bullet}(\mathcal I_{\mathcal Z}, \mathcal I_{\mathcal Z})\\&=&H^\bullet (X, \mathcal O_X(D))\otimes \mathcal L -\text{Ext}^{\bullet}(\mathcal O_{\mathcal Z}\otimes \mathcal O_X(-D), \mathcal O)\otimes \mathcal L -\text{Ext}^{\bullet}(\mathcal I_{\mathcal Z}, \mathcal I_{\mathcal Z})\, .
                                                                                                                                                               \end{eqnarray*}
                                                                                                                                                               Two further calculations are needed. First, $$\text{Ext}^{\bullet}(I_Z, I_Z) = \text{Ext}^0(I_Z, I_Z)-\text{Ext}^1(I_Z, I_Z)+\text{Ext}^2(I_Z, I_Z)=\mathbb C-TX^{[m]}+H^0(K_X)^{\vee}\, ,$$ where we have used that $X$ is simply connected and Serre duality in the second equality. Second,
                                                                                                                                                               \begin{eqnarray*}
   \text{Ext}^{\bullet}(\mathcal O_{Z}\otimes \mathcal O_X(-D), \mathcal O)&=&\text{Ext}^{0}(\mathcal O_{Z}\otimes \mathcal O_X(-D), \mathcal O)
                                                               -\text{Ext}^{1}(\mathcal O_{Z}\otimes \mathcal O_X(-D), \mathcal O)\\
                                                                                                                                                                 & &\hspace{40pt} +\text{Ext}^{2}(\mathcal O_{Z}\otimes \mathcal O_X(-D), \mathcal O)\\
                                                                                                                                                                 &=&H^0(K_X(-D)\otimes \mathcal O_Z)^{\vee}\, ,\end{eqnarray*}
                                                                                                                                                               where we used vanishing for dimension reasons and Serre duality. 
Substituting, we find 
\begin{eqnarray*}\text{Tan}-\text{Obs}&=&H^\bullet (X, \mathcal O_X(D))\otimes \mathcal L-  \left((K_X(-D))^{[m]}\right)^{\vee}\otimes \mathcal L+T{X^{[m]}}-\mathbb C-H^0(K_X)^{\vee}\\&=& T\,{\mathbb P}  - H^1(X, \mathcal O_X(D))\otimes \mathcal L + H^2(X, \mathcal O_X(D))\otimes \mathcal L\\
  & & \hspace{20pt} -  \left((K_X(-D))^{[m]}\right)^{\vee}\otimes \mathcal L +T{X^{[m]}}  - H^0(K_X)^{\vee}.
\end{eqnarray*} For the second equality, we have also 
used the Euler sequence $$0\to \mathcal O\to H^0(X, \mathcal O(D)) \otimes \mathcal L\to \text{Tan\,}_{\mathbb P}\to 0\, .$$ In conclusion, we see that the $K$-theory class of the obstruction bundle equals $$H^1(X, \mathcal O_X(D))\otimes \mathcal L - H^2(X, \mathcal O_X(D))\otimes \mathcal L +  \left((K_X(-D))^{[m]}\right)^{\vee}\otimes \mathcal L+H^0(K_X)^{\vee}.$$ After setting $M=K_X-D$, we can rewrite the obstruction bundle as \begin{equation}\label{obun}\text{Obs}= (H^1(M)-H^0(M)+M^{[m]})^{\vee}\otimes \mathcal L + H^0(K_X)^{\vee}.\end{equation}
By the definition of the virtual Euler characteristic,
\begin{eqnarray*}\label{spl}\mathsf e^{\mathrm{vir}} (\mathsf {Quot}_{\,X}(\mathbb C^1, n,D))&=&\int_{X^{[m]}\times\mathbb P}\mathsf e(\text{Obs}) \,\frac{c(TX^{[m]}) \,c(T\mathbb P)}{c(\text{Obs})}\, .\end{eqnarray*} 
\vskip.1in
\subsection{Examples ($N=1$)} We illustrate the calculations above by examples corresponding to several different geometries. \vskip.1in

\subsubsection {\it Rational surfaces} A rich theory is obtained when $X$ is a rational surface. Since $H^0(K_X)=0$
for rational surfaces,  the obstruction bundle simplifies to $$\text{Obs}= (H^1(M)-H^0(M)+M^{[m]})^{\vee}\otimes \mathcal L.$$

\noindent {\it Proof of Proposition \ref{t6}}. Let $X$ be the blow-up of a rational surface at one point
with exceptional divisor $E$. Take $D=E$ so that $$\text{Obs} = \left(M^{[m]}\right)^{\vee}.$$ Thus for $n=m+1$, 
$$\mathsf e^{\mathrm{vir}} (\mathsf {Quot}_{\,X}(\mathbb C^1, n, E))=\int_{X^{[m]}} \mathsf e\left(\left(M^{[m]}\right)^{\vee}\right) \frac{c(TX^{[m]})}{c\left(\left(M^{[m]}\right)^{\vee}\right)}.$$ Such integrals have been computed in equation \eqref {xm} of Corollary \ref{pro2}. We find 
$$\sum_{n=1}^{\infty} q^{n-1} \mathsf e^{\mathrm{vir}} (\mathsf {Quot}_{\,X}(\mathbb C^1, n, E))=\left((1-q)^{2} (1-2q)^{-1}\right)^{K_X^2+1}.$$ \vskip.1in

\subsubsection {K3 surfaces.}\label{k3k3k3} Let $X$ be a $K3$ surface, and let $D$ be a primitive big and nef curve class. In particular, we have 
$$H^0(M)=H^1(M)=0.$$
We write $g$ for the genus of $D$.
The obstruction bundle $$\text{Obs} = \left(M^{[m]}\right)^{\vee}\otimes \mathcal L + H^0(K_X)^{\vee}$$ has a trivial summand. As a result, all virtual invariants vanish. 

A reduced obstruction bundle can be defined by removing the trivial factor. With the new obstruction theory, we find \begin{equation}\label{eul}\mathsf e^{\mathrm{red}} (\mathsf {Quot}_{\,X}(\mathbb C^1, n, D))= \int_{X^{[m]}\times\mathbb P}\mathsf e\left(\left(M^{[m]}\right)^{\vee}\otimes \mathcal L\right) \,\frac{c(TX^{[m]}) \,c(T\mathbb P)}{c(\left(M^{[m]}\right)^{\vee}\otimes \mathcal L)},\end{equation} for $M=\mathcal O_X(-D)$. \vskip.1in

\noindent {\it Proof of Theorem \ref{t7}.} Without the dual placed on the tautological bundle $M^{[m]}$, integrals similar to \eqref{eul} also
appear in G\"ottsche's conjecture 
 and are computed by the  Kawai-Yoshioka formula \eqref{kawaiy}:  \begin{equation}\label{got} N_ {g, n}=\int_{X^{[m]}\times\mathbb P}\mathsf e\left(D^{[m]}\otimes \mathcal L\right) \,\frac{c(TX^{[m]}) \,c(T\mathbb P)}{c(D^{[m]}\otimes \mathcal L)}\, .\end{equation}  This was noted in \cite [Section 4]{KST}.

To prove the claim of Theorem \ref{t7}, $$\mathsf e^{\mathrm{red}} (\mathsf {Quot}_{\,X}(\mathbb C^1, n, D))=N_{g, n}\, ,$$ we will use formulas \eqref{eul} and \eqref{got}. Since $X^{[m]}$ is holomorphic symplectic, we can replace the tangent bundle in \eqref{eul} with the isomorphic cotangent bundle. 
Thus, we must show
\begin{multline*}\int_{X^{[m]}\times\mathbb P}\mathsf e\left(\left(M^{[m]}\right)^{\vee}\otimes \mathcal L\right) \,\frac{c\left(\left(TX^{[m]}\right)^{\vee}\right) \,c(T\mathbb P)}{c(\left(M^{[m]}\right)^{\vee}\otimes \mathcal L)}\\
  =\int_{X^{[m]}\times\mathbb P}\mathsf e\left(D^{[m]}\otimes \mathcal L\right) \,\frac{c(TX^{[m]}) \,c(T\mathbb P)}{c(D^{[m]}\otimes \mathcal L)}\, .
\end{multline*} 
After integrating out the hyperplane class on $\mathbb P$, we are led to the statement 
$$\int_{X^{[m]}} \mathsf P\left(c_i\left(\left(M^{[m]}\right)^{\vee}\right), c_j\left(\left(TX^{[m]}\right)^{\vee}\right)\right)=\int_{X^{[m]}} \mathsf P\left(c_i(D^{[m]}), c_j(TX^{[m]})\right)$$ where $\mathsf P$ is a uniquely defined universal polynomial in the Chern classes of various tautological bundles on the Hilbert scheme $X^{[m]}$. After removing the duals (since $X^{[m]}$ is even dimensional), we must show \begin{equation}\label{symmetric} \int_{X^{[m]}} \mathsf P\left(c_i(M^{[m]}), c_j(TX^{[m]})\right)=\int_{X^{[m]}} \mathsf P\left(c_i(D^{[m]}), c_j(TX^{[m]})\right).\end{equation}

Equality \eqref{symmetric} is then a consequence of \cite[Theorem 4.1]{EGL}.
Expressions such as the ones in \eqref{symmetric} are given by universal formulas in the Chern numbers. For the left hand side, these Chern numbers are $$c_1(M)^2\, , \ \ K_X^2\, , \ \ c_1(M)\cdot K_X\, ,\ \  c_2(X)\, .$$
The right hand side is similar, with the relevant numbers being
$$c_1(D)^2\, , \ \ K_X^2\, ,  \ \ c_1(D)\cdot K_X\, ,\ \ c_2(X)\, .$$
Since $X$ is a $K3$ surface, all the Chern numbers match, including 
$$c_1(M)\cdot K_X=c_1(D)\cdot K_X=0\, ,$$
which may in general sign change.
\qed

\vspace{10pt}
The case $D=0$ is not covered by Theorem \ref{t7}. However, in the $K3$ case, we can consider the reduced theory of the Hilbert scheme of points $X^{[n]}$ obtained by removing the canonical trivial factor from the obstruction bundle $\left(\mathcal O^{[n]}\right)^{\vee}.$ The reduced virtual dimension is
$n+1$, and the obstruction bundle equals $$\text{Obs}=\left(\mathcal O^{[n]}-\mathcal O\right)^{\vee}\to X^{[n]}\, .$$
While the question does not involve any curve classes, the calculation below makes use of Theorem \ref{t7} for curves of genus $1$. 

\begin{proposition} We have $$\sum_{n=1}^{\infty} q^n \mathsf e^{\mathrm{red}} (X^{[n]})= \frac{24q}{(1-q)^2}\, .$$
\end{proposition} 
\proof We have $$\mathsf e^{{red}}(X^{[n]})=\int_{X^{[n]}} \mathsf e(\text{Obs}) \cdot \frac{c(TX^{[n]})}{c(\text{Obs})}=\int_{X^{[n]}} \mathsf e\left(\left(\mathcal O^{[n]}-\mathcal O\right)^{\vee}\right)\frac{c(TX^{[n]})}{c\left(\left(\mathcal O^{[n]}-\mathcal O\right)^{\vee}\right)}\, .$$
For $n>0$, we will prove
\begin{equation}\label{ccllmm}
  \mathsf e^{\mathrm{red}}(X^{[n]})=N_{1, n}\, .
  \end{equation}

We start by writing $\alpha_1, \ldots, \alpha_n$ for the roots of $\mathcal O^{[n]}$
with the convention that $\alpha_1=0$ corresponds to the trivial summand of the obstruction bundle.
Then, claim \eqref{ccllmm}
becomes \begin{equation}\label{claimeq}\int_{X^{[n]}} \prod_{i=2}^{n} \frac{-\alpha_i}{1-\alpha_i}\cdot c(TX^{[n]})=N_{1, n}\, .\end{equation}
 
By deformation invariance, we may assume $X$ is an elliptically fibered $K3$ with fiber class $f$. We apply Theorem \ref{t7} for the curve class $D=f$. The associated line bundle $D$ has no higher cohomology, and the proof of Theorem \ref{t7} applies even though
$D$ is not big. 
We find  $$ N_{1,n}=
\int_{X^{[n]}\times \mathbb P^1}\mathsf e\left(\left(M^{[n]}\right)^{\vee}\otimes \mathcal L\right) \,\frac{c(TX^{[n]}) \,c(T\mathbb P^1)}{c(\left(M^{[n]}\right)^{\vee}\otimes \mathcal L)}\, ,$$
where $M=\mathcal O_X(-f)$.

We write $\mu_1, \ldots, \mu_n$ for the roots of $M^{[n]}$, and let $\zeta$ be the hyperplane class on the projective line. The above integral becomes
\begin{eqnarray*}N_{1, n}&=&\int_{X^{[n]}\times \mathbb P^1}
                             \prod_{i=1}^{n} \frac{\zeta-\mu_i}{1+\zeta-\mu_i} \cdot c(TX^{[n]}) (1+2\zeta)\\
                         &=&\int_{X^{[n]}} \left(2 \prod_{i=1}^{n} \frac{-\mu_i}{1-\mu_i} + \sum_{i=1}^{n} \frac{1}{(1-\mu_i)^2}\prod_{j\neq i} \frac{-\mu_j} {1-\mu_j} \right) c(TX^{[n]})\, ,
\end{eqnarray*}
where, in  the second equality, we have integrated out the hyperplane class on $\mathbb P^1$.
The resulting integral is a universal polynomial in the quantities
\begin{equation}\label{kk339}
  M^2\, , \ \ M\cdot K_X\, , \ \ K_X^2\, , \ \ c_2(X)\, .
  \end{equation} Indeed, the 
expression
$$2 \prod_{i=1}^{n} \frac{-\mu_i}{1-\mu_i} + \sum_{i=1}^{n} \frac{1}{(1-\mu_i)^2}\prod_{j\neq i} \frac{-\mu_j} {1-\mu_j} $$ can be written in terms of the Chern classes of $M^{[n]}$. The claimed universality then follows from \cite[Theorem 4.1]{EGL}.

Since the four numerical invariants \eqref{kk339}  are the same if $M=-f$ or $M=0$,
we are free to replace the $\mu_i$'s by the $\alpha_i$'s without changing the answer. Therefore,
$$N_{1, n}=\int_{X^{[n]}} \left(2 \prod_{i=1}^{n} \frac{-\alpha_i}{1-\alpha_i} + \sum_{i=1}^{n} \frac{1}{(1-\alpha_i)^2}\prod_{j\neq i} \frac{-\alpha_j} {1-\alpha_j} \right) c(TX^{[n]})\, .$$
Since $\alpha_1=0$, we obtain
$$N_{1, n}=\int_{X^{[n]}} \prod_{j\neq 1} \frac{-\alpha_j} {1-\alpha_j}\cdot c(TX^{[n]})\, ,$$ as claimed in \eqref{claimeq}.

Finally, using the Kawai-Yoshioka formula \eqref{kawaiy}, we find $$N_{1, n} =[q\cdot y^n]\left(\sqrt{y}-\frac{1}{\sqrt{y}}\right)^{-2}\prod_{n=1}^{\infty} \frac{1}{(1-q^n)^{20}(1-yq^n)^2(1-y^{-1}q^n)^2}=24n\, ,$$ for $n>0$.  
\qed

\subsubsection {Surfaces of general type.} 
\label{ddssgg}
Let $X$ be a nonsingular, simply connected, projective
surface of general type with $p_g(X)>0$. 
For $D=K_X$, 
the obstruction bundle \eqref{obun} takes the form
$$\text{Obs}= (\mathcal O^{[m]}-\mathcal O)^{\vee}\otimes \mathcal L + H^0(K_X)^{\vee}\, .$$ Due to the presence of a trivial summand, the virtual Euler characteristic vanishes $$\mathsf e^{\mathrm{vir}} (\mathsf {Quot}_{\,X}(\mathbb C^1, n,K_X))=0$$ for $m>0$. The case $$m=0\iff n=-K_X^2$$ is special, yielding the answer $$\mathsf e^{\mathrm{vir}} (\mathsf {Quot}_{\,X}(\mathbb C^1, n,K_X))=\int_{\mathbb P} \frac{1}{c(\mathcal L)}= (-1)^{p_g+1}=(-1)^{\chi(\mathcal O_X)},$$ in agreement with \cite {CK, DKO}. 
\vskip.1in

\noindent {\it Proof of Proposition \ref{pro22}.} Let $D$ be an arbitrary
effective curve class. To start, we take  $N=1$ and assume 
$$D\neq 0\, , \ \ \ D\neq K_X\, ,$$
since these cases have already been considered. Recall  $$\text{Obs}= (H^1(M)-H^0(M)+M^{[m]})^{\vee}\otimes \mathcal L + H^0(K_X)^{\vee}.$$ If $M=K_X-D \text{ is not effective}$, then $H^0(M)=0$. The virtual class is then forced
to vanish by the trivial summand $H^0(K_X)^{\vee}$ of the obstruction
bundle.

We may therefore assume $M$ to be effective. 
By Serre duality,
 $$\text{rank}\, \text{Obs}=h^1(D)-h^2(D)+m+p_g\, .$$ 
We also have
 $$\left[\mathsf {Quot}_X(\mathbb C^1, n, D)\right]^{\mathrm{vir}}=\mathsf e(\text{Obs})\cap \left[\mathbb P\times X^{[m]}\right]\, ,$$ 
where $$\dim \mathbb P=h^0(D)-1\, .$$ 
We write $h\in A^1(\mathbb P)$ for the hyperplane class 
and $\alpha_i$ for the Chern roots of $M^{[m]}$ on $X^{[m]}$.
The virtual class then equals the degree $p_g+m+h^1(D)-h^2(D)$ part of  
\begin{equation}\label{exxe}
c(\text{Obs})=(1+h)^{h^1(D)-h^2(D)}\cdot \prod_{i=1}^{m}(1+h-\alpha_i)\, .
\end{equation} The expression \eqref{exxe}
 contains terms of the form $$h^{k}\cdot \text{symmetric polynomial of degree at most }m\text{ in the roots }\alpha_i\, ,$$ where $k\leq h^0(D)-1$ for dimension reasons. All terms therefore have degree bounded by \begin{equation}\label{eeee}h^0(D)-1+m<h^1(D)-h^2(D)+m+p_g\, .\end{equation} Consequently, the Euler class vanishes. 

To justify inequality \eqref{eeee}, we use the following chain of equivalences:
\begin{eqnarray*}
h^0(D)-1<h^1(D)-h^2(D)+p_g&\iff& \chi(D)<1+p_g=\chi(\mathcal O_X)\\
&\iff& D\cdot (D-K_X)<0\\
&\iff& D\cdot M>0\, .
\end{eqnarray*}
 The inequality $D\cdot M>0$ holds
 since the pair $(D, M)$ is a nontrivial effective splitting of $K_X$ 
(the canonical class is $1$-connected for minimal surfaces of general type \cite[Proposition 6.1] {BPV}). The proof of Proposition \ref{pro22} for
 $N=1$ is complete.

For $N>1$, we use $\mathbb C^\star$-equivariant localization. 
The natural $\mathbb C^*$-action on $\mathsf {Quot}_{\,X}(\mathbb C^N, n, D)$ has fixed loci $$\mathsf F[(n_1, D_1), \ldots, (n_N, D_N)]$$ indexed by all possible effective splittings $$n_1+\ldots+n_N=n\, , \ \ \ D_1+\ldots+D_N=D\, .$$ 
The corresponding subsheaves are $$S=\bigoplus_{i=1}^{N} I_{Z_i}\hookrightarrow \mathbb C^N\otimes \mathcal O_X\, , \ \ \ c_1(Z_i)=D_i\, , \ \ \
 \chi(\mathcal O_{Z_i})=n_i\, .$$ 
The induced virtual class of $$\mathsf F[(n_1, D_1), \ldots, (n_N, D_N)]=\mathsf {Quot}_{\,X}(\mathbb C^1, n_1,D_1)\times \cdots \times  \mathsf {Quot}_{\,X}(\mathbb C^1, n_N,D_N)$$ is determined by the fixed part of $$\text{Ext}^{\bullet}(S, Q)^{\text{fix}}=\bigoplus \text{Ext}^{\bullet}(I_{Z_i}, \mathcal O_{Z_i})$$ and, therefore, splits over the factors. Using the case $N=1$ already established, 
in order to obtain a nontrivial virtual fundamental class on the ${i}^{\text{th}}$ factor, we must have $$D_i=0 \ \text{ or }\  D_i=K_X\ \ \ \implies\ \ \ D=\ell K_X \ \text{ for }\ 0\leq \ell\leq N\, .$$ 
By the paragraph preceding the proof of Proposition \ref{pro22}, 
the choice $D_i=K_X$ forces $Z_i$ to be supported only on canonical curves, without any point contributions. 
\qed

\subsection {Proof of Theorem \ref{genttt}.} 
\label{jjj999}
The $N=1$ case of Theorem \ref{genttt} is a consequence of the calculations
of Section \ref{ddssgg}. In the $N=2$ case, Theorem \ref{genttt} can be derived
from Theorem \ref{thm1}: the localization contributions can be expressed as integrals over the symmetric product with $7$ Segre factors.{\footnote{We leave the
argument to the intrepid reader.}} However, we will treat all the
cases $N\geq 1$ together using the strategy of the
the proof Theorem \ref{rattt}.

Let $X$ be a nonsingular,
simply connected, minimal surface of general type admitting a nonsingular canonical curve 
$C\subset X$ of genus $$g=K_X^2+1\, .$$ 
Let $0\leq \ell\leq N$. Let 
$$\mathsf Z_{X, N, \ell K_X}^{\mathcal E}(q)=\sum_{n\in \mathbb Z} q^n \mathsf e^{\mathrm{vir}}(\mathsf {Quot}_{X}(\mathbb C^N, n, \ell K_X))\, .$$ 
The formula of Theorem \ref{genttt} is
$$\mathsf Z_{X, N, \ell K_X}^{\mathcal E}(q)=(-1)^{\ell\cdot \chi(\mathcal O_X)} \ q^{\ell(1-g)}\cdot \sum_{1\leq i_1<\ldots<i_{N-\ell}\leq N} \mathsf A(r_{i_1}, \ldots, r_{i_{N-\ell}})^{1-g}\, .$$ 
The sum is taken over all $\binom{N}{N-\ell}$
 choices of  $N-\ell$ distinct roots of the equation 
$$z^N=q(z-1)^N\, .$$ Furthermore,  
$$\mathsf A(x_1, \ldots, x_{N-\ell})= \frac{(-1)^{\binom{N-\ell}{2}}}{N^{N-\ell}} \cdot \prod_{i=1}^{N-\ell} \frac{(1+x_i)^N (1-x_i)}{x_i^{N-1}}\cdot \prod_{i<j} \frac{(x_i-x_j)^2}{1-(x_i-x_j)^2}\, .$$ 
In case $\ell=N$, the formula is interpreted as
$$\mathsf Z_{X, \,N, \,N K_X}^{\mathcal E}(q)=(-1)^{N\cdot \chi(\mathcal O_X)} \, q^{N(1-g)}\, .$$

To prove the claimed evaluation, 
we consider the $\mathbb{C}^{\star}$-action on $\mathsf {Quot}_{X} (\mathbb C^N,n, \ell K_X)$ with weights $w_1,\ldots, w_N$ on the middle term of the sequence $$0\to S\to\mathbb C^N\otimes \mathcal O_X\to Q\to 0\, .$$ 
We write $$n=m+\ell(1-g)\, .$$ 
For convenience, we set $$k=N-\ell\, .$$ By the last sentence in the proof of Proposition \ref{pro22}, the contributing fixed loci correspond to kernels of the form $$S=\bigoplus_{i=1}^{\ell} \mathcal O_X(-D_i)\oplus \bigoplus_{j=1}^{k} I_{Z_j}\hookrightarrow \mathbb C^N\otimes \mathcal O_X\, ,$$ where $D_i\in |K_X|$ and $Z_j$ is a $0$-dimensional scheme of length $m_j$. Of course, we have $$\sum_{j=1}^{k} m_j=m\, .$$ 
 The weights $w_1, \ldots, w_N$ are distributed over the summands of $S$ in $\binom{N}{k}$ possible ways, depending on the location of the curves and points. 
The fixed loci are therefore
indexed by tuples $(m_1, \ldots, m_k)$ as well as choices of $\binom{N}{k}$ 
summands of $\mathbb C^N$. For a fixed partition $(m_1, \ldots, m_k)$, 
there are  $\binom{N}{k}$ fixed loci all isomorphic to
 $$\mathsf F[m_1, \ldots, m_k]=\left(\prod_{i=1}^{\ell} \mathbb P\right)\times \left(\prod_{j=1}^{k}X^{[m_j]}\right)\, .$$ 
Here, $\mathbb P$ denotes the linear series $|K_X|$. The obstruction bundle splits into obstruction bundles over the  factors, $$\text{Obs}=\sum_{i=1}^{\ell} \text{pr}^{\star}_i \left(H^0(K_X)^{\vee}-\mathcal L\right) + \left(\sum_{j=1}^{k}K_X^{[m_j]}\right)^{\vee}\,.$$ We therefore obtain 
\begin{eqnarray}\label{ss99}\nonumber
\left[\mathsf F[m_1, \ldots, m_k] \right]^{\mathrm{vir}}&=&\mathsf e\left(\sum_{i=1}^{\ell} \text{pr}^{\star}_i \left(H^0(K_X)^{\vee}-\mathcal L\right) + \left(\sum_{j=1}^{k}K_X^{[m_j]}\right)^{\vee}\right)\\ \nonumber
&=&\left[\prod_{i=1}^{\ell} \text{pr}_i^{\star} \frac{1}{1+c_1(\mathcal L)}\right]_{\left(\sum_{i=1}^{\ell}\dim \mathbb P\right)} \cdot \prod_{j=1}^{k}\mathsf e\left(\left(K_X^{[m_j]}\right)^{\vee} \right)
\\ \nonumber
&=&(-1)^{\ell \chi} \cdot (-1)^m \iota_\star\left( \text{[pt]} \times \cdots \times \text{[pt]} \times \left[C^{[m_1]}\times \cdots \times C^{[m_k]}\right]\right).
\end{eqnarray}
The subscript on the second line indicates that we need to select the terms of top degree. In addition, $\chi=\chi(\mathcal{O}_X)$, and,
for the canonical curve $C\subset X$, we have written
 $$\iota: \text{[pt]} \times \cdots \times \text{[pt]} \times \left(C^{[m_1]}\times \cdots \times C^{[m_k]} \right) \hookrightarrow \left(\prod_{i=1}^{\ell} \mathbb P\right)\times \left(X^{[m_1]}\times \cdots \times X^{[m_k]}\right)$$ for the natural morphism.  

We write $$j :\mathsf F[m_1, \ldots, m_k]\hookrightarrow \mathsf {Quot}_{X} (\mathbb C^N, n,\ell K_X)$$ for the natural inclusion. The integral
$$\mathsf e^{\mathrm{vir}}  \left(\mathsf {Quot}_{X} (\mathbb C^N, n,\ell K_X)\right)= \int_{\left[\mathsf {Quot}_{X} (\mathbb C^N, n,\ell K_X)\right]^{\mathrm{vir}}} c(T^{\mathrm{vir}} \mathsf {Quot}_{X})$$ can be calculated by $\mathbb C^\star$-equivariant localization. Each fixed locus $\mathsf F= \mathsf F[m_1, \ldots, m_k]$ yields
 a contribution 
\begin{equation}\label{eevv}
\int_{\left[\mathsf F[m_1, \ldots, m_k]\right]^{\mathrm{vir}}} \frac{c(j^{\star} \,T^{\mathrm{vir}}{\mathsf {Quot}_X})}{\mathsf e(\mathsf N^{\mathrm{vir}})}= (-1)^{m+\ell \chi} \int_{C^{[m_1]}\times \cdots \times C^{[m_k]}} \iota^{\star} \left(\frac{c(T^{\mathrm{vir}}\mathsf F)\,c(\mathsf N^{\mathrm{vir}})}{\mathsf e(\mathsf N^{\mathrm{vir})}}\right).\end{equation}

We will analyze these contributions separately. We assume
 the weights $w_1, \ldots, w_{\ell}$ are distributed on the curve summands and
the weights  $w_{\ell+1}, \ldots, w_N$ are distributed on the point summands. In other words, the kernels are $$S=\bigoplus_{i=1}^{\ell} \mathcal O_X(-D_i)[w_i] \oplus \bigoplus_{j=1}^{k} I_{Z_j} [w_{j+\ell}]\, .$$ We will 
use the indices $i, i'$ to refer to the curve summands, while the
indices $j, j'$ will be reserved for the point summands. We obtain 
\begin{eqnarray*}
T^{\mathrm{vir}}\mathsf F&=& \sum_{i=1}^{\ell} T\mathbb P+ \sum_{j=1}^{k} TX^{[m_j]}-\text{Obs}\\
&=&\sum_{i=1}^{\ell} T\mathbb P+ \sum_{j=1}^{k} TX^{[m_j]}-\left(\sum_{i=1}^{\ell} \text{pr}_i^{\star} \left(H^0(K_X)^{\vee}-\mathcal L\right) + \sum_{j=1}^{k} \left((K_X)^{[m_j]}\right)^{\vee}\right)\, ,
\end{eqnarray*}
 which yields 
\begin{eqnarray*}
\iota^{\star} T^{\mathrm{vir}}\mathsf F
&=& \sum_{i=1}^{\ell}  \mathbb C^{p_g-1} + \sum_{j=1}^{k} \iota^{\star} TX^{[m_j]} - \left(\sum_{i=1}^{\ell} \mathbb C^{p_g-1} +\sum_{j=1}^{k} \iota^{\star} \left((K_X)^{[m_j]}\right)^{\vee}\right)\\
&=&\sum_{j=1}^{k} \left(TC^{[m_j]}+\Theta^{[m_j]} - \left(\Theta^{[m_j]}\right)^{\vee}\right)\\
&=&\sum_{j=1}^{k} \left(\left(K_C^{[m_j]}\right)^{\vee}+\Theta^{[m_j]} - \left(\Theta^{[m_j]}\right)^{\vee}\right)\,.
\end{eqnarray*}
Here, $\Theta=\mathcal O_X(C)|_{C}$ is the theta characteristic. The last equality was shown in the proof of Theorem \ref{rattt}, see \eqref{qpp22}. 
There are no equivariant weights for $\iota^{\star} T^{\mathrm{vir}}\mathsf F$. 
  
 The virtual normal bundle splits into four terms $$\mathsf N^{\mathrm{vir}}=\mathsf N_1+ \mathsf N_2+\mathsf N_3+\mathsf N_4$$ where  
 \begin{eqnarray*}\mathsf N_1&=& \sum_{i=1}^{\ell}\sum_{j=1}^{k} \text{Ext}^{\bullet}(\mathcal O_X(-D_i), \mathcal O_{Z_j})[w_{j+\ell}-w_i]\,,\\ 
 \mathsf N_2 &=&  \sum_{i=1}^{\ell}\sum_{j=1}^{k} \text{Ext}^{\bullet}(I_{Z_j}, \mathcal O_{D_i})[w_i-w_{j+\ell}]\,,\\ \mathsf N_3&=& \sum_{i=1}^{\ell}\sum_{i'\neq i} \text{Ext}^{\bullet}(\mathcal O_X(-D_i), \mathcal O_{D_{i'}})[w_{i'}-w_i]\,,\\ \mathsf N_4&=&\sum_{j=1}^{k} \sum_{j'\neq j}\text{Ext}^{\bullet}(\mathcal I_{Z_j}, \mathcal O_{Z_{j'}})[w_{j'+\ell}-w_{j+\ell}]\,.\end{eqnarray*} 
We would normally include the tautological line bundle $\mathcal L$ in the expression of the subsheaf, but, since we are in the end restricting to a point via $\iota$, there is no need.

We write $\mathsf N_1^{ij}$ for the $ij$-summand of $\mathsf N_1$. We find 
\begin{eqnarray*}
\iota^{\star} \mathsf N^{ij}_{1}&=&\iota^{\star}(K_X)^{[m_j]}[w_{j+\ell}-w_i] =\Theta^{[m_j]}[w_{j+\ell}-w_i]\, .\end{eqnarray*}
Similarly
\begin{eqnarray*}
\iota^{\star} \mathsf N^{ij}_{2}&=&\iota^{\star}\left(\text{Ext}^{\bullet}(\mathcal O-\mathcal O_{Z_j}, \mathcal O-K_X^{-1})[w_i-w_{j+\ell}]\right)\\
&=&\iota^{\star} \left(H^{\bullet}(\mathcal O)-H^{\bullet}(K_X^{-1})-\left((K_X)^{[m_j]}\right)^{\vee}+\left((K_X^{\otimes {2}})^{[m_j]}\right)^{\vee}\right)[w_i-w_{j+\ell}]\, ,
\end{eqnarray*}
where we have used, suppressing indices, that $$\text{Ext}^{\bullet}(\mathcal O_Z, \mathcal O)=\text{Ext}^{2-\bullet}(\mathcal O, K_X\otimes \mathcal O_Z)^{\vee}=\left((K_X)^{[m]}\right)^{\vee}\, ,$$ $$\text{Ext}^{\bullet}(\mathcal O_Z, K_X^{-1})=\text{Ext}^{2-\bullet}(K_X^{-1}, K_X\otimes \mathcal O_Z)^{\vee}=\left((K_X^{\otimes {2}})^{[m]}\right)^{\vee}\, .$$ 
Since $$H^{\bullet}(\mathcal O)-H^{\bullet}(K_X^{-1})=-\mathbb C^{g-1}\, ,$$ 
we have $$\iota^{\star}  \mathsf N_2^{ij}=-\mathbb C^{g-1}[w_i-w_{j+\ell}]-\left(\Theta^{[m_j]}\right)^{\vee}[w_i-w_{j+\ell}]+\left(K_C^{[m_j]}\right)^{\vee}[w_i-w_{j+\ell}]\, .$$ 
For the third term of the virtual normal bundle,
$$\iota^{\star}\mathsf N_3^{ii'}=\left(H^{\bullet}(\mathcal O(D_i))-H^{\bullet}(\mathcal O(D_i-D_i'))\right)[w_{i'}-w_i]=\left(H^{\bullet}(K_X)-H^{\bullet}(\mathcal O_X)\right)[w_{i'}-w_i]=0.$$ For
the fourth term, we have already
computed in equation \eqref{no} of the proof of Lemma \ref{l8},  
for  $j\neq j'$,
 $$\iota^{\star}\left( \mathsf N_4^{jj'}\right)=\mathsf T_{jj'}+\mathsf N_{jj'}\, ,$$ where \begin{eqnarray*}\mathsf T_{jj'}&=&\text{Ext}^{\bullet}_C(\mathcal O_{Z_j}, \mathcal O_{Z_{j'}}\otimes \Theta)[w_{j'+\ell}-w_{j+\ell}]\, ,\\ \mathsf N_{jj'}&=&\text{Ext}^{\bullet}_C(\mathcal I_{Z_j}, \mathcal O_{Z_{j'}})[w_{j'+\ell}-w_{j+\ell}]
\, .\end{eqnarray*} 
We also had observed there
 that, as a consequence of Serre duality,
 $$\mathsf e(\mathsf T_{jj'}+\mathsf T_{j'j})=1\,.$$ 
Moreover, $\mathsf N=\sum_{j\neq j'}\mathsf N_{jj'}$ is identified with the 
 normal bundle of the fixed locus $$C^{[m_1]}\times \cdots \times C^{[m_k]}\hookrightarrow \mathsf {Quot}_C(\mathbb C^{k}, m)\, ,$$ where the $\mathbb C^\star$-action has weights $w_{\ell+1}, \ldots, w_N$ on $\mathbb C^k$. 
 
After collecting all terms, the fixed locus contribution becomes 
$$(-1)^{m+\ell \chi}\cdot \int_{C^{[m_1]}\times \ldots \times C^{[m_k]}}\prod_{j=1}^{k} \frac{
c\left(\left(K_C^{[m_j]}\right)^{\vee}\right) \cdot c(\Theta^{[m_j]})}{c\left(\left(\Theta^{[m_j]}\right)^{\vee}\right)}\cdot \prod_{i=1}^{\ell}\prod_{j=1}^{k} \frac{c(\Theta^{[m_j]}[w_{j+\ell}-w_i])}{\mathsf e(\Theta^{[m_j]}[w_{j+\ell}-w_i])}$$ $$\cdot \prod_{i=1}^{\ell}\prod_{j=1}^{k} \left(\frac{\mathsf e\left(\left(\Theta^{[m_j]}\right)^{\vee}[w_i-w_{j+\ell}]\right)}{c\left(\left(\Theta^{[m_j]}\right)^{\vee}[w_i-w_{j+\ell}]\right)}\cdot \frac{c\left(\left(K_C^{[m_j]}\right)^{\vee}[w_i-w_{j+\ell}]\right)}{\mathsf e\left(\left(K_C^{[m_j]}\right)^{\vee}[w_i-w_{j+\ell}]\right)} \cdot \frac{(w_i-w_{j+\ell})^{g-1}}{(1+w_i-w_{j+\ell})^{g-1}}\right)$$ $$\cdot \prod_{1\leq j\neq j'\leq k} c(\mathsf T_{jj'})c(\mathsf N_{jj'})\cdot \frac{1}{\mathsf e(\mathsf N)}\, .$$
We note a cancellation between the Euler classes in the denominator of the second product and the numerator of the third product, yielding the answer \begin{equation}\label{expres}(-1)^{m+\ell\chi}\cdot (-1)^{\ell m}\cdot \int_{C^{[m_1}\times \ldots \times C^{[m_k]}}\prod_{j=1}^{k} \frac{
c\left(\left(K_C^{[m_j]}\right)^{\vee}\right) \cdot c(\Theta^{[m_j]})}{c\left(\left(\Theta^{[m_j]}\right)^{\vee}\right)} \end{equation} $$\cdot
\prod_{i=1}^{\ell}\prod_{j=1}^{k} \left(\frac{c(\Theta^{[m_j]}[w_{j+\ell}-w_i])}{c\left(\left(\Theta^{[m_j]}\right)^{\vee}[w_i-w_{j+\ell}]\right)}\cdot \frac{c\left(\left(K_C^{[m_j]}\right)^{\vee}[w_i-w_{j+\ell}]\right)}{\mathsf e\left(\left(K_C^{[m_j]}\right)^{\vee}[w_i-w_{j+\ell}]\right)} \cdot \frac{(w_i-w_{j+\ell})^{g-1}}{(1+w_i-w_{j+\ell})^{g-1}}\right)$$ $$\cdot \prod_{1\leq j\neq j'\leq k} c(\mathsf T_{jj'})c(\mathsf N_{jj'})\cdot \frac{1}{\mathsf e(\mathsf N)}\, .$$

Let $\mathsf {Contr }[m_1, \ldots, m_k]\in \mathbb Q((w))$
 denote the integral thus obtained (without including the sign
$(-1)^{(\ell+1)m+\ell\chi }$). We have
\begin{equation}\label{pp233}
\mathsf Z^{\mathcal E}_{X, N, \ell K_X}(q)=\sum_{n\in \mathbb Z} q^n \mathsf e^{\mathrm{vir}} \left(\mathsf {Quot}_{X} (\mathbb C^N, \ell K_X, n)\right)=\sum \mathsf Z[m_1, \ldots, m_k]\, ,
\end{equation}
 where $$\mathsf Z[m_1, \ldots, m_k] 
= (-1)^{(\ell+1)m+\ell\chi}  \, q^{\ell(1-g)+m}\cdot \mathsf{ Contr }[m_1, \ldots, m_k]\, .$$ As usual, the sum on the right
in the equation \eqref{pp233} has $\binom{N}{k}$ terms depending on the placement of the weights (and is also over $m_1,\ldots,m_k$). We will set $w=0$ at the end. 

We will transform the above contribution formulas
 into integrals over the Quot scheme $\mathsf {Quot}_C(\mathbb C^{k}, m)$. 
Recall, from the proof of Theorem \ref{rattt}, the virtual bundle $$\mathcal T_{m}=\text{Ext}^{\bullet}_C(Q, Q\otimes \Theta)$$ on $\mathsf {Quot}_C(\mathbb C^{k}, m)$.
The tautological bundle $$L^{[m]}\to \mathsf {Quot}_C(\mathbb C^{k}, m)$$ 
associated to a line bundle $L$ on $C$ was
defined in Section \ref{cccvvv}. We define 
\begin{multline*}\mathsf Z_{C, k}(q\,|\,w_1, \ldots, w_\ell\,|\, w_{\ell+1}, \ldots, w_N)=\\
\sum_{m=0}^\infty q^m \int_{\mathsf {Quot}_C(\mathbb C^{k}, m)} c(T\mathsf {Quot}_C(\mathbb C^{k}, m))\cdot c(\mathcal T_{m})\cdot \prod_{i=1}^{\ell}\left(\frac{c(\Theta^{[m]}[-w_i])}{c\left(\left(\Theta^{[m]}\right)^{\vee}[w_i]\right)}\cdot \frac{c\left(\left(K_C^{[m]}\right)^{\vee}[w_i]\right)}{\mathsf e\left(\left(K_C^{[m]}\right)^{\vee}[w_i]\right)}\,\right) .
\end{multline*}
 In the integrand, twists by trivial bundles with nontrivial equivariant weights are included. We consider the function above as a $\mathbb C^\star$-equivariant integral given by the $\mathbb C^\star$-action on the Quot scheme with weights $w_{\ell+1}, \ldots, w_{N}$.
 The function $\mathsf Z_{C, k}$ depends on $q$ and on the weights $w$. 
By an algebraic cobordism argument, we see 
 $$\mathsf Z_{C, k}=\mathsf A^{1-g}$$ where 
$$\mathsf A=\mathsf A(q\,|\,w_1, \ldots, w_\ell\,|\, w_{\ell+1}, \ldots, w_N)$$ is a universal function which does not depend on the genus $g$ of $C$.  

We will apply $\mathbb C^\star$-equivariant localization to the integrals appearing in the formula for $\mathsf Z_{C, k}$. The result is related to \eqref{expres}:
 each integral in $\mathsf Z_{C, k}$ becomes a sum of contributions of the fixed loci $$\iota:C^{[m_1]}\times \cdots\times C^{[m_k]}\hookrightarrow \mathsf {Quot}_C(\mathbb C^{k}, m)\, .$$
 We note the restrictions  
$$\iota^{\star}T\mathsf {Quot}_C(\mathbb C^{k}, m)=\sum_{j=1}^{k} TC^{[m_j]} + \sum_{j\neq j'}\mathsf N_{jj'}\, ,$$ $$\iota^{\star} \mathcal T_m=\sum_{j\neq j'} \mathsf T_{jj'}+\sum_{j=1}^{k} (\Theta^{[m_j]}-(\Theta^{[m_j]})^{\vee})\, .$$ Here, for $j=j'$, we have used  $$\text{Ext}^{\bullet}(\mathcal O_Z, \mathcal O_Z\otimes \Theta)=\Theta^{[m]}-\left(\Theta^{[m]}\right)^{\vee}.$$ 
Furthermore, the $\mathbb C^\star$-equivariant restrictions of the tautological bundles $\Theta^{[m]}$ on the Quot scheme to the fixed loci are given by $$\iota^{\star} \Theta^{[m]}[-w_i]=\sum_{j=1}^{k} \Theta^{[m_j]}[w_{j+\ell}-w_i]\, ,$$ $$\iota^{\star}\left(\Theta^{[m]}\right)^{\vee}[w_i]=\sum_{j=1}^{k} (\Theta^{[m_j]})^{\vee}[w_i-w_{j+\ell}]\, ,$$ where the sign $-w_{j+\ell}$ on the second line appears because the dual was taken. Finally, $$\iota^{\star}\left(K_C^{[m]}\right)^{\vee}[w_i]=\sum_{j=1}^{k} (K_C^{[m_j]})^{\vee}[w_i-w_{j+\ell}]\, .$$ 
The above $\mathbb C^\star$-equivariant localization terms of $\mathsf Z_{C, k}$ match expression \eqref{expres} up to a common factor and signs. Summarizing, we find:  
\begin{multline*}\mathsf Z^{\mathcal E}_{X, N, \ell K_X}(q)=(-1)^{\ell \chi}\, 
q^{\ell(1-g)}\cdot \sum  \left(\prod_{i=1}^{\ell}\prod_{j=1}^{k} \frac{1+w_i-w_{j+\ell}}{w_i-w_{j+\ell}}\right)^{1-g} \\
\cdot \mathsf A((-1)^{\ell+1}q \,|w_1, \ldots, w_{\ell}\, \,|w_{\ell+1}, \ldots, w_{N}\,)^{1-g}\, .
\end{multline*}
We write 
\begin{multline*}
\widetilde {\mathsf A}(q\,|w_1, \ldots, w_{\ell}\,|\, w_{\ell+1}, \ldots, w_N)= \prod_{i=1}^{\ell}\prod_{j=1}^{k}\frac{1+w_i-w_{j+\ell}}{w_i-w_{j+\ell}}\\
\cdot \mathsf A((-1)^{\ell+1}q \,|w_1, \ldots, w_{\ell}\, \,|w_{\ell+1}, \ldots, w_{N}\,)\, ,
\end{multline*} so we have
 $$\mathsf Z^{\mathcal E}_{X, N, \ell K_X}(q)=(-1)^{\ell \chi}\, q^{\ell(1-g)}\cdot \sum \widetilde {\mathsf A}(q\,|w_1, \ldots, w_{\ell}\,|\, w_{\ell+1}, \ldots, w_N)^{1-g}.$$ 

The last remaining step is to determine the function $\mathsf A$. 
After specializing the curve $C=\mathbb P^1$, 
we have
\begin{multline*}
\mathsf A(q\,|w_1, \ldots, w_{\ell}\,|\, w_{\ell+1}, \ldots, w_N)=\\
\sum_{m=0}^\infty q^m \int_{\mathsf {Quot}_{\mathbb P^1}(\mathbb C^{k}, m)} c(T\mathsf {Quot}_C(\mathbb C^{k}, m))\cdot c(\mathcal T_{m})\cdot \prod_{i=1}^{\ell}\left(\frac{c(\Theta^{[m]}[-w_i])}{c\left(\left(\Theta^{[m]}\right)^{\vee}[w_i]\right)}\cdot \frac{c\left(\left(K_C^{[m]}\right)^{\vee}[w_i]\right)}{\mathsf e\left(\left(K_C^{[m]}\right)^{\vee}[w_i]\right)}\right)\, .
\end{multline*}
All tautological structures in the above integral have been understood in the proof of Theorem \ref{rattt}. In fact, compared to the integrals which appear in the proof of Theorem \ref{rattt}, the only new terms are $$\prod_{i=1}^{\ell}\prod_{j=1}^{k} \left(\frac{c(\Theta^{[m_j]}[w_{j+\ell}-w_i])}{c\left(\left(\Theta^{[m_j]}\right)^{\vee}[w_i-w_{j+\ell}]\right)}\cdot \frac{c\left(\left(K_C^{[m_j]}\right)^{\vee}[w_i-w_{j+\ell}]\right)}{\mathsf e\left(\left(K_C^{[m_j]}\right)^{\vee}[w_i-w_{j+\ell}]\right)}\right)\,$$ considered over the product 
\begin{equation}\label{prpr}
\mathbb P^{m_1}\times \cdots \times \mathbb P^{m_k}\, .
\end{equation}

As before, we write  $h_1, \ldots, h_k$ for the hyperplane classes on
the respective projective spaces in the product \eqref{prpr}. 
Using Lemma \ref{segr}, we obtain $$c\left(\Theta^{[m]}[w]\right)=(1-h+w)^{m}\,,
 \ \ \ c\left((\Theta^{[m]})^{\vee}[-w]\right)=(1+h-w)^{m}\,,$$ 
$$c\left(\left(K_C^{[m}\right)^{\vee}[-w]\right)=\frac{(1+h-w)^{m+1}}{1-w}\,, 
\ \ \  \mathsf e\left(\left(K_C^{[m]}\right)^{\vee}[-w]\right)=\frac{(h-w)^{m+1}}{-w}\, .$$  
The new terms contribute the expression 
$$\prod_{i=1}^{\ell}\prod_{j=1}^{k} \left(\frac{(1-h_j+w_{j+\ell}-w_i)^{m_j}}{(h_j+w_i-w_{j+\ell})^{m_j}}\cdot \frac{1+h_j+w_i-w_{j+\ell}}{1+w_i-w_{j+\ell}}\cdot \frac{w_i-w_{j+\ell}}{h_j+w_i-w_{j+\ell}}\right)\,.$$ Therefore, using \eqref{pxpxpx},
 the contribution of the fixed locus of $\mathsf {Quot}_\mathbb P^1(\mathbb C^{k}, m)$ corresponding to the partition $(m_1, \ldots, m_k)$ equals $$(-1)^{m(k-1)+\binom{k}{2}} \int_{\mathbb P^{m_1}\times \ldots \times \mathbb P^{m_k}} \Phi_1(h_1)^{m_1} \cdots \Phi_k(h_k)^{m_k} \cdot \Psi(h_1, \ldots, h_k)\, $$ where 
\begin{multline*}
\Phi_j(h_j)=\prod_{j'=1}^{k} (1-h_{j}+w_{j+\ell}-w_{j'+\ell})\cdot \prod_{j'\neq j} (h_j+w_{j'+\ell}-w_{j+\ell})^{-1}\cdot\prod_{i=1}^{\ell} \frac{1-h_j+w_{j+\ell}-w_i}{h_j+w_i-w_{j+\ell}} \,,
\end{multline*}
and 
\begin{multline*}
 \Psi= \prod_{j'<j} (h_{j}-h_{j'}+w_{j'+\ell}-w_{j+\ell})^2 
\cdot \prod_{j, j'} (1+h_j+w_{j'+\ell}-w_{j+\ell})\cdot (1+h_j-h_{j'}+w_{j'+\ell}-w_{j+\ell})^{-1}\\ \cdot 
\prod_{j\neq j'}(h_j+w_{j'+\ell}-w_{j+\ell})^{-1}\cdot \prod_{i=1}^{\ell} \prod_{j=1}^{k}\left(\frac{1+h_j+w_i-w_{j+\ell}}{h_j+w_i-w_{j+\ell}}\cdot \frac{w_i-w_{j+\ell}}{1+w_i-w_{j+\ell}}\right)\,.
\end{multline*}
 Only the products involving $i$ are different from the expressions written in the proof of Theorem \ref{rattt}.

We now apply the Lagrange-B\"urmann formula
for the change of variables 
\begin{eqnarray*}
t_j&=&\frac{h_j}{\Phi_j(h_j)}=\prod_{\alpha=1}^{N} \frac{h_j+w_{\alpha}-w_{j+\ell}}{1-h_{j}+w_{j+\ell}-w_{\alpha}}.\end{eqnarray*} In the above product, the index $\alpha$ collects the terms in $\Phi_j$ corresponding to both $i$ and $j'$ into a uniform expression. 
We find 
$$\mathsf A=(-1)^{\binom{k}{2}}\cdot \frac{\Psi}{K} (h_1, \ldots, h_k)$$ where 
$$
q(-1)^{k-1}=t_j=\prod_{\alpha=1}^{N} \frac{h_j+w_{\alpha}-w_{j+\ell}}{1-h_{j}+w_{j+\ell}-w_{\alpha}}\, ,
$$
Let $$\widetilde \Psi=\prod_{j'<j} (h_{j}-h_{j'}+w_{j'+\ell}-w_{j+\ell})^2 \\
\cdot \prod_{j, j'} (1+h_j+w_{j'+\ell}-w_{j+\ell})\cdot (1+h_j-h_{j'}+w_{j'+\ell}-w_{j+\ell})^{-1}\cdot $$ $$\prod_{j\neq j'}(h_j+w_{j'+\ell}-w_{j+\ell})^{-1}\cdot \prod_{i=1}^{\ell} \prod_{j=1}^{k}\frac{1+h_j+w_i-w_{j+\ell}}{h_j+w_i-w_{j+\ell}}\,.$$ 
We find $$\widetilde{\mathsf A}(q\,|w_1, \ldots, w_{\ell}\,|\, w_{\ell+1}, \ldots, w_N)=(-1)^{\binom{k}{2}}\frac{\widetilde {\Psi}}{K}(h_1, \ldots, h_k)$$ where, 
taking all signs into account, we have 

$$q(-1)^{(k-1)+(\ell+1)}=\prod_{\alpha=1}^{N} \frac{h_j+w_{\alpha}-w_{j+\ell}}{1-h_{j}+w_{j+\ell}-w_{\alpha}}\,.$$ 
 In the limit $w\to 0$, the above
 equation becomes $$q(-1)^{N}=h^N(1-h)^{-N}\, .$$ 

The limit is justified as in the proof of Theorem \ref{rattt}: 
we let $H_1, \ldots, H_N$ be the roots of the single equation $$q(-1)^{N}=\prod_{\alpha=1}^{N} \frac{h+w_{\alpha}-w_{1}}{1-h+w_{1}-w_{\alpha}}\, ,$$ and then
 we have $$h_j=H_{j+\ell}+w_{j+\ell}-w_1\, .$$ 
The final answer is a sum of $\binom{N}{k}$ terms corresponding to choices of subsets of $k$ roots out of $H_1, \ldots, H_N$. Using the explicit expressions for $\tilde{\Psi}$ and $K$, the answer is seen to be symmetric in the $H$'s
and, therefore,  expressible in terms of the elementary symmetric functions which are polynomials in $w$.

We find $\frac{\widetilde {\Psi}}{K}$ simplifies in the limit to the expression $$\prod_{j<j'} (h_j-h_{j'})^{2}\cdot (1+h_j)^{k}\cdot \prod_{j,j'}(1-(h_{j'}-h_{j}))^{-1}\cdot \prod_{j=1}^{k} h_j^{-(k-1)}\cdot \prod_{j=1}^{k} \frac{(1+h_j)^{\ell}}{h_j^{\ell}}\cdot \prod_{j=1}^{k} \frac{1-h_j}{N}$$ where the last product comes from the $K$-term. 
Further simplification yields $$\frac{1}{N^k}\prod_{j<j'}\frac{(h_j-h_{j'})^2}{1-(h_j-h_{j'})^2} \cdot \prod_{j=1}^{k} \frac{(1+h_j)^N\cdot (1-h_j)}{h_j^{N-1}}\, ,$$ 
which is precisely the formula stated in Theorem \ref{genttt}.\qed

\appendix
\section{A combinatorial proof of Theorem \ref{thm3}} 
We present here
 a purely combinatorial argument for Theorem \ref{thm3}. For simplicity of notation, we consider trees whose edges are painted in only two colors denoted $A$ and $B$. The generalization to several colors does not require additional ideas. 

We write $a$ for the total number of $A$ edges, $b$ for the number of $B$ edges, and $n$ for the number of vertices of a tree $T$. 
Clearly $$a+b=n-1\, .$$ 
For each vertex $v$, we write $a_v$ and $b_v$ for the number of outgoing edges colored $A$ and $B$ respectively.  Therefore,
 $$\mathsf{wt}(T)=\frac{1}{(n-1)!} \prod_{v} a_v\,!\, b_v\,!\ .$$ 
We set $$w_n (a, b)= \sum_T \mathsf{wt}(T)\, .$$ Let $$t_n(a,b)=\frac{1}{n} \binom{2a+b}{a} \binom{a+2b}{b}\, .$$ 
\vspace{5pt}

\noindent The claim of Theorem \ref{thm3} in the case of two colors is
 \begin{equation}\label{ssxxss}
w_{n}(a, b)=t_n(a, b)\, .
\end{equation}

Define the generating series $$\mathsf W(q\,|\,x, y)=\sum_{n=1}^{\infty} \sum_{a+b=n-1}w_{n}(a, b)\cdot x^a y^b q^n\, ,$$ $$\mathsf T(q\,|\,x, y)=\sum_{n=1}^{\infty} \sum_{a+b=n-1}t_n(a, b)\cdot x^a y^b q^n\, .$$ 
By Lemmas \ref{kkk1} and \ref{kkk2} below,
 both $\mathsf W$ and $\mathsf T$ satisfy the cubic equation
\begin{equation}\label{gppg} 
\mathsf Z \cdot (1-x\mathsf Z)\cdot (1-y\mathsf Z)=q\, ,\ \ \ \mathsf Z|_{q=0}=0\, .
\end{equation}
 Since the solution of \eqref{gppg} is
unique, we obtain $$\mathsf W=\mathsf T\, $$ 
which implies \eqref{ssxxss} and completes the proof of Theorem \ref{thm3}. 

\begin{lemma} \label{kkk1} 
We have $$\mathsf T\cdot (1-x\mathsf T)\cdot (1-y\mathsf T)
=q\, .$$ \end{lemma}
\proof The argument exactly follows the proof
of Lemma \ref{loga}. Set $$f(t)=(1-xt)^{-1}(1-yt)^{-1}\, .$$ 
For $a+b=n-1$,
 we have $$t_n(a, b)=\frac{(-1)^{n-1}}{n} \binom{-n}{a}\binom{-n}{b}\, .$$ 
Therefore, \begin{eqnarray*}\mathsf T(q) &=& \sum_{n=1}^{\infty} \frac{q^n}{n} \cdot \left(\left[t^{n-1}\right] (1-xt)^{-n} (1-yt)^{-n}\right)\\
&=& \sum_{n=1}^{\infty} \frac{q^n}{n} \cdot \left(\left[t^{n-1}\right] f(t)^{n-1}\cdot f(t)\right)\, .\end{eqnarray*} 
Then, $$\frac{d\mathsf T}{dq} = \sum_{n=1}^{\infty} {q^{n-1}} \cdot \left(\left[t^{n-1}\right] f(t)^{n-1}\cdot f(t)\right)= \frac{dt}{dq}\, ,$$ where equation \eqref{lb} was used above for the change of variables $q=\frac{t}{f(t)}$.
Hence, we obtain 
$$\mathsf T=t\, ,$$ and  the change of variables proves the Lemma. \qed

\begin{lemma} \label{kkk2}
We have $$\mathsf W\cdot (1-x\mathsf W)\cdot (1-y\mathsf W)
=q\, .$$
\end{lemma}

\proof We will prove a recursion for $w_n(a, b)$ which implies
 the cubic equation of the Lemma. For convenience, we set $w_n(a, b)=0$ whenever $a+b\neq n+1$. 

Fix a labelled 2-colored  tree $T$ with $n$ vertices. 
Consider the vertex $\star$ with the {\it highest} label $n$. 
After removing the vertex $\star$ and all its incident edges from the tree $T$,
 we obtain
disjoint subtrees $T_1, \ldots, T_{\ell}$. 
We set up the following notation:
\begin{itemize}
\item[$\bullet$] $r$ and $s$ denote the number of edges incident
to the vertex $\star$ which are 
 colored  $A$ and $B$ respectively   (where $\ell=r+s$),
\item[$\bullet$] $n_1, \ldots, n_{\ell}$ are the number of vertices of the subtrees $T_1, \ldots, T_{\ell}$ respectively,
\item[$\bullet$] $(a_1, b_1), \ldots, (a_\ell, b_\ell)$ are the numbers 
of edges of each color for subtrees $T_1, \ldots, T_{\ell}$. 
\end{itemize}

The above quantities satisfy various constraints which
are most easily expressed using partitions. 
We denote an ordered partitions by
 $$\alpha^{\bullet}=(\alpha_1, \ldots, \alpha_\ell)\, ,$$
 and we write $|\alpha^{\bullet}|$ for the sum of parts. Then 
$$|n^{\bullet}|=n-1\, , \ \ \ |a^{\bullet}|+r=a\,\, \text{(counting }A\, \text{edges)} \, , \ \ \  |b^{\bullet}|+s=b\,\, \text{(counting }B\, \text{edges)}\, .$$ 

The removal of the vertex $\star$ yields the
following recursion:
 \begin{equation}\label{wre}w_n(a, b)=\sum \eta_{n^{\bullet}}\cdot w_{n_1}(a_1, b_1)\cdots w_{n_{\ell}}(a_\ell, b_\ell)\end{equation} with
 the combinatorial factor $$\eta_{n^{\bullet}}= \frac{r!}{\text{Aut}(n_1, \ldots, n_r)} \cdot \frac{s!}{\text{Aut}(n_{r+1}, \ldots, n_{\ell})}\, .$$ Here, $\text{Aut}(n^{\bullet})$ counts the automorphisms of the partition $n^{\bullet}$, and thus can be expressed as a product of factorials determined by the repetitions amongst the parts of $n^{\bullet}$.
To justify equation \eqref{wre}, note that the vertex $\star$ contributes $ r!\,s!$ to the weight of $T$, while the other vertices are contained in one of the trees $T_1, \ldots, T_\ell$. Therefore $$\mathsf{wt}(T)= \frac{1}{(n-1)!} r\,!\,s! \cdot \prod_{j=1}^{\ell} (n_j-1)! \,\text{wt}(T_j)\, .$$ 
After summing over all trees, we obtain \begin{eqnarray*} w_n(a, b)&=&\sum_{T} \mathsf{wt}(T)=\sum c_{n^{\bullet}} \cdot \frac{1}{(n-1)!} r\,!\,s! \cdot \prod_{j=1}^{\ell} (n_j-1)!\, w_{n_j}(a_j, b_j)\, .\end{eqnarray*} 
The combinatorial factor $$c_{n^{\bullet}}=\left(n_1\cdots n_{\ell}\right)\cdot \binom{n-1}{n_1, \ldots, n_{\ell}}\cdot \frac{1}{\text{Aut}(n_1, \ldots, n_r)}\cdot \frac{1}{\text{Aut}(n_{r+1}, \ldots, n_{\ell})}$$ arises as follows
\begin{itemize}
\item [$\bullet$] the term $n_1\cdots n_{\ell}$ counts all possible ways to attach the vertex $\star$ to one of the $n_j$ vertices of the tree $T_j$, for $1\leq j\leq \ell$,
\item [$\bullet$] $\binom{n-1}{n_1, \ldots, n_{\ell}}$ counts all possible ways of distributing the labels $\{1, \ldots, n-1\}$ to the trees $T_1, \ldots, T_{\ell}$, 
\item [$\bullet$] the last two terms account for automorphisms.  
\end{itemize}
Equation \eqref{wre} then follows by collecting terms. 

For notational convenience, we define the
  relabelling $$n'_{j}=n_{j+r}\, ,\ \ \
 a_j'=a_{j+r}\, ,  \ \ \ b'_j=b_{j+r},\,\,\, 1\leq j\leq s\, .$$
In the new notation, 
recursion \eqref{wre} takes the form:
 $$w_n(a, b)=\sum \frac{r!}{\text{Aut}(n^{\bullet})}\cdot \frac{s!}{\text{Aut}(n'^{\bullet})}\cdot \prod_{j=1}^{r} w_{n_j}(a_j, b_j)\cdot \prod_{j=1}^{s} w_{n'_j}(a'_j, b'_j)\, .$$ 
We define $$\mathsf W_n=\sum_{a+b=n-1} w_{n}(a, b) \cdot x^{a}y^{b}$$ satisfying 
$$\mathsf W=
\sum_{n=1}^{\infty} \sum_{a+b=n-1}w_{n}(a, b)\cdot x^a y^b q^n=\sum_{n=1}^{\infty} q^n \mathsf W_n\, .$$ 
We compute \begin{eqnarray*} \frac{\mathsf W}{q}&=&\sum_{n=1}^\infty 
w_n(a, b) x^{a} y^{b} q^{n-1}\\&=&\sum \frac{r!}{\text{Aut}(n^{\bullet})}\cdot \frac{s!}{\text{Aut}(n'^{\bullet})}\cdot \prod_{j=1}^{r} w_{n_j}(a_j, b_j)\cdot \prod_{j=1}^{s} w_{n'_j}(a'_j, b'_j)\cdot x^{|a^{\bullet}|+|a'^{\bullet}|} y^{|b^{\bullet}|+|b'^{\bullet}|}x^{r}y^{s} q^{|n^{\bullet|}+|n'^{\bullet}|}\\&=& \sum  \frac{r!}{\text{Aut}(n^{\bullet})}\cdot \frac{s!}{\text{Aut}(n'^{\bullet})} \cdot \prod_{j=1}^{r} \mathsf W_{n_j} \cdot \prod_{j=1}^{s} \mathsf W_{n'_j} \cdot x^r y^s q^{|n^{\bullet|}+|n'^{\bullet}|}\\&=& \left(\sum  \frac{r!}{\text{Aut}(n^{\bullet})} \prod_{j=1}^{r} \mathsf W_{n_j} \cdot x^r q^{|n^{\bullet}|}\right)\cdot \left(\sum  \frac{s!}{\text{Aut}(n'^{\bullet})} \prod_{j=1}^{s} \mathsf W_{n'_j} \cdot y^s q^{|n'^{\bullet}|}\right)\\&=& (1-x\mathsf W)^{-1}\cdot (1-y\mathsf W)^{-1}\, ,\end{eqnarray*}  
where, on the third line, we have summed over the $a$'s and $b$'s. 

For the last line, we have
used the identity \begin{equation}\label{bin}\frac{1}{1-x\mathsf W}=\sum  \frac{r!}{\text{Aut}(n^{\bullet})} \prod_{j=1}^{r} \mathsf W_{n_j} \cdot x^r q^{|n^{\bullet}|}\, \end{equation} which is easily derived from the Binomial Theorem. 
Indeed, after setting $$\alpha_n=\mathsf W_n \cdot x q^n\,, \  \ \  
\alpha=\sum \alpha_n=x\mathsf W\, ,$$ 
equation \eqref{bin} becomes $$\frac{1}{1-\alpha}=\sum  \frac{r!}{\text{Aut}(n^{\bullet})} \prod_{j=1}^{r} \alpha_{n_j}\, ,$$ 
which is true since the two sides are different ways of expressing $\sum_r \alpha^r$.
\qed


\begin{thebibliography}{99}

\bibitem {AJLOP}

N. Arbesfeld, D. Johnson, W. Lim, D. Oprea, R. Pandharipande, {\it The virtual $K$-theory of Quot schemes of surfaces}, J. Geom. Phys, \texttt{DOI: 10.1016/j.geomphys.2021.104154}


\bibitem {A}
J.C. Aval, {\it Multivariate Fuss-Catalan numbers}, Discrete Math. {\bf 308} (2008), 4660--4669

\bibitem{Ber} A. Bertram, {\it Towards a Schubert calculus for maps
from a Riemann surface to a Grassmannian}, Internat. J. Math {\bf 5}
(1994), 811--825

\bibitem {BDW}
A. Bertram, G. Daskalopoulos and R. Wentworth, {\it Gromov invariants for
holomorphic maps
from Riemann surfaces to Grassmannians}, J. Amer. Math. Soc. {\bf 9} (1996), 529--571

\bibitem {BPV}

W. Barth, C. Peters,  A. Van de Ven, {\it Compact complex surfaces}, Springer-Verlag, 1984

\bibitem{CK}

H.L. Chang, Y.H. Kiem, {\it Poincar\'e invariants are Seiberg-Witten invariants}, 
Geom. Topol. {\bf 17} (2013), 1149--1163

\bibitem{C}

E. Cotterill, {\it Geometry of curves with exceptional secant planes: linear
series along the general curve}, Math. Z. {\bf 267} (2011), 549--582


\bibitem {DKO}

M. D\"urr, A. Kabanov, C. Okonek, {\it Poincar\'e invariants},  Topology 
{\bf 46} (2007), 225--294

\bibitem {EGL}

G. Ellingsrud, L. G\"ottsche, M. Lehn, {\it On the
cobordism class of the Hilbert scheme of a surface}, J. Alg. Geom. 
{\bf 10} (2001), 81--100

\bibitem {EL} 

G. Ellingsrud, M. Lehn, {\it The irreducibility of the punctual Quot scheme}, preprint 1997

\bibitem {GF}

B. Fantechi, L. G\"ottsche, {\it Riemann-Roch theorems and elliptic genus for virtually smooth schemes}, Geom. Topol. {\bf 14} (2010), 83--115

\bibitem {G}

I. Gessel, {\it A combinatorial proof of the multivariable Lagrange inversion formula}, J. Comb. Theory, Ser. A,  {\bf 45} (1987), 178--195


\bibitem {GSY}

A. Gholampour, A. Sheshmani, S.T. Yau, {\it Nested Hilbert schemes on surfaces: Virtual fundamental class}, Advances in Math, {\bf 365} (2020), \texttt{DOI: j.aim.2020.107046}

\bibitem {GT}

A. Gholampour, R.P. Thomas, {\it Degeneracy loci, virtual cycles and nested Hilbert schemes},  Tunisian J. Math. {\bf 2} (2020), 633--665

\bibitem {GK1}

L. G\"ottsche, M. Kool, {\it Virtual refinements of the Vafa-Witten formula}, Commun. Math. Phys. {\bf 376} (2020), 1--49

\bibitem {GK2}

L. G\"ottsche, M. Kool, {\it Refined SU(3) Vafa-Witten invariants and modularity}, Pure Appl. Math. Quart. {\bf 14} (2018), 467--513

\bibitem {GNY} 

L. G\"ottsche, H. Nakajima, K. Yoshioka, {\it Instanton counting and Donaldson invariants}, J. Diff. Geom. {\bf 80} (2008), 343--390

\bibitem {GP} 

T. Graber, R. Pandharipande, 
{\it Localization of virtual classes}, Invent. Math. {\bf 135} (1999), 487--518

\bibitem {I} 

K. Intriligator, {\it Fusion residues}, Modern Physics Letters {\bf A6} (1991), 
3543--3556

\bibitem {JOP}

D. Johnson, D. Oprea, R. Pandharipande, {\it Rationality of descendent series for Hilbert and Quot schemes of surfaces}, Selecta Math., to appear, \texttt{arXiv:2002.08787}


\bibitem {KY}

T. Kawai, K. Yoshioka, {\it String partition functions and infinite products}, Adv. Theor.
Math. Phys. {\bf 4} (2000), 397--485


\bibitem {KL}

Y. H. Kiem, J. Li, {\it Localizing Virtual Cycles by Cosections},  J. Amer. Math. Soc. {\bf 26} (2013), 1025--1050

\bibitem {K}

M. Kool, {\it Stable pair invariants of surfaces and Seiberg-Witten invariants}, Quart. J. Math. Oxford {\bf 6} (2016), 365--386


\bibitem{KST}

M. Kool, V. Shende, R.P. Thomas, {\it A short proof of the G\"ottsche conjecture}, Geom. Topol. {\bf 15} (2011), 397--406


\bibitem {La}

T. Laarakker, {\it Monopole contributions to refined Vafa-Witten invariants,} \texttt{arXiv:1810.00385}

\bibitem{LeB} 

P. Le Barz, {\it Sur une formule de Castelnuovo pour les espaces multisecants},  Boll. Unione Mat. Ital. Sez. B {\bf 10} (2007), 381--387

\bibitem{LP} 

J. Lee, T. Parker, {\it A structure theorem for the Gromov-Witten invariants of K\"ahler surfaces},  J. Diff. Geom. {\bf 77} (2007),  483--513

\bibitem {L}

M. Lehn, {\it Chern classes of tautological sheaves on Hilbert schemes of points on surfaces},  Invent. Math. 
{\bf 136} (1999), 157--207


\bibitem{LM} 

M. Levine, F. Morel, {\it Algebraic cobordism}, Springer Monographs in Mathematics. Springer, Berlin, 2007


\bibitem{LevP} 

M. Levine, R. Pandharipande, {\it Algebraic cobordism revisited}, Invent. Math. {\bf 176} (2009),  63--130


\bibitem {Li}

J. Li, {\it Algebraic geometric interpretation of Donaldson's polynomial invariants of algebraic surfaces}, J. Diff. Geom. {\bf 37} (1993), 416--466


\bibitem {Lim}

W. Lim, {\it Virtual $\chi_{-y}$-genera of Quot schemes on surfaces}, preprint, \texttt{arXiv:2003.04429}

\bibitem {MarianCr} 

A. Marian,  {\it On the intersection theory of the Quot
schemes and moduli of bundles with sections}, J. Reine Angew. Math. {\bf 610} (2007), 13--27

\bibitem {MarOp} 

A. Marian, D. Oprea, {\it Counts of maps to Grassmannians
and intersections on the moduli space of bundles}, J. Diff. Geom. {\bf 76} (2007), 155--175


\bibitem {mo} 

A. Marian, D. Oprea, {\it Virtual intersections on the Quot scheme and Vafa-Intriligator formulas}, Duke Math. Journal {\bf 136} (2007), 81--131


\bibitem {quotients}

A. Marian, D. Oprea, R. Pandharipande, {\it The moduli space of stable quotients}, Geom. Topol. {\bf 15} (2011), 1651--1706

\bibitem {MOP2}

A. Marian, D. Oprea, R. Pandharipande, {\it Higher rank Segre integrals over the Hilbert scheme of points}, J. Eur. Math. Soc., to appear, \texttt{arXiv:1712.02382}

\bibitem {segre}

A. Marian, D. Oprea, R. Pandharipande, {\it Segre classes and Hilbert schemes of points}, Annales Scientifiques de l'ENS {\bf 50} (2017), 239--267

\bibitem {MOP}

A. Marian, D. Oprea, R. Pandharipande, {\it The combinatorics of Lehn's conjecture}, J. Math. Soc. Japan {\bf 71} (2019), 299--308


\bibitem{MNOP1} 

D. Maulik, N. Nekrasov, A. Okounkov, R. Pandharipande, {\it Gromov-Witten theory and Donaldson-Thomas theory I}, 
Compos. Math. {\bf 142} (2006), 126--1285

\bibitem{MNOP2} 

D. Maulik, N. Nekrasov, A. Okounkov, R. Pandharipande, {\it Gromov-Witten theory and Donaldson-Thomas theory II}, 
Compos. Math. {\bf 142} (2006), 1286--1304



\bibitem {MP} 
D. Maulik, R. Pandharipande, {\it New calculations in Gromov-Witten theory}, Pure Appl.
Math. Q. {\bf 4} (2008), 469--500

\bibitem {MPT} 
D. Maulik, R. Pandharipande, R. P. Thomas, {\it Curves on K3 surfaces and modular forms}, J. Topol. {\bf 3} (2010), 937--996



\bibitem {P}
R. Pandharipande,  {\it A calculus for the moduli space of curves}, in Algebraic geometry: Salt Lake City 2015,  
Proc. Sympos. Pure Math. {\bf 97}, Providence, RI, 2018, 459--487 

\bibitem {PP}
R. Pandharipande, A. Pixton, {\it Relations in the tautological ring of the moduli
space of curves}, \texttt{arXiv:1301.4561}


\bibitem {PP1}
R. Pandharipande, A. Pixton, {\it Descendents on local curves: Rationality},
Comp. Math. {\bf 149} (2013), 81--124

\bibitem {PP2}
R. Pandharipande, A. Pixton, {\it Descendent theory for stable
pairs on toric 3-folds}, J. Math. Soc. Japan {\bf 65} (2013), 1337--1372


\bibitem {PT22}

R. Pandharipande, R. P. Thomas, {\it Stable pairs and BPS invariants}, 
JAMS {\bf 23} (2010), 267--297


\bibitem {PT}

R. Pandharipande, R. P. Thomas, {\it The Katz-Klemm-Vafa conjecture for K3 surfaces}, Forum of Math {\bf 4} (2016), 1--111

\bibitem {Sc}

D. Schultheis, {\it Virtual invariants of Quot schemes over del Pezzo surfaces}, Ph.D. Thesis, UC San Diego (2012)

\bibitem{Sh} 

J. Shen, {\it Cobordism invariants of the moduli space of stable pairs}, J. Lond. Math. Soc. {\bf 94} (2016), 427--446

\bibitem {ST}

B. Siebert, G. Tian, {\it On quantum cohomology rings of Fano manifolds and
a formula of Vafa and Intriligator,} Asian J. Math. {\bf 1} (1997),  679--695

\bibitem {S}

R. Stanley, {\it Catalan numbers}, Cambridge University Press, 2015


\bibitem {TT1}

Y. Tanaka, R. P. Thomas, {\it Vafa-Witten invariants for projective surfaces I: stable case}, J. Alg. Geom. (to appear)

\bibitem {TT2}

Y. Tanaka, R. P. Thomas, {\it Vafa-Witten invariants for projective surfaces II: semistable case}, Pure Appl. Math. Quart. {\bf 13} (2017), 517--562

\bibitem {Th}

M. Thaddeus, {\it Stable pairs, linear systems and the Verlinde formula}, Invent. Math. {\bf 117} (1994), 317--353

\bibitem {T}

R. P. Thomas, {\it A holomorphic Casson invariant for Calabi-Yau 3-folds, and bundles on K3 fibrations}, J. Diff. Geom. {\bf 54} (2000), 367--438


\bibitem{VW} 

C. Vafa, E. Witten, {\it A strong coupling test of S-duality}, Nuclear Phys. B {\bf 431} (1994), 3--77


\bibitem {voisin}

C. Voisin, {\it Segre classes of tautological bundles on Hilbert schemes of surfaces}, Algebr. Geom. {\bf 6} (2019), 186--195


\bibitem {Wang} 

Z. Wang, {\it Tautological integrals on Hilbert schemes
of points on curves},  Acta Math. Sin. {\bf 32} (2016), 901--910

\bibitem {WZ}

Z. Wang, J. Zhou, {\it Generating series of intersection numbers on
Hilbert schemes of points}, Front. Math. China {\bf 12} (2017), 1247--1264


\bibitem {WW}

E. Whittaker, G. Watson, {\it A course of modern analysis}, Cambridge University Press, 1927

\end{thebibliography}
\end{document}